\documentclass[10pt,twoside,a4paper]{article}

\newcommand\sibrouillon[1]{}
\newcommand\Hum[1]{}

\usepackage[utf8]{inputenc} 
\usepackage[T1]{fontenc}
\usepackage{lmodern}
\usepackage{amsfonts}
\usepackage{mathrsfs}
\usepackage{mathtools}
\usepackage{stmaryrd}
\usepackage[all]{xy}
\usepackage[french]{babel}
\usepackage{amssymb}
\usepackage{amsthm}
\usepackage[pdftex]{graphicx}
\usepackage{natbib}
\usepackage{xspace}

\usepackage[bookmarksopen=false,breaklinks=true,%
      backref=page,pagebackref=true,plainpages=false,%
      hyperindex=true,pdfstartview=FitH,colorlinks=true,%
      pdfpagelabels=true,colorlinks=true,linkcolor=blue,%
      citecolor=red,urlcolor=red,hypertexnames=false%
      ]%
   {hyperref}

\usepackage{etoolbox}

\makeatletter\def\th@plain{\slshape}\makeatother
\makeatletter\patchcmd{\th@remark}{\itshape}{\slshape}{}{}\makeatother

\newcounter{bidon}
\newcommand{\rdb}{\refstepcounter{bidon}}

\newtheorem{theorem}{Théorème}[section]

\newtheorem{proposition}[theorem]{Proposition}
\newtheorem{lemma}[theorem]{Lemme}
\newtheorem{corollary}[theorem]{Corolaire}

\newtheorem{fact}[theorem]{Fait}

\theoremstyle{definition}
\newtheorem{definition}[theorem]{Définition}
\newtheorem{definitions}[theorem]{Définitions}

\newtheorem{notation}[theorem]{Notation}

\theoremstyle{remark}

\newcommand {\junk}[1]{}

\newcommand {\rem}{\noindent \textsl{Remarque. }  }
\newcommand {\rems}{\noindent \textsl{Remarques. }  }
\newcommand {\comm}{\noindent \textsl{Commentaire. }  }
\newcommand\exl{\noindent \textsl{Exemple. }}

\newcommand \eoe {\hbox{}\nobreak\hfill
\vrule width .5em height .5em depth 0mm \par \smallskip}

\newenvironment{Proof}[1]{
\trivlist \item[\hskip \labelsep{\sl #1.}]\hskip 0pt\\}{\hfill 
\mbox{$\Box$}
\endtrivlist}

\newcommand\hum[1] {\sibrouillon{{\begin{flushleft}\tt\small Hors-texte:  #1
\end{flushleft}}}}

\newcommand\oge{\leavevmode\raise.3ex\hbox{$\scriptscriptstyle\langle\!\langle\,$}}
\newcommand\feg{\leavevmode\raise.3ex\hbox{$\scriptscriptstyle\,\rangle\!\rangle$}}
\newcommand\gui[1]{\oge{#1}\feg}

\renewcommand \leq{\leqslant}
\renewcommand \preceq{\preccurlyeq}

\renewcommand \geq{\geqslant}

\newcommand \noi {\noindent}
\renewcommand \ss {\smallskip}
\newcommand \sni {\smallskip\noi}
\newcommand \ms {\medskip}
\newcommand \mni {\ms\noi}

\newcommand \Ex {{\exists }}
\newcommand \Tt {{\forall }}

\renewcommand \cir {{^\circ}}
\newcommand\equidef{\buildrel{{\rm 
def}}\over{\quad\Longleftrightarrow\quad}}
\newcommand\eqdefi{\buildrel{\rm def}\over {\;=\;}}

\newcommand\vers[1]{\buildrel{#1}\over \longrightarrow }

\newcommand\gen[1]{\left\langle{#1}\right\rangle}
\newcommand\so[1]{\{\,{#1}\,\}}

\newcommand \sotq[2]{\so{#1\,\vert\,#2}}

\newcommand\bidule[1]{\left\lceil{#1}\right\rceil}

\newcommand \snic[1] {\sni\centerline{$#1$}\ss}

\newcommand \aqo[2]{#1\!\left/\gen{#2}\right.}

\newcommand \ov[1] {\overline{#1}}
\newcommand \ovs[1] {\overline{\{#1\}}}

\newcommand \wi[1] {\widetilde{#1} }

\newcommand \dsp {\displaystyle}

\newcommand \CC{\mathbb {C}}
\newcommand \NN{\mathbb {N}}
\newcommand \ZZ{\mathbb {Z}}

\newcommand \gA{\mathbf{A}}
\newcommand \gB{\mathbf{B}}

\newcommand \gH{\mathbf{H}}
\newcommand \gK{\mathbf{K}}
\newcommand \gL{\mathbf{L}}

\newcommand \gT{\mathbf{T}}

\newcommand \cI {{\cal I}}

\newcommand \cF {{\cal F}}

\newcommand \cP {{\cal P}}

\newcommand \cM {{\cal M}}

\newcommand \rI {\mathrm{I}}
\newcommand \rD {\mathrm{D}}
\newcommand \rF {\mathrm{F}}
\newcommand \rH {\mathrm{H}}
\newcommand \rJ {\mathrm{J}}
\newcommand \rK {\mathrm{K}}

\newcommand \rS {\mathrm{S}}
\newcommand \rV {\mathrm{V}}
\newcommand \IZ {\mathrm{IZ}}
\newcommand \FZ {\mathrm{FZ}}
\newcommand \IZA {\IZ_\gA}
\newcommand \DA {\rD_\gA}
\newcommand \VA {\rV_{\!\gA}}

\newcommand \JA {\rJ_\gA}
\newcommand \JT {\rJ_\gT}
\newcommand \Fmin {\rF_{\mathrm{min}}}
\newcommand \Imax {\rI_{\mathrm{max}}}

\newcommand\fa{\mathfrak{a}}
\newcommand\fb{\mathfrak{b}}
\newcommand\fc{\mathfrak{c}}
\newcommand\fA{\mathfrak{A}}
\newcommand\fB{\mathfrak{B}}
\newcommand\fD{\mathfrak{D}}
\newcommand \DT {\fD_\gT}
\newcommand \VT {\fV_\gT}
\newcommand \DTo {\fD_{\gT\cir}}
\newcommand \VTo {\fV_{\gT\cir}}

\newcommand\fII{\mathfrak{I}}
\newcommand\fj{\mathfrak{j}}
\newcommand\fJ{\mathfrak{J}}
\newcommand\fF{\mathfrak{F}}
\newcommand\ff{\mathfrak{f}}
\newcommand\ffg{\mathfrak{g}}

\newcommand\fm{\mathfrak{m}}

\newcommand\fp{\mathfrak{p}}
\newcommand\fP{\mathfrak{P}}
\newcommand\fq{\mathfrak{q}}
\newcommand\fV{\mathfrak{V}}

\newcommand \vu {\vee} 
\newcommand \vi {\wedge} 
\newcommand \Vu {\bigvee}
\newcommand \Vi {\bigwedge}
\newcommand \im {\rightarrow} 
\newcommand \dar { \,\downarrow\! }
\newcommand \uar { \,\uparrow\! }
\newcommand \bal[1] {^\rK_{#1}}
\newcommand \ul[1] {_\rK^{#1}}

\newcommand \Un {{\bf 1}}
\newcommand \Deux {{\bf 2}}
\newcommand \Trois {{\bf 3}}

\newcommand \Pf {{{\cal P}_{\mathrm{f}}}}

\newcommand \Bd {\mathrm{Bd}}
\newcommand \Frac {\,\mathrm{Frac}}
\newcommand \He {\mathrm{He}}
\newcommand \Heit {\mathrm{Heit}}
\newcommand \HeA {{\Heit\,\gA}}

\newcommand \Jspec {\mathsf{Jspec}}
\newcommand \jspec {\mathsf{jspec}}

\newcommand \Hdim {\mathsf{Hdim}}
\newcommand \Kdim {\mathsf{Kdim}}
\newcommand \Jdim {\mathsf{Jdim}}
\newcommand \jdim {\mathsf{jdim}}

\newcommand \Id {\mathrm{Id}}

\newcommand \Ann {\mathrm{Ann}}

\newcommand \Max {\mathsf{Max}}
\newcommand \Min {\mathsf{Min}}
\newcommand \Patch {\mathsf{Patch}}
\newcommand \Spec {\mathsf{Spec}}
\newcommand \Zar {\,\mathsf{Zar}}
\newcommand \ZarA {{\Zar\,\gA}}

\newcommand \OQC {\mathsf{Oqc}}
\renewcommand \mod {\,\mathrm{mod}}

\newcommand \recu {récur\-rence }
\newcommand \recuz {récur\-rence}
\newcommand \hdr {hypo\-thèse de \recu }

\newcommand \cad {c'est-à-dire }
\newcommand \Cad {C'est-à-dire }
\newcommand \cade {c'est-à-dire encore }
\newcommand \ssi {si, et seu\-le\-ment si, }

\newcommand \spdg {sans perte de géné\-ra\-lité }
\newcommand \Propeq {Les pro\-priétés sui\-van\-tes sont
équiva\-lentes:}
\newcommand \propeq {les pro\-priétés sui\-van\-tes sont
équiva\-lentes:}

\newcommand \Amo {$\gA$-mo\-du\-le }
\newcommand \Amos {$\gA$-mo\-du\-les }

\newcommand \Amosz {$\gA$-mo\-du\-les}

\newcommand \agB {\alg de Boole }

\newcommand \agBz {\alg de Boole}

\newcommand \agH {\alg de Heyting }
\newcommand \agHs {\algs de Heyting }
\newcommand \agHz {\alg de Heyting}
\newcommand \agHsz {\algs de Heyting}

\newcommand \alg {al\-gè\-bre }
\newcommand \algs {al\-gè\-bres }

\newcommand \auto {auto\-mor\-phisme }
\newcommand \autos {auto\-mor\-phismes }

\newcommand \com {co\-ma\-xi\-maux }

\newcommand \dfn{défi\-nition }
\newcommand \dfns{défi\-nitions }
\newcommand \dfnz{défi\-nition}
\newcommand \dfnsz{défi\-nitions}

\newcommand \egmt {éga\-lement }

\newcommand \egt {éga\-lité }
\newcommand \egts {éga\-lités }
\newcommand \egtz {éga\-lité}
\newcommand \egtsz {éga\-lités}

\newcommand \elr{élé\-men\-tai\-re }
\newcommand \elrs{élé\-men\-tai\-res }
\newcommand \elrz{élé\-men\-tai\-re}
\newcommand \elrsz{élé\-men\-tai\-res}

\newcommand \elt{élé\-ment }
\newcommand \elts{élé\-ments }

\newcommand \eltsz{élé\-ments}

\newcommand \eqvc {équivalence }

\newcommand \gtr{gé\-né\-ra\-teur }
\newcommand \gtrs{gé\-né\-ra\-teurs }
\newcommand \gtrz{gé\-né\-ra\-teur}
\newcommand \gtrsz{gé\-né\-ra\-teurs}

\newcommand \homo {ho\-mo\-mor\-phisme }

\newcommand \homos {ho\-mo\-mor\-phismes }

\newcommand \id {idéal }
\newcommand \ids {idéaux }
\newcommand \idz {idéal}
\newcommand \idsz {idéaux}

\newcommand \idema {\id ma\-xi\-mal }
\newcommand \idemaz {\id ma\-xi\-mal}
\newcommand \idemas {\ids ma\-xi\-maux }
\newcommand \idemasz {\ids ma\-xi\-maux}

\newcommand \idf {idéal de Fitting }

\newcommand \idep {\id pre\-mier }
\newcommand \idepz {\id pre\-mier}
\newcommand \ideps {\ids pre\-miers }
\newcommand \idepsz {\ids pre\-miers}

\newcommand \iso {iso\-mor\-phisme }

\newcommand \itf {\id \tf}
\newcommand \itfs {\ids \tf}
\newcommand \itfz {\id \tfz}
\newcommand \itfsz {\ids \tfz}

\newcommand \mo {mo\-no\"{\i}de }
\newcommand \moco {\mos\com}

\newcommand \mos {mo\-no\"{\i}des }

\newcommand \oqc {ouvert \qc}
\newcommand \oqcs {ouverts \qcs}
\newcommand \oqcz {ouvert \qcz}
\newcommand \oqcsz {ouverts \qcsz}

\newcommand \pf {de pré\-sen\-ta\-tion finie }
\newcommand \pfz {de pré\-sen\-ta\-tion finie}

\newcommand \qc {quasi-compact }
\newcommand \qcs {quasi-compacts }
\newcommand \qcz {quasi-compact}
\newcommand \qcsz {quasi-compacts}

\newcommand \tf {de type fini }
\newcommand \tfz {de type fini}
\newcommand \trdi {treillis distri\-butif }
\newcommand \trdis {treillis distri\-butifs }
\newcommand \trdiz {treillis distri\-butif}
\newcommand \trdisz {treillis distri\-butifs}

\newcommand \Tho {Théo\-rè\-me }
\newcommand \tho {théo\-rè\-me }
\newcommand \thos {théo\-rè\-mes }
\newcommand \thoz {théo\-rè\-me}

\newcommand \cof {cons\-truc\-tif }

\newcommand \cofz {cons\-truc\-tif}
\newcommand \cofsz {cons\-truc\-tifs}

\newcommand \cov {cons\-truc\-ti\-ve }
\newcommand \covz {cons\-truc\-ti\-ve}
\newcommand \covsz {cons\-truc\-ti\-ves}
\newcommand \covs {cons\-truc\-ti\-ves }

\newcommand \coma {\maths\covs}
\newcommand \comaz {\maths\covsz}
\newcommand \clama {\maths classiques }
\newcommand \clamaz {\maths classiques}

\renewcommand \cot {cons\-truc\-ti\-ve\-ment }
\newcommand \cotz {cons\-truc\-ti\-ve\-ment}

\newcommand \maths {ma\-thé\-ma\-ti\-ques }

\newcommand \prco {preuve \cov}

\newdimen\xyrowsp
\newcommand{\SCO}[6]{
\xymatrix @R = \xyrowsp {
                                  &1 \ar@{-}[dl] \ar@{-}[dr] \\
#3 \ar@{-}[ddr]                   &   & #6 \ar@{-}[ddl] \\
                                  &\bullet\ar@{-}[d] \\
                                  &\bullet   \\
#2 \ar@{-}[ddr] \ar@{-}[uur]      &   & #5 \ar@{-}[ddl] \ar@{-}[uul] \\
                                  &\bullet \ar@{-}[d] \\
                                  &\bullet  \\
#1 \ar@{-}[uur]                   &   & #4 \ar@{-}[uul] \\
                                  & 0 \ar@{-}[ul] \ar@{-}[ur] \\
}
}

\newcounter{MF}
\newcommand\stMF{\stepcounter{MF}}

\newcommand{\lec}{\stMF\ifodd\value{MF}lecteur \else 
lectrice \fi}
\newcommand{\lecz}{\stMF\ifodd\value{MF}lecteur\else lectrice\fi}

\newcommand{\lecs}{\stMF\ifodd\value{MF}lecteurs \else 
lectrices \fi}
\newcommand{\lecsz}{\stMF\ifodd\value{MF}lecteurs\else 
lectrices\fi}

\newcommand{\alec}{\stMF\ifodd\value{MF}au lecteur \else%
à la lectrice \fi}
\newcommand{\alecz}{\stMF\ifodd\value{MF}au lecteur\else%
à la lectrice\fi}

\newcommand{\dlec}{\stMF\ifodd\value{MF}du lecteur \else%
de la lectrice \fi}
\newcommand{\dlecz}{\stMF\ifodd\value{MF}du lecteur\else%
de la lectrice\fi}

\newcommand{\llec}{\stMF\ifodd\value{MF}le lecteur \else la lectrice \fi}
\newcommand{\llecz}{\stMF\ifodd\value{MF}le lecteur\else la lectrice\fi}

\newcommand{\Llec}{\stMF\ifodd\value{MF}Le lecteur \else La lectrice \fi}

\newcommand{\lui}{\ifodd\value{MF}lui \else
elle \fi}
\newcommand{\luiz}{\ifodd\value{MF}lui\else
elle\fi}

\newcommand{\celui}{\ifodd\value{MF}celui \else
celle \fi}

\newcommand{\ceux}{\ifodd\value{MF}ceux \else
celles \fi}

\newcommand{\er}{\ifodd\value{MF}er \else
ère \fi}

\newcommand{\eux}{\ifodd\value{MF}eux \else
elles \fi}

\newcommand{\eUx}{\ifodd\value{MF}eux \else
euse \fi}

\newcommand{\leux}{\ifodd\value{MF}leux \else
leuse \fi}

\newcommand{\il}{\ifodd\value{MF}il \else
elle \fi}

\newcommand{\ien}{\ifodd\value{MF}ien \else
ienne \fi}

\newcommand{\e}{\ifodd\value{MF} \else e \fi}
\newcommand{\ez}{\ifodd\value{MF}\else e\fi}

\newcommand{\n}{\ifodd\value{MF}n \else nne \fi}
\newcommand{\nz}{\ifodd\value{MF}n\else nne\fi}

\makeatletter
\newcommand{\la}{\@ifstar{\ifodd\value{MF}le\else
la\fi}{\stMF\ifodd\value{MF}le \else la \fi}}
\makeatother

\RequirePackage{etoolbox}
\patchcmd{\sectionmark}{\MakeUppercase}{}{}{}
\begin{document}
\title{Dimension de Heitmann des \trdis et des anneaux commutatifs}

\author{
Thierry Coquand
(\thanks {~
Chalmers, University of G\"oteborg, Sweden,
email: coquand@cs.chalmers.se}~),
Henri Lombardi
(\thanks{~
Laboratoire de Mathématiques de Besançon, CNRS UMR 6623, 
Université Bourgogne Franche-Comté, 25030 BESANCON cedex, FRANCE,
email: henri.lombardi@univ-fcomte.fr, \url{http://hlombardi.free.fr}.
}~),
Claude Quitté
(\thanks {~
Laboratoire de Mathématiques,
SP2MI, Boulevard 3, Teleport 2, BP 179,
86960 FUTUROSCOPE Cedex, FRANCE,
email: quitte@math.univ-poitiers.fr}~)
}
\date{
Version corrig\'ee: \today.\\[.3em]
Article original: 
\textsl{Publications math\'ematiques de Besançon. \\[.2em]
Alg\`ebre et Th\'eorie des Nombres.}
(2006), pages 57--100.
}
\maketitle


\begin{abstract}
Nous \'etudions la notion de dimension introduite par Heitmann dans
son article remarquable de 1984 \cite{Hei84}, ainsi qu'une notion
voisine, seulement implicite dans ses preuves.  Nous d\'eveloppons
ceci d'abord dans le cadre g\'en\'eral de la th\'eorie des treillis distributifs et
des espaces spectraux.  Nous obtenons ensuite des versions constructives de
certains th\'eor\`emes importants d'alg\`ebre commutative.  Les versions constructives de
 ces th\'eor\`emes s'av\`erent en fin de compte plus simples, et parfois
plus g\'en\'erales, que les versions classiques abstraites
correspondantes.
\end{abstract}
\mni MSC 2000: 13C15, 03F65, 13A15, 13E05

\sni Mots clés :  Dimension de Krull, Dimension de Heitmann, Bord 
d'une
sous-variété, \Tho \gui{Stable range} de Bass, \Tho 
\gui{Splitting off} de
Serre, \Tho de Forster-Swan, \Tho de Kronecker, Nombre de \gtrs d'un 
module,
Mathématiques constructives.

\sni Key words: Krull dimension, Boundary of a subvariety,
Constructive Mathematics, Kronecker theorem, Heitmann dimension, Forster-Swan Theorem, Serre's Splitting Off, Bass cancellation theorem.

\begin{center}
{\bf \large Avertissement} 
\end{center}

L'article original est paru aux \textsl{Publications mathématiques de Besançon. Algèbre et Théorie des Nombres.}
(2006), pages 57--100.

Nous corrigeons ici un certain nombre d'erreurs, répertoriées en post-scriptum à la fin du texte.
Nous donnons aussi des références bibliographiques supplémentaires.

Signalons enfin que les sections \ref{secKroBass} à \ref{secSwan} ont fait, bien après la parution de l'article original, l'objet d'un exposé détaillé dans l'ouvrage
\cite[Alg\`ebre Commutative. Méthodes constructives. 
 Calvage \& Mounet]{ACMC}.
 
Nous nous conformons à l'orthographe nouvelle recommandée (par exemple: à priori, corolaire, connaitre), et les \lecs subissent l'alternance des sexes.

\newpage
\rdb \thispagestyle{empty}
{\small
\tableofcontents
}

\markboth{Introduction}{ }
\section*{Introduction} \label{sec Introduction}
\addcontentsline{toc}{section}{Introduction}

Nous étudions la notion de dimension introduite par Heitmann dans
son article de 1984 \cite{Hei84}, ainsi qu'une notion voisine,
seulement implicite dans ses preuves.  Nous développons ceci d'abord
dans le cadre général de la théorie des \trdis et des espaces
spectraux.  Nous appliquons ensuite cette problématique dans le
cadre de l'\alg commutative.

Dans la dualité entre \trdis et espaces spectraux, le spectre de
Zariski d'un anneau commutatif correspond (comme l'a indiqué André Joyal dans
\cite{Joy76}) au treillis des \ids qui sont radicaux d'\itfsz.  Nous
montrons que l'espace spectral défini par Heitmann pour sa notion de
dimension correspond au treillis formé par les idéaux qui sont
radicaux de Jacobson d'\itfsz.  Ceci nous permet d'obtenir une \dfn
\cov \elr de la dimension définie par Heitmann (que nous notons
$\Jdim$).  Nous introduisons une autre dimension, que nous appelons
dimension de Heitmann (et que nous notons $\Hdim$), qui est
\gui{meilleure} en ce sens que $\Hdim\leq \Jdim$ et qu'elle permet des
preuves par \recu naturelles.

Nous obtenons alors des versions \covs de certains \thos classiques
importants, dans leur version non noethérienne (en général due
à Heitmann).

Les versions \covs de ce ces \thos s'avèrent en fin de compte plus
simples, et parfois plus générales, que les versions classiques
abstraites correspondantes.

En particulier nous rappelons les versions non noethériennes des
\thos de Swan et de Serre (splitting off) obtenues récemment pour la
première fois dans \cite{CLQ2004} et \cite{Duc2006}.

\smallskip
Naturellement, le principal avantage que nous voyons dans notre
traitement est son caractère tout à fait \elrz.  En particulier
nous n'utilisons pas d'hypothèses \gui{non nécessaires} comme
l'axiome du choix et le principe du tiers exclu, inévitables pour
faire fonctionner les preuves classiques antérieures.

Enfin, le fait de s'être débarrassé de toute hypothèse
noethérienne est aussi non négligeable, et permet de mieux voir
l'essence des choses.

\medskip En conclusion cet article peut être vu pour l'essentiel
comme une mise au point \cov de la théorie des espaces spectraux via
celle des \trdisz, avec une insistance particulière sur la dimension
de Heitmann et quelques applications marquantes en algèbre
commutative.

\medskip \rem Nous avons résolu de la manière suivante un
problème de terminologie qui se pose en rédigeant cet article.  Le
mot \gui{dualité} apparaît à priori dans le contexte des \trdis
avec deux significations différentes.  Il y a d'une part la
dualité qui correspond au renversement de la relation d'ordre dans
un treillis.  D'autre part il y a la dualité entre treillis
distributifs et espaces spectraux, qui correspond à une
antiéquivalence de catégorie.  Nous avons décidé de réserver
\gui{dualité} pour ce dernier usage.  Le terme \gui{treillis dual} a
donc été systématiquement remplacé par \gui{treillis
opposé}.  De même on a remplacé \gui{la notion duale} par
\gui{la notion renversée} ou par \gui{la notion opposée}, et
\gui{par dualité} par \gui{par renversement de l'ordre}.

\newpage
\section{Treillis distributifs}
\label{secTRDI}

Les axiomes des \trdis peuvent être formulés avec des \egts
universelles concernant uniquement les deux lois $\vi$ et $\vu$ et les
deux constantes $0_\gT$ (l'\elt minimum du \trdi $\gT$) et $1_\gT$ (le
maximum).  La relation d'ordre est alors définie par $a\leq_\gT
b\;\Leftrightarrow\;a\vi b=a$.  On obtient ainsi une théorie
purement équationnelle, avec toutes les facilités afférentes. 
Par exemple on peut définir un \trdi par générateurs et
relations, la catégorie comporte des limites inductives (qu'on peut
définir par générateurs et relations) et des limites projectives
(qui ont pour ensembles sous-jacents les limites projectives
ensemblistes correspondantes).

Un ensemble totalement ordonné est un \trdi s'il possède un
maximum et un minimum.  On note ${\bf n}$ un ensemble totalement
ordonné à $n$ éléments, c'est un \trdi si $n\neq 0$.  Le
treillis $\Deux$ est le \trdi libre à 0 \gtrz, et $\Trois$ celui à
un \gtrz.

Pour tout \trdi $\gT$, si l'on remplace la relation d'ordre $x\leq_\gT
y$ par la relation symétrique $y\leq_\gT x$ on obtient le
\textsl{treillis opposé} $\gT\cir$ avec échange de $\vi$ et $\vu$
(on dit parfois \textsl{treillis dual}).

\subsection{Idéaux, filtres}

Si $\varphi :\gT\rightarrow \gT'$ est un morphisme de \trdisz, 
$\varphi^{-1}(0)$
est appelé un \textsl{\id de $T$}. Un \id $\fII $ de $\gT$ est une 
partie de
$\gT$ soumise
aux  contraintes suivantes:
\begin{equation} \label{eqIdeal}
\left.
\begin{array}{rcl}
   & &  0 \in \fII    \\
x,y\in \fII & \Longrightarrow   &  x\vu y \in \fII    \\
x\in \fII ,\; z\in \gT& \Longrightarrow   &  x\vi z \in \fII    \\
\end{array}
\right\}
\end{equation}
(la dernière se réécrit $(x\in \fII ,\;y\leq x)\Rightarrow y\in 
\fII $).
Un \textsl{\id principal} est un \id engendré par un seul \elt $a$:
il est égal à
\begin{equation} \label{eqda}
\,\dar a=\sotq{x\in \gT}{x\leq a}
\end{equation}
L'\id $\,\dar a$, muni des lois $\vi$ et $\vu$ de $\gT$ est un \trdi
dans lequel l'\elt maximum est~$a$.  L'injection canonique $\,\dar
a\rightarrow \gT$ \textsl{n'est pas} un morphisme de \trdis parce que
l'image de $a$ n'est pas égale à $1_\gT$.  Par contre
l'application surjective $\gT\rightarrow \,\dar a,\;x\mapsto x\vi a$
est un morphisme surjectif, qui définit donc $\dar a$ comme une
structure quotient.

\smallskip La notion opposée à celle d'\id est la notion de {\em
filtre}.  Le filtre principal engendré par $a$ est noté $\,\uar a$.

\smallskip L'\textsl{\id engendré} par une partie $J$ de $\gT$ est
$\cI_\gT(J)=\sotq{x\in\gT}{\Ex J_0\in \Pf(J),\,x\leq \Vu J_0}$. 
En conséquence \textsl{tout \itf est principal}.

Si $A$ et $B$ sont deux parties de $\gT$ on note
\begin{equation} \label{eqvuvi}
A\vu B=\sotq{a\vu b}{a\in A,\,b\in B}  \quad \mathrm{et}\quad  A\vi
B=\sotq{a\vi b}{a\in A,\,b\in B}.
\end{equation}

Alors l'\id engendré par deux \ids $\fa$ et $\fb$ est égal à
\begin{equation} \label{eqSupId}
\cI_\gT(\fa\cup \fb) = \fa\vu\fb =\sotq{z}{\exists 
x\in\fa,\,\exists
y\in\fb,\,z\leq x\vu y}\,.
\end{equation}

L'ensemble des \ids de $\gT$ forme lui même un \trdi pour
l'inclusion, avec pour inf de $\fa$ et $\fb$ l'\idz:
\begin{equation} \label{eqInfId}
\fa\cap \fb=\fa\vi\fb.
\end{equation}

Ainsi les opérations $\vu$ et $\vi$ définies en (\ref{eqvuvi})
correspondent au sup et au inf dans le treillis des \idsz.

On notera $\cF_\gT(S)$ le filtre de $\gT$ engendré par le sous
ensemble $S$.  
Quand on considère le treillis des filtres il faut
faire attention à ce que produit le renversement de la relation
d'ordre: $\ff\cap\ffg=\ff\vu\ffg$ est le inf de $\ff$ et $\ffg$,
tandis que leur sup est égal à $\cF_\gT(\ff\cup \ffg)=\ff\vi
\ffg$.

\medskip Le \textsl{treillis quotient de $\gT$ par l'\id $\fJ$}, noté
$\gT/(\fJ=0)$ est défini comme le \trdi engendré par les \elts de
$\gT$ avec pour relations, les relations vraies dans $\gT$ d'une part,
et les relations $x=0$ pour les $x\in \fJ$ d'autre part.  Il peut
aussi être défini par la relation de préordre
$$ a\preceq b\quad\Longleftrightarrow\quad a\leq_{\gT/(\fJ=0)}b 
\equidef \quad
\exists x\in \fJ \;\;a\leq  x\vu b
$$
Ceci donne
$$ a\equiv b\;\;\mod\;(\fJ=0)\quad \Longleftrightarrow  \quad \exists 
x\in \fJ
\;\;a\vu x=b\vu x
$$
et dans le cas du quotient par un \id principal $\,\dar a$ on obtient
$\gT/(a=0)\simeq\,\uar a$ avec le morphisme $y\mapsto y\vu a$ de $\gT$ 
vers
$\,\uar a$.

\subsubsection*{Transporteur, différence}
\addcontentsline{toc}{subsubsection}{Transporteur, différence}

Par analogie avec l'\alg commutative, si $\fb$  est un \id et  $A$ 
une partie de
$\gT$ on  notera
\begin{equation} \label{eqTrans}
\fb:A\eqdefi\sotq{x\in\gT}{\Tt a\in A\quad a\vi x\in  \fb}
\end{equation}
Si $\fa$ est l'\id engendré par $A$ on a $\fb:A=\fb:\fa$,
on l'appelle le \textsl{transporteur de $\fa$ dans $\fb$.}

On note aussi
$b:a$ l'\id
$ (\dar b):(\dar a)=\sotq{x\in\gT}{x\in\gT\;|\;x\vi a\leq  b}$.

La notion opposée est celle de \textsl{filtre différence de deux 
filtres}
\begin{equation} \label{eqDiff}
\ff\setminus\ff'\eqdefi\sotq{x\in\gT}{\Tt a\in \ff'\quad a\vu 
x\in  \ff}
\end{equation}
On note aussi
$b\setminus a$ le filtre
$ (\uar b)\setminus(\uar a)=\sotq{x\in\gT}{b\leq  x\vu a}$.

\subsubsection*{Radical de Jacobson}
\addcontentsline{toc}{subsubsection}{Radical de Jacobson}

Un \id $\fm$ d'un \trdi $\gT$ non trivial (i.e. distinct de $\Un$) est
dit \textsl{maximal} \hbox{si $\gT/(\fm=0)\,=\,\Deux$}, \cad si $1\notin\fm$ et
$\Tt x\in\gT\;(x\in\fm$ ou $\Ex y\in \fm \;x\vu y=1)$.

Il revient au même de dire qu'il s'agit d'un \id \gui{maximal
parmi les \ids stricts}.

En \clama on a le lemme suivant.
\begin{lemma}
\label{lemHspec1}
Dans un \trdi $\gT\neq\Un$ l'intersection des \idemas est égale à 
l'\id
$$ \sotq{a\in\gT}{\forall x\in\gT \;(a\vu x = 1 \Rightarrow  
x=1)}.
$$
On l'appelle le \textsl{radical de Jacobson de $\gT$}. On le note 
$\JT(0)$. \\
Plus généralement l'intersection des \idemas contenant un \id 
strict $\fJ$
est égale à l'\id
\begin{equation} \label{eqRJJ}
\JT(\fJ)\,=\, \sotq{a\in\gT}{\forall x\in\gT \;(a\vu x  = 1 
\Rightarrow \Ex
z\in \fJ \;\;z\vu x=1)}
\end{equation}
On l'appelle le \textsl{radical de Jacobson de l'\id $\fJ$}. En 
particulier:
\begin{equation} \label{eqRJb}
\JT(\dar b)=\sotq{a\in\gT}{\forall x\in\gT \;(\,a\vu x = 1\; 
\Rightarrow
\;\;b\vu x=1\,)}
\end{equation}
\end{lemma}
\begin{proof}
La deuxième affirmation résulte de la première en passant au
treillis quotient $\gT/(\fJ=0)$.  Voyons la première.  On montre que
$a$ est en dehors d'au moins un \idema \ssi $\Ex x\neq 1$ tel que
$a\vu x=1$.  Si c'est le cas, un \idema qui contient $x$ (il en existe
puisque $x\neq 1$) ne peut pas contenir $a$ car il contiendrait $a\vu
x$.  Inversement, si $\fm$ est un \idema ne contenant pas $a$, l'\id
engendré par $\fm$ et $a$ contient $1$.  Or cet \id est l'ensemble
des \elts majorés par au moins un $a\vu x$ où $x$ parcourt $\fm$.
\end{proof}

En \clama un \trdi est appelé \textsl{treillis de Jacobson} si tout
\idep est égal à son radical de Jacobson.  Comme tout \id est
intersection des \ideps qui le contiennent, cela implique que tout \id
est égal à son radical de Jacobson.

En \coma on  adopte les \textsl{\dfnsz} suivantes.

\medskip 
\begin{definitions}
\label{defJac} \textsl{(radical de Jacobson, treillis faiblement 
Jacobson)}
\begin{enumerate}
\item Si $\fJ$ est un \id de $\gT$  son \textsl{radical de Jacobson} 
est défini
par l'\egt $(\ref{eqRJJ})$
(on ne fait pas l'hypothèse que  $\gT\neq\Un$).
On notera $\JT(a)$ pour  $\JT(\dar a)$.
\item Un \trdi est appelé un \textsl{treillis faiblement Jacobson} si 
tout \id
principal est égal à son radical de Jacobson, \cade si
\begin{equation} \label{eqTJac}
\Tt a,b\in\gT\;\;[\,(\forall x\in\gT\, (a\vu 
x=1 \;\Rightarrow\;b\vu 
x=1))\;\Rightarrow
\;a\leq b\,]
\end{equation}
\end{enumerate}
\end{definitions}

On vérifie sans difficulté que $\JT(\fJ)$ est un \id et que
$1\in\JT(\fJ)\Leftrightarrow 1\in\fJ$.

\hum{Cela serait bien d'avoir une formulation \cov
  sans quantification sur \gui{l'ensemble} des \ids
  pour la notion de treillis
  de Jacobson. Encore que la notion la plus pertiente
  semble être celle de treillis faiblement Jacobson
  puisque le $\jspec$ doit être remplacé par le $\Jspec$}

\subsection{Quotients}

Un \textsl{\trdi quotient $\gT'$ de $\gT$} est donné par une relation 
binaire
$\preceq$ sur $\gT$ vérifiant les propriétés suivantes:
\begin{equation} \label{eqPreceq}
\left.
\begin{array}{rcl}
a\leq b&  \Longrightarrow  & a\preceq b   \\
a\preceq b,\,b\preceq c&  \Longrightarrow  & a\preceq c   \\
a\preceq b,\,a\preceq c&  \Longrightarrow  & a\preceq b\vi c   \\
b\preceq a,\,c\preceq a&  \Longrightarrow  & b\vu c\preceq a
\end{array}
\right\}
\end{equation}

\begin{proposition}
\label{propIdealFiltre} Soit $\gT$ un \trdi et
$(J,U)$ un couple de parties de $\gT$.
On considère le quotient $\gT'$ de $\gT$ défini par les
relations $x=0$ pour les $x\in J$ et $y=1$ pour les $y\in U$. Alors
  on a $a\leq_{\gT'}b$ \ssi
il existe une partie finie $J_0$ de $J$ et une partie finie $U_0$ de
$U$ telles que:
\begin{equation} \label{eqpropIdealFiltre}
a \vi \Vi U_0 \; \leq_\gT\; b \vu \Vu J_0
\end{equation}
Nous noterons $\gT/(J=0,U=1)$ ce treillis quotient $\gT'.$
\end{proposition}

\subsubsection*{Idéaux dans un quotient}
\addcontentsline{toc}{subsubsection}{Idéaux dans un quotient}

Le fait suivant résulte des \egts (\ref{eqIdeal}), (\ref{eqvuvi}),
(\ref{eqSupId}) et (\ref{eqInfId}).
\begin{fact}
\label{factIdDansQuo}
Soit  $\pi:\gT\to\gL$ un treillis quotient.
\begin{itemize}
\item L'image réciproque d'un \id de $\gL$ par $\pi^{-1}$ est un 
\id de $\gT$,
ceci donne un morphisme pour $\vu$ et $\vi$ (mais pas nécessairement
pour $0$ et $1$).
\item L'image d'un \id de $\gT$ par $\pi$ est un \id de $\gL$,
ceci donne un \homo surjectif de treillis.
\item Un \id de $\gT$ est de la forme $\pi^{-1}(\fa)$ \ssi il est 
saturé pour
la relation~$=_\gL$.
\item Résultats analogues pour les filtres.
\end{itemize}
\end{fact}

Notez qu'en algèbre commutative, le morphisme de passage au 
quotient par un
idéal ne se comporte pas aussi bien pour les \ids dans le cas d'une
intersection puisqu'on peut très bien avoir 
$\fa+(\fb\cap\fc)\varsubsetneq
(\fa+\fb)\cap(\fa+\fc)$.

  Le lemme suivant donne quelques renseignements complémentaires 
pour les
quotients par un \id et par un filtre.

\begin{lemma}
\label{lemIQT}
Soit $\fa$ un \id et $\ff$ un filtre de $\gT$
\begin{enumerate}
\item Si $\gL=\gT/(\fa=0)$ alors la projection canonique 
$\pi:\gT\to\gL$
établit une bijection croissante entre les \ids  de $\gT$ contenant 
$\fa$ et
les \ids de $\gL$. La bijection réciproque est fournie par 
$\fj\mapsto\pi^{-
1}(\fj)$. En outre si $\fj$ est un \id de $\gL$, on obtient~\hbox{$\pi^{-1}(\rJ_\gL(\fj)) = \JT(\pi^{-1}(\fj))$}.
\item  Si $\gL=\gT/(\ff=1)$ alors la projection canonique 
$\pi:\gT\to\gL$
établit une bijection croissante entre les \ids $\fJ$ de $\gT$ 
vérifiant
\gui{$\Tt f \in \ff,\;\;\fJ:f=\fJ$} et les \ids de~$\gL$.
\end{enumerate}
\end{lemma}

\rem On notera que $\gT\mapsto\JT(0)$ n'est pas une opération 
fonctorielle. La
deuxième affirmation du point \textsl{1} du lemme précédent, qui admet 
une preuve
\cov directe,  s'explique facilement en \clama par le fait que, dans 
le cas
très particulier du quotient par un \idz, les \idemas de $\gL$ 
contenant $\fj$
correspondent par $\pi^{-1}$ aux \idemas de~$\gT$
contenant~$\pi^{-1}(\fJ)$.

\subsubsection*{Recollement de treillis quotients}
\addcontentsline{toc}{subsubsection}{Recollement de treillis quotients}

En \alg commutative, si $\fa$ et $\fb$ sont deux \ids d'un anneau 
$\gA$
on a une \gui{suite exacte} de \Amos (avec $j$ et $p$ des \homos 
d'anneaux)
$$0\to\gA/(\fa\cap\fb)\vers{j}(\gA/\fa) \times 
(\gA/\fb)\vers{p}\gA/(\fa+\fb)\to 0$$
qu'on peut lire en langage courant: le système de congruences  
$x\equiv
a\;\mod\;\fa$, $x\equiv b\;\mod\;\fb$ admet une solution \ssi $a\equiv
b\;\mod\;\fa+\fb$ et dans ce cas la solution est unique modulo 
$\fa\cap\fb$.
Il est remarquable que ce \gui{\tho des restes chinois} se 
généralise à un
système \textsl{quelconque} de congruences \ssi l'anneau est
\textsl{arithmétique} \cite[Théorème XII-1.6]{ACMC}, \cad si le treillis des \ids est distributif.
Le \tho des restes chinois \gui{contemporain} concerne le cas 
particulier d'une
famille d'\ids deux à deux comaximaux, et il fonctionne sans 
hypothèse sur
l'anneau de base.

D'autres épimorphismes de la catégorie des anneaux commutatifs 
sont les
localisations. Et il y a un principe de recollement analogue au \tho 
des restes
chinois pour les localisations, extrêmement fécond (le principe 
local-global).

\smallskip De la même manière on peut récupérer un \trdi 
à partir
d'un nombre fini de ses quotients,
si l'information qu'ils contiennent est \gui{suffisante}. On peut 
voir ceci au
choix comme une procédure de recollement (de passage du local au 
global), ou
comme une version du \tho des restes chinois pour les \trdisz. Voyons 
les choses
plus précisément.

\begin{definition}
\label{defRecolTD}
Soit $\gT$ un \trdiz, $(\fa_i)_{i=1,\ldots n}$ (resp. 
$(\ff_i)_{i=1,\ldots n}$)
une famille finie d'\ids (resp. de filtres)  de $\gT$.  On dit que 
les \ids
$\fa_i$ \textsl{recouvrent $\gT$} si $\bigcap_i\fa_i=\so{0}$. De 
même on dit
que les filtres $\ff_i$ \textsl{recouvrent $\gT$} si 
$\bigcap_i\ff_i=\so{1}$.
\end{definition}

Pour un \id $\fb$ nous écrivons $x\equiv y\;\mod\;\fb$ comme 
abréviation
pour  $x\equiv y\;\mod$ \hbox{$(\fb=0)$}.
\begin{fact}
\label{factRecolTD}
Soit $\gT$ un \trdiz, $(\fa_i)_{i=1,\ldots n}$ une famille finie 
d'\ids principaux ($\fa_i=\dar s_i$)  de
$\gT$ et $\fa=\bigcap_i\fa_i$.
\begin{enumerate}
\item Si $(x_i)$ est une famille d'\elts de $\gT$ telle que pour 
chaque $i,j$ on
a $x_i\equiv x_j\;\mod\;\fa_i\vu\fa_j$, alors il existe un unique $x$ 
modulo
$\fa$ vérifiant:  $x\equiv x_i\;\mod\;\fa_i\;(i=1,\ldots ,n)$.
\item Notons $\gT_i=\gT/(\fa_i=0)$, 
$\gT_{ij}=\gT_{ji}=\gT/(\fa_i\vu\fa_j=0)$,
$\pi_i:\gT\to\gT_i$ et $\pi_{ij}:\gT_i\to\gT_{ij}$ les projections 
canoniques.
Si les $\fa_i$ recouvrent $\gT$, $(\gT,(\pi_i)_{i=1,\ldots n})$  est 
la limite
projective du diagramme $((\gT_i)_{1\leq i\leq n},(\gT_{ij})_{1\leq 
i<j\leq
n};(\pi_{ij})_{1\leq i\neq j\leq n})$ (voir la figure ci-après)
\item Soit maintenant $(\ff_i)_{i=1,\ldots n}$ une famille finie de 
filtres principaux,
notons $\gT_i=\gT/(\ff_i=1)$, 
$\gT_{ij}=\gT_{ji}=\gT/(\ff_i\cup\ff_j=1)$,
$\pi_i:\gT\to\gT_i$ et $\pi_{ij}:\gT_i\to\gT_{ij}$ les projections 
canoniques.
Si les $\ff_i$ recouvrent $\gT$, $(\gT,(\pi_i)_{i=1,\ldots n})$  est 
la limite
projective du diagramme $((\gT_i)_{1\leq i\leq n},(\gT_{ij})_{1\leq 
i<j\leq
n};(\pi_{ij})_{1\leq i\neq j\leq n})$.
\end{enumerate}
\end{fact}
 
 
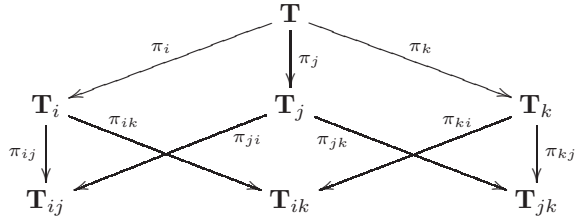
{\hspace*{10em}{
\xymatrix @R=2em @C=7em{
          &  \gT \ar[rd]^{\pi _{k}}\ar[d]^{\pi _{j}}\ar[ld]_{\pi _{i}}\\
 \gT _i\ar[d]_{\pi _{ij}}\ar@/-0.75cm/[dr]^{\pi _{ik}} &
     \gT _j\ar@/-.8cm/[dl]_{\pi _{ji}}\ar@/-.8cm/[dr]^{\pi _{jk}} &
        \gT _k\ar@/-0.75cm/[dl]_{\pi _{ki}}\ar[d]^{\pi _{kj}} &
\\
 \gT _{ij}  & 
    \gT _{ik}   & 
      \gT _{jk}   
}
}}

\begin{proof}
\textsl{1}. Il suffit de le démontrer avec $\fa=0$, ce qui est le point \textsl{2}.

\smallskip \noindent \textsl{2}. Soit $(\gH,(\psi_i)_{i\in I})$ la limite projective du diagramme.
On a un unique morphisme $$\varphi:\gT\to \gH$$ tel que $\varphi\circ \psi_i=\pi_i$
pour chaque~$i$. Et~$\varphi$ est injectif par hypothèse: $\varphi(x)=\varphi(y)$ implique $\varphi(x)\equiv\varphi(y) \mod\, (\fa_i=0)$ pour chaque $i$, et on a $\bigcap_i\fa_i=0$. 
On doit montrer qu'il est surjectif. \hbox{Soit $x=(x_i)_{i\in I}$} un \elt de $\gH$: \hbox{on a $x_i\in \gT_i$} pour chaque~$i$ et~\hbox{$\pi_{ij}(x_i)=\pi_{ji}(x_j)$} \hbox{pour $i\neq j$}. \hbox{Si $x_i=\pi_i(y_i)$}  on a donc dans $\gT$ la congruence 

\snic{y_i\equiv y_j \mod \dar (s_i\vu s_j).}

\noindent L'injectivité de $\varphi$ signifie que $\Vi_{i=1}^ns_i=0$. On a $\pi_i(y_i)=\pi_i(y_i\vu s_i)$ donc on peut supposer \hbox{que $y_i\geq s_i$}.
Les \egts \hbox{$\pi_{ij}(x_i)=\pi_{ji}(x_j)$} s'écrivent   

\snic{y_i\vu s_j\vu s_i=y_j\vu s_j\vu s_i,}

\noindent \cad $y_i\vu s_j=y_j\vu s_i$. Posons $y=\Vi_{i=1}^ny_i$.
\\
Alors, avec par exemple $j=1$, on obtient

\snic{\dsp y\vu s_1=y_1\vu\Vi_{i=2}^n(y_i\vu s_1)=y_1\vi\Vi_{i=2}^n(y_1\vu s_i)=y_1}

\noindent  (car $a\vi(a\vu b))=a$). Ainsi $\pi_j(y)=x_j$ pour chaque $j$.
Et $\varphi$ est bien surjective.
\end{proof}

Il y a aussi une procédure de recollement proprement dit.
Pour l'établir nous avons besoin du lemme suivant.

Rappelons que pour $s\in\gT$ le quotient $\gT/(s=0)$ est isomorphe au filtre principal $\uar s$ que l'on voit comme un \trdi dont l'\elt $0$ est $s$.
\begin{lemma}[Dans un \trdiz, les quotients principaux sont \gui{scindés}] \label{lemquoprinctrdi} ~ 
\\
Soit $\pi:\gT  \to \gT'$ un morphisme de \trdis et $s\in \gT$.
\Propeq
\begin{enumerate}
\item $\pi$ est un morphisme de passage au   quotient de $\gT$ par l'idéal principal $\fa=\dar s$.
\item Il existe un morphisme $\varphi:\gT'\to\,\uar s$ tel que
$\pi\circ \varphi=\Id_{\gT'}$.
\end{enumerate}
Dans ce cas $\varphi$ est uniquement déterminé par $\pi$ et $s$.\\
Naturellement, l'énoncé analogue \gui{renversé} est valable pour un quotient par un filtre principal. 
\end{lemma}
%
\begin{proof} 
\textsl{1} $\Rightarrow$ \textsl{2.} Soit $y\in \gT'$. On a $y=\pi(x)$ pour un $x\in \gT $. 
\\
On veut définir  $\varphi:\gT'\to\uar s$ par l'\egt $\varphi(y)=x\vu s$. Tout d'abord c'est bien défini: si $\pi(x)=\pi(x')$, alors $x\vu s=x'\vu s$ d'après le rappel précédent. Ensuite il est immédiat que $\varphi$ est un morphisme de \trdis et que $\pi\circ \varphi=\Id_{\gT'}$.

\smallskip \noindent 
\textsl{2} $\Rightarrow$ \textsl{1.} L'\egt $\pi\circ \varphi=\Id_{\gT'}$ implique
que $\pi$ est surjectif et que $\varphi$ est un \iso de $\gT'$ sur $\uar s$
avec la restriction de $\pi$ pour \iso réciproque. 
Ceci montre que $\varphi$ est uniquement déterminé par $\pi$ et $s$. On doit montrer l'\eqvc 

\snic{\pi(x_1)=\pi(x_2)\; \Leftrightarrow \;x_1\vu s=x_2\vu s.}

\noindent Comme $\varphi(0)=s$, on a $\pi(s)=0$, et $x_1\vu s=x_2\vu s$ implique $\pi(x_1)=\pi(x_2).$\\
Si $\pi(x_1)=\pi(x_2)$ alors  $\pi(x_1\vu s)=\pi(x_2\vu s)$, et puisque la restriction de $\pi$ à $\uar s$ est injective, cela implique
$x_1\vu s=x_2\vu s$. 
\end{proof}
\rem Nous avons utilisé en titre du lemme l'expression \gui{les quotients principaux sont scindés} par analogie avec les surjections scindées 
\hbox{entre \Amosz}, 
vue l'\egt $\pi\circ \varphi=\Id_{\gT'}$, mais l'analogie est limitée. Ici la \gui{section} $\varphi$ de $\pi$ est unique
(une différence importante),
et ce n'est \gui{pas tout à fait} un morphisme de~$\gT'$ dans~$\gT $ (une autre différence importante).

\begin{proposition}
\label{propRecolTD} \emph{(Recollement de \trdisz)}
  Supposons donnés un ensemble fini totalement ordonné~$I$ et dans la catégorie des \trdis  un diagramme

\snic{\big((\gT_i)_{i\in I},(\gT_{ij})_{i<j\in I},(\gT_{ijk})_{i<j<k\in I};
(\pi_{ij})_{i\neq j},(\pi_{ijk})_{i< j, j\neq k\neq i}\big)}

\noindent 
comme dans la figure ci-après, 
ainsi qu'une famille d'\elts 

\snic
{(s_{ij})_{i\neq j\in I}\in \prod\nolimits_{i\neq j\in I}\gT_{i}}

\noindent satisfaisant les conditions suivantes:
\begin{itemize}
\item le diagrammme est commutatif ($\pi_{ijk}\circ \pi_{ij}=\pi_{ikj}\circ \pi_{ik}$ pour tous $i$, $j$, $k$ distincts), 
\item pour $i\neq j$, $\pi_{ij}$ est un morphisme de passage au quotient par l'\id $\dar s_{ij}$,
\item pour $i$, $j$, $k$ distincts, $\pi_{ij}(s_{ik})=\pi_{ji}(s_{jk})$ et  $\pi_{ijk}$ est un morphisme de passage au quotient par \hbox{l'\id $\dar\pi_{ij}(s_{ik})$}.   
\end{itemize}

\smallskip {\hspace*{10em}
\xymatrix @R=2em @C=7em{
 \gT_i\ar[d]_{\pi _{ij}}\ar@/-0.75cm/[dr]^{\pi _{ik}} &
     \gT_j\ar@/-.8cm/[dl]_{\pi _{ji}}\ar@/-.8cm/[dr]^{\pi _{jk}} &
        \gT_k\ar@/-0.75cm/[dl]_{\pi _{ki}}\ar[d]^{\pi _{kj}} &
\\
 ~\gT_{ij}~ \ar[rd]_{\pi _{ijk}} & 
    ~\gT_{ik}~  \ar[d]^{\pi _{ikj}} & 
      ~\gT_{jk}~  \ar[ld]^{\pi _{jki}} 
\\
   &  ~\gT_{ijk}~ 
\\
}
}

\smallskip \noindent Alors si $\big(\gT\,;\,(\pi_i)_{i\in I}\big)$ est la limite projective du diagramme, les~\hbox{$\pi_i:\gT\to \gT_i$} forment un recouvrement par quotients principaux de $\gT$, et le diagramme est isomorphe à celui obtenu
dans le fait~\ref{factRecolTD}.
Plus précisément, il existe des $s_i\in\gT$ tels que chaque~$\pi_i$ est un morphisme de passage au quotient par l'\id $\dar s_i$ et $\pi_i(s_j)=s_{ij}$ pour tous $i\neq j$.

\noindent Le résultat analoque est valable pour les quotients par des filtres principaux.
\end{proposition}

\begin{proof}
Nous posons $s_{ii}=0$, $\gT_{ii}=\gT_i$, $\varphi_{ii}=\pi_{ii}=\Id_{\gT_i}$. Le lemme \ref{lemquoprinctrdi} nous donne des \gui{sections}  $\varphi_{ij}:\gT_{ij}\to \gT_i$ et~$\varphi_{ijk}:\gT_{ijk}\to \gT_{ij}$. 
\\
Les conditions imposées impliquent que les  
\idsz~\hbox{$\dar\pi_{jk}(s_{ji})$} \hbox{et $\dar\pi_{kj}(s_{ki})$} sont égaux, i.e.~\hbox{$\pi_{jk}(s_{ji})=\pi_{kj}(s_{ki})$}. 
\\
Pour $i\in I$, on définit $s_i\in \prod_k\gT_k$ par 
${s_i=(s_{ji})_{j\in I}}$, de sorte que $\pi_j(s_i)=s_{ji}$. Les coordonnées de $s_i$ sont compatibles (i.e. $s\in \gT$) car
$\pi_{jk}(s_{ji})=\pi_{kj}(s_{ki})$.
\\
Nous définissons ensuite une application $\varphi_i=\gT_i\to \prod_k\gT_k$ par

\snic{\varphi_i(x)=y=(y_j)_{j\in I} \hbox{ avec } y_j=\varphi_{ji}(x_j)=\varphi_{ji}\big(\pi_{ij}(x)\big).}

\noindent Montrons que les coordonnées de $y$ sont compatibles (i.e. $y\in \gT$). En effet

\snic{y_j=s_{ji}\vu y_j, \hbox{ donc } 
\pi_{jk}(y_j)= \pi_{jk}(s_{ji}\vu y_j)=\pi_{jk}(s_{ji})\vu\pi_{jk}(y_j), 
}

\noindent de m\^eme $\pi_{kj}(y_k)=\pi_{kj}(s_{ki})\vu\pi_{kj}(y_k)$.
Et puisque $\pi_{jki}$ est un morphisme de passage au quotient par l'\id $\dar\pi_{jk}(s_{ji})=\dar\pi_{kj}(s_{ki})$, l'\egt 

\snic{\pi_{jk}(y_j)=\pi_{kj}(y_k)}

\noindent peut \^etre testée en prenant les images par $\pi_{jki}$. \\
Or, puisque $\pi_{ji}(y_j)=\pi_{ji}\big(\varphi_{ji}(x_j\big)=x_j=\pi_{ij}(x)$, on obtient en utilisant  la commutativité du diagramme

\snic{\pi_{jki}\big(\pi_{jk}(y_j)\big)=\pi_{ijk}\big(\pi_{ji}(y_j)\big)=\pi_{ijk}\big(\pi_{ij}(x)\big).}

\noindent De m\^eme $\pi_{kji}\big(\pi_{kj}(y_k)\big)=\pi_{ikj}\big(\pi_{ik}(x)\big)$. Et nous concluons en utilisant une deuxième fois la commutativité du diagramme.\\
Une fois établi que $\varphi_i$ est bien une application $\gT_i\to \gT$, nous constatons facilement \hum{écrire les détails?} que $\pi_i\circ \varphi_i=\Id_{\gT_i}$, que l'image de
$\varphi_i$ est le filtre $\uar s_i$ de $\gT$ et que~$\varphi_i$ est un morphisme de
\trdis de $\gT_i$ sur le filtre $\uar s_i$. Donc, par le lemme \ref{lemquoprinctrdi}, $\pi_i$ est un morphisme de passage au quotient
par $\dar s_i$.
\end{proof}

\subsubsection*{Treillis de Heitmann}
\addcontentsline{toc}{subsubsection}{Treillis de Heitmann}

Un quotient intéressant, qui n'est ni un quotient par un \id ni un 
quotient
par un filtre, est le treillis de Heitmann.

\begin{lemma}
\label{lemHeT}
Sur un \trdi arbitraire $\gT$ la relation $\JT(a)\subseteq\JT(b)$
est une relation de préordre $a\preceq b$ qui définit un quotient 
de $\gT$.
On a aussi:
\begin{equation} \label{eqJaJb}
a\preceq b \quad \Longleftrightarrow\quad  a\in \JT(b)
\quad \Longleftrightarrow\quad \forall x\in\gT \; (a\vu x = 1 \Rightarrow 
b\vu x=1)
\end{equation}
\end{lemma}
\begin{proof}
Les équivalences $a\in \JT(b) \;\Leftrightarrow\;
\JT(a)\subseteq\JT(b)\;\Leftrightarrow\;\forall x\in\gT \;(a\vu x = 1
\Rightarrow b\vu x=1)$ résultent de ce qui a été dit page 
\pageref{eqRJb}
concernant le radical de Jacobson d'un \id (voir 
l'\egtz~(\ref{eqRJb})).\\
Par ailleurs on vérifie sans difficulté les relations 
(\ref{eqPreceq})
nécessaires pour qu'un préordre définisse un quotient.
\end{proof}

\begin{definition}
\label{defHeT}
On appelle \textsl{treillis de Heitmann de $\gT$} et on note $\He(\gT)$ 
le
treillis quotient de~$\gT$ obtenu en remplaçant sur $\gT$ la 
relation
d'ordre $\leq_\gT $  par la relation de préordre 
$\preceq_{\He(\gT)}$
définie comme suit

\vspace{-1em}
\begin{equation} \label{eqdefHeT}
\begin{array}{rcl}\qquad 
a\preceq_{\He(\gT)} b & \equidef  &   \JT(a)\subseteq\JT(b)  \quad \hbox{(cf. \dfn \ref{defJac})}
  \end{array}
  \end{equation}
Ce treillis quotient peut être identifié à l'ensemble des 
\ids $\JT(a)$,
avec la projection canonique
$$ \gT\longrightarrow \He(\gT),\quad a\longmapsto \JT(a)$$

\end{definition}

Notez qu'avec l'identification précédente on a les \egtsz:
\begin{equation} \label{eqHeT2}
\JT(a\vi b)=\JT(a)\vi_{\He(\gT)}\JT(b),\quad
\JT(a\vu b)=\JT(a)\vu_{\He(\gT)}\JT(b)
\end{equation}

Dire que le treillis $\gT$ est faiblement Jacobson revient à dire 
que
$\gT=\He(\gT)$.

\smallskip Le lemme suivant est une précision (et une généralisation) de la
première \egt ci-dessus. Il nous sera utile dans la suite.

\begin{lemma}
\label{lemJacInter}
Si $\fa$ et $\fb$ sont deux \ids de $\gT$, on a
$\JT(\fa\cap\fb)=\JT(\fa)\cap\JT(\fb)$.
\end{lemma}
\begin{proof}
Il suffit de montrer que si $z\in\JT(\fa)\cap\JT(\fb)$ alors
$z\in\JT(\fa\cap\fb)$. Soit $t\in\gT$ tel que $z\vu t=1$, nous 
cherchons
$c\in\fa\cap\fb$ tel que $c\vu t=1$.
Or nous avons un  $a\in\fa$ tel que $a\vu t=1$ et  un  $b\in\fb$ tel 
que $b\vu
t=1$. Il suffit donc de prendre $c=a\vi b$.
\end{proof}

On notera que la preuve ne marcherait pas pour une intersection 
infinie d'\idsz.

\hum{Il est probable que les \ids de $\He(\gT)$ sont exactement les
$\pi(\JT(\fa))$. }

\begin{fact}
\label{factHeHe} Soit $\gT$ un \trdiz, $\gT'=\gT/(\JT(0)=0)$, 
$x\in\gT$ et $\fa$
un idéal.
\begin{enumerate}
\item $x=_{\He(\gT)}1\;\Longleftrightarrow\; x=1$.
\item $x=_{\He(\gT)}0\;\Longleftrightarrow\; x\in\JT(0)$.
\item $\He(\He(\gT))\;=\;\He(\gT')\;=\;\He(\gT)$.
\item Si $\gL=\gT/(\fa=0)$,
$\He(\gL)$ s'identifie à $\He(\gT)/(\JT(\fa)=0)$.
\end{enumerate}
\end{fact}

\rem On notera cependant que $\He$ ne définit pas un foncteur.

\begin{proof}
Les points \textsl{1} et \textsl{2} sont immédiats.  \\
Le point \textsl{4} est laissé \alec.  Il implique
$\He(\gT')\;=\;\He(\gT)$.\\
Dans le point \textsl{3} les treillis $\He(\He(\gT))$ et $\He(\gT')$ sont identifiés à
des quotients de $\gT$.  Montrons l'égalité $\He(\He(\gT))=\He(\gT)$,
\cad que pour $a,b\in\gT$, $a\preceq_{\He(\He(\gT))}b\Rightarrow
a\preceq_{\He(\gT)}b$.  Par \dfn l'hypothèse signifie: $\Tt x \in
\gT\;(a\vu x=_{\He(\gT)}1\;\Rightarrow \;b\vu x=_{\He(\gT)}1)$.  Or
d'après le point~\textsl{1} cela veut dire $\Tt x \in \gT\;(a\vu
x=1\;\Rightarrow \;b\vu x=1)$, \cad $a\preceq_{\He(\gT)}b$.
\end{proof}

\subsection{Algèbres de Heyting, de Brouwer, de Boole}\label{subsecAgHagB}

\subsubsection*{Algèbres de Heyting}
\addcontentsline{toc}{subsubsection}{Algèbres de Heyting}

Un \trdi $\gT$ est appelé un {\sl treillis implicatif} (\cite{Cur73}) 
ou une {\sl \agHz} (\cite{Joh1986}) lorsqu'il existe
une opération binaire  $\im$ vérifiant pour tous $a,\,b,\,c$:
\begin{equation} \label{eqAgHey}
a\vi b \leq c \;\;\Longleftrightarrow \;\; a \leq  (b\im c)\,
\end{equation}
Ceci signifie que pour tous $b, c\in\gT$, l'\id $c:b$  est principal, 
son \gtr
étant noté $b\im c$.
Donc si elle existe, l'opération $\im$ est déterminée de 
manière unique
par la structure du treillis.
On définit alors la loi unaire  $\neg x = x\im 0$.
La structure d'\agH peut être définie comme purement 
équationnelle en
donnant de bons axiomes. Précisément un treillis $\gT$ (non 
supposé
distributif) muni d'une loi $\im$ est une \agH \ssi les axiomes 
suivants sont
vérifiés (cf. \cite{Joh1986}):
$$\begin{array}{rcl}
a\im a&=   &1    \\
a\vi(a\im b)&=   &a\vi b    \\
b\vi(a\im b)&=   & b   \\
a\im(b\vi c)&=   &(a\im b)\vi(a\im c)
\end{array}$$
Notons aussi les faits importants suivants:

$$\begin{array}{rcl}
(a\vu b)\im c &=& (a\im c)\vi(b\im c)    \\
\neg(a\vu b)&=   & \neg a\vi \lnot b   \\
 a&\leq    &\neg\neg a   \\
\neg a\vu b&\leq    & a\im b   \\
a\leq b&\Leftrightarrow& a\im b =1
\end{array}$$

  Tout \trdi fini est une \agHz, car tout \itf est principal.

\smallskip Un cas particulier important d'\agH est une \textsl{\agBz}:
c'est un \trdi dans lequel tout \elt $x$ possède \textsl{un
complément}, \cad un \elt $y$ vérifiant $y\vi x=0$ et $y\vu x=1$
($y$ est noté $\lnot x$ et l'on a $a\im b=\lnot a\vu b$).

\smallskip Un \textsl{\homo d'\agHsz} est un \homo $\varphi :\gT\to\gT'$
de \trdis qui vérifie $\varphi(a\im b)=\varphi(a)\im\varphi(b)$ pour
tous $a,b\in\gT$.

\smallskip Le fait suivant est immédiat.

\begin{fact}
\label{factQuoAgH}
Soit $\pi:\gT\to\gT'$ un \homo de \trdisz. Supposons que $\gT$ et 
$\gT'$ sont
deux \agHs et notons $a\preceq b$ pour 
$\varphi(a)\leq_{\gT'}\varphi(b)$. Alors
$\pi$ est un \homo d'\agHs \ssi on a pour tous $a,a',b,b'\in\gT$:
$$
a\preceq a'\Rightarrow (a'\im b)\preceq(a\im b)  \qquad 
\mathrm{et}\qquad
b\preceq b'\Rightarrow (a\im b)\preceq(a\im b')
$$
\end{fact}

On a aussi:

\begin{fact}
\label{factQuoAgH2}
Si $\gT$ est une \agH tout quotient $\gT/(y=0)$ (\cad tout quotient 
par un \id
principal) est aussi une \agHz.
\end{fact}
\begin{proof}
Soit $\pi:\gT\to\gT'=\gT/(y=0)$ la projection canonique. On a
$\pi(x)\vi\pi(a)\,\leq_{\gT'}\, \pi(b)\;\Leftrightarrow\; \pi(x \vi
a)\,\leq_{\gT'}\, \pi(b)\;\Leftrightarrow\; x\vi a \,\leq\, b\vu
y\;\Leftrightarrow\; x\,\leq\, a\im(b\vu y)$. Or $y\,\leq\, b\vu 
y\,\leq\,
a\im(b\vu y)$, donc $\pi(x)\vi\pi(a)\,\leq_{\gT'}\, 
\pi(b)\;\Leftrightarrow\;
x\,\leq\, (a\im(b\vu y))\vu y$, \cad $\pi(x)\leq_{\gT'}\pi(a\im(b\vu 
y))$, ce
qui montre que $\pi(a\im(b\vu y))$ vaut pour $\pi(a)\im\pi(b)$ dans 
$\gT'$.
\end{proof}

\hum{

1.  Cependant il ne semble pas que $\pi$ soit en général un \homo
d'\agHsz.

2.  Il serait bon d'avoir un exemple d'un treillis quotient d'une \agH
qui ne serait pas une \agHz.

3.  De manière générale, il serait bon d'avoir des exemples
variés d'\agH non noethériennes à notre disposition.

}

\medskip \rem 
La notion d'\agH est reminiscente de la notion d'anneau
cohérent en \alg commutative.  En effet un anneau cohérent peut
être caractérisé comme suit: l'intersection de deux \itfs est un
\itf et le transporteur d'un \itf dans un \itf est un \itfz.  Si l'on
\gui{relit} ceci pour un \trdi en se rappelant que tout \itf est
principal on obtient une \agHz.

\medskip \rem 
Tout \trdi $\gT$ engendre une \agH de façon
naturelle.  Autrement dit on peut rajouter formellement un \gtr pour
tout \id $b:c$.  Mais si on part d'un \trdi qui se trouve être une
\agHz, l'\agH qu'il engendre est strictement plus grande.  Prenons par
exemple le treillis $\Trois$ (fini donc c'est une \agHz), qui est le \trdi libre à un \gtrz. 
L'\agH qu'il engendre est donc l'\agH libre à un générateur.  Or
celle-ci est infinie (cf.  \cite[section 4.11]{Joh1986}).  A contrario le treillis
booléen engendré par $\gT$ (cf.  \cite{CC00}, \cite[Théorème XI-1.8]{ACMC}) reste égal à $\gT$
lorsque celui-ci est booléen.

\subsubsection*{Treillis avec négation}
\addcontentsline{toc}{subsubsection}{Treillis avec négation}

Un \trdi \textsl{possède une négation} si pour tout $x$ l'idéal 
$(0:x)$ est
principal, engendré par un \elt que l'on note  $\lnot x$.
Les règles suivantes sont immédiates.

$$\begin{array}{rclcrcl}
x\vi y =0&\Leftrightarrow& y\leq \lnot x&,&a\leq b&  \Rightarrow  & \lnot b   \leq     \lnot a \\
a& \leq   & \lnot\lnot a  &,   & \lnot a  &=& \lnot\lnot\lnot a \\
\lnot(a\vu b)& =   &  \lnot a\vi\lnot b &  , & \lnot a\vu\lnot b 
&\leq &
\lnot(a\vi b)\\
\lnot(x\vu\lnot x)&=&0&,&\lnot\lnot(x\vu\lnot x)&=&1
\end{array}$$

Si pour tout $a$, $\lnot\lnot a=a$, le treillis est une \agB parce 
qu'alors
$x\vu\lnot x=1$.

\begin{fact}
\label{factSpecMin}
Si $\gT$ possède une négation,   notons $\Fmin(\gT)=\ff$ le filtre
engendré par les $x\vu\lnot x$. Alors 
$\He(\gT\cir)=(\gT/(\ff=1))\cir$, et ce
treillis est une \agBz.
\end{fact}
\begin{proof}
Il est clair que $\lnot x$ est un complément de $x$ dans 
$\gT/(\ff=1)$, ce
treillis est donc une \agBz. En présence de la négation, la 
relation
$a\leq_{\He(\gT\cir)}b$ est équivalente à $\lnot a\leq \lnot b$
et ceci est facilement équivalent à $b\leq a\;\mod\;(\ff=1)$.
\end{proof}

\begin{fact}
\label{factWJavecneg}
Si $\gT$ est un treillis avec négation, le treillis $\gT\cir$ est 
faiblement
Jacobson \ssi $\gT$ est une \agBz.
\end{fact}
\begin{proof}
En présence de négation, les équations (\ref{eqJaJb}) et 
(\ref{eqdefHeT})
donnent pour $a\leq_{\He(\gT\cir)}b$ la condition équivalente 
$\lnot b\leq
\lnot a$.
Le treillis  $\gT\cir$ est donc faiblement Jacobson \ssi $\lnot b\leq 
\lnot a$
implique  $a\leq b$. En particulier on doit avoir $a=\lnot\lnot a$ 
pour tout
$a$.
\end{proof}

\subsubsection*{Algèbres de Brouwer}
\addcontentsline{toc}{subsubsection}{Algèbres de Brouwer}

  Un \trdi tel que le treillis opposé est une \agH est appelé une
  \textsl{\alg de Brouwer}.  C'est un \trdi dans lequel tous les filtres
  $c\setminus b$ sont principaux.  On note alors $c-b$ le \gtr de
  $c\setminus b$.

En passant au treillis opposé le fait suivant dit la même chose
que le fait~\ref{factSpecMin}.

\begin{fact}
\label{factSpecMax}
On dit que \textsl{le treillis $\gT$ possède un complément de Brouwer} lorsque pour tout $x$ 
le filtre
$(1\setminus x)$ est principal. Il est alors engendré par un unique \elt que l'on note  
$1- x$. Dans ce cas,
notons $\Imax(\gT)=\fa$ l'\id engendré par les $x\vi(1-x)$. Alors 
$\He(\gT)$ est égal à l'\agB $\gT/(\fa=0)$.
\end{fact}

Nous laissons \alec le soin de traduire le fait 
\ref{factWJavecneg} lorsque l'on renverse la relation d'ordre.

\subsection{Treillis distributifs noethériens}

En \clamaz, pour un \trdi $\gT$ \propeq
\begin{itemize}
\item [$(1)$] Tout \id de $\gT$ est principal.
\item [$(2)$] Toute suite croissante d'\elts de $\gT$ est stationnaire.
\item [$(3)$] Toute suite croissante d'\ids  de $\gT$ est stationnaire.
\end{itemize}

Un tel treillis est appelé \textsl{noethérien} (par analogie avec
l'\alg commutative, on pourrait aussi l'appeler \textsl{principal}).  C'est
clairement est une \agH (en \clamaz).

Tout sous-treillis et tout treillis quotient d'un treillis
noethérien est noethérien.

En \coma la notion est plus délicate.  Aucun treillis non trivial ne
vérifie le point (2) (qui est à priori la formulation la plus faible
des trois).  On pourrait définir un \trdi noethérien comme un
treillis vérifiant une condition \gui{ACC \covz} du style: toute
suite croissante admet deux termes consécutifs égaux.  Cette
condition est équivalente à (2) en \clamaz.  Mais il y a à priori
plusieurs variantes intéressantes.

En pratique, on est en général intéressé par le
fait que certains \ids bien précis sont principaux, comme dans le
cas des \agHsz.  Or le fait qu'un treillis est une \agH ne résulte
pas \cot de la condition ACC \cov (de la même manière, en \alg
commutative, la cohérence, qui est souvent plus importante que la
noethérianité, ne résulte d'aucune variante \cov connue de la
noethérianité).  Voir à ce sujet la proposition~\ref{propZarHeyt}.

\medskip \rem
Montrons en \clama que si $\gT$ et $\gT\cir$ sont noethériens 
alors~$\gT$ est fini.  Les \idemas sont des $\dar x$ où $x$ est un
prédécesseur immédiat de $1$.  Et le spectre maximal est fini,
parce que si $(\fm_n)=(\dar x_n)$ est une suite infinie d'\idemasz, la
suite $(\Vi_{i\leq n}x_i)$ est strictement décroissante.  On peut
ensuite appliquer le résultat à chacun des treillis quotients par
les \idemasz.  On termine par le lemme de K\"onig.  Rendre cette
preuve \covz, avec une \dfn \cov suffisamment forte de la
noethérianité est un défi intéressant.

\section{Espaces spectraux}
\label{secESSP}

\subsection{Généralités}

\subsubsection*{En \clamaz }
\addcontentsline{toc}{subsubsection}{En \clamaz }

Un {\sl \id premier} $\fp$ d'un treillis $\gT\neq \Un$ est un \id dont
le complé\-mentaire $\ff$ est un filtre (qui est alors un {\em
filtre premier}).  On a alors $\gT/(\fp=0,\ff=1)\simeq\Deux$.  Il
revient au même de se donner un \idep de $\gT$ ou un morphisme de
\trdis $\gT\rightarrow \Deux$.

Dans cette section, nous noterons $\theta_\fp:\gT\to\Deux$ l'\homo
associé à l'\idep $\fp$.

On vérifie facilement que si $S$ est une partie génératrice du
\trdi $\gT$, un \idepz~$\fp$ de $\gT$ est complètement
caractérisé par sa trace sur $S$ (cf.  \cite{CC00}).

Un \textsl{\idema} (resp.  \textsl{premier minimal}) est un \id maximal
parmi les \ids stricts (resp.  minimal parmi les \idepsz).  Il revient
au même de dire que $\fm$ est maximal ou que
$\gT/(\fm=0)\simeq\Deux$, les \idemas sont donc premiers.  Il revient
au même de dire que $\fp$ est un \idep minimal ou que son
complémentire est un filtre maximal.

En \clama tout \id strict est contenu dans un \idema et (par
renversement) tout filtre strict est contenu dans le complémentaire d'un
\idep minimal.

\smallskip Le \textsl{spectre d'un \trdi $\gT$} est l'ensemble $\Spec
\,\gT$ de ses \idepsz, muni de la topologie suivante: une base
d'ouverts est donnée par les 
$$
\DT(a)\eqdefi\sotq{\fp\in\Spec
\,\gT}{a\notin\fp},\quad a\in \gT.
$$
On vérifie que
\begin{equation} \label{eqDa}
\left.\begin{array}{rclcrcl}
  \DT(a\vi b)   & =  & \DT(a)\cap \DT(b) ,&\quad & \DT(0)  & =  & 
\emptyset  ,\\
  \DT(a\vu b)   & =  & \DT(a)\cup \DT(b) ,&&  \DT(1) & =  &  
\Spec\,\gT.
  \end{array}
\right\}
\end{equation}

Le complémentaire de $\DT(a)$ est un fermé qu'on note $\VT(a)$.

On étend la notation $\VT(a)$ comme suit: si $I\subseteq\gT$, on
pose $\VT(I)\eqdefi\bigcap_{x\in I}\VT(x)$.  Si $\cI_\gT(I)=\fII$, on
a $\VT(I)=\VT(\fII)$.  On dit parfois que $\VT(I)$ est \textsl{la
variété associée à $I$}.

\medskip\noindent
{\bf Définition.} 
Un espace topologique homéomorphe à un espace $\Spec(\gT)$
est appelé un \textsl{espace spectral}. Les espaces spectraux proviennent de l'étude de Stone (\cite{Sto37}).

\medskip Johnstone les appelle des \textsl{espaces cohérents} (\cite{Joh1986}). C'est Hochster qui les a baptisés dans \cite{Hoc1969}.

Avec la logique classique et l'axiome du choix, l'espace $\Spec \,\gT$
a \gui{suffisamment de points}: on peut retrouver le treillis $\gT$
à partir de son spectre.  Voici comment.

Tout d'abord on a le


\medskip\noindent
{\bf \Tho de Krull }\label{ThKrull} (en \clamaz)\\ 
{\sl Supposons que $\fJ$ est un \idz, $\fF$ un filtre et 
$\fJ\cap\fF=\emptyset$.
Alors il existe un \idep $\fP$ tel \hbox{que $\fJ\subseteq\fP$} et
$\fP\cap\fF=\emptyset$.
  }

\medskip
On en déduit que:
\begin{itemize}
\item L'application $a\in\gT\,\mapsto\,\DT(a)\in\cP(\Spec\,\gT)$
est injective: elle identifie $\gT$ à un treillis d'ensembles 
(\textsl{\tho de représentation de Birkhoff}).
\item Si $\varphi : \gT\to\gT'$ est un \homo injectif l'application
$\varphi^\star:\Spec\,\gT'\to\Spec\,\gT$ obtenue par dualité est 
surjective.
\item Tout \id de $\gT$ est intersection des \ideps qui le 
contiennent.
\item L'application $\fII\mapsto \VT(\fII)$, des \ids de $\gT$ vers 
les fermés
de $\Spec\,\gT$ est un \iso d'ensembles ordonnés (pour l'inclusion 
et
l'inclusion renversée).
\end{itemize}

On montre aussi que les \oqcs de $\Spec \,\gT$ sont exactement les
$\DT(a)$.  D'après les \egts (\ref{eqDa}) les \oqcs de $\Spec \,\gT$
forment un \trdi de parties de~$\Spec \,\gT$, isomorphe à $\gT$.

\`A partir d'un espace spectral $X$ on peut considérer le \trdi
$\OQC(X)$ formé par ses \oqcsz.  Puisque pour tout \trdi $\gT$,
$\OQC(\Spec(\gT))$ est canoniquement isomorphe à $\gT$, pour tout
espace spectral $X$, $\Spec(\OQC(X))$ est canoniquement 
homéo\-morphe~à~$X$.

\smallskip Tout \homo $\varphi :\gT\rightarrow \gT'$ de \trdis fournit
par dualité une application continue $\varphi^\star:\Spec
\,\gT'\rightarrow \Spec \,\gT$, qui est appelée une
\textsl{application spectrale}.  Pour qu'une application continue entre
espaces spectraux soit spectrale il faut et il suffit que l'image
réciproque de tout \oqc soit un \oqcz.

L'article fondateur  \cite{Sto37} démontre pour l'essentiel que la catégorie spectrale ainsi définie est antiéquivalente à
celle des \trdis \cite[{II-3.3}, coherent locales]{Joh1986}.
Plus précisément, cet énoncé qui semble ici tautologique devient non trivial lorsque l'on donne une \dfn des espaces spectraux en termes purement d'espaces topologiques, comme dans la remarque qui suit.
Pour plus de détails sur cette antiéquivalence, on peut se reporter au théorème de Krull page~\pageref{ThKrull}, à \cite[\hbox{section V-8}]{BW74},
à \cite{CL2001-2018} et à l'article de synthèse \cite{Lom2020}.

\medskip 
\rem
Une \dfn purement topologique des espaces spectraux 
est la sui\-vante~\cite{Sto37}.
\begin{itemize}
\item L'espace est de Kolmogoroff (i.e., de type $\mathrm{T}_0$): 
étant donnés deux points il existe un voisinage de l'un des deux qui ne contient pas l'autre.
\item L'espace est \qcz.
\item L'intersection de deux \oqcs est un \oqcz.
\item Tout ouvert est réunion d'\oqcsz.
\item Pour tout fermé $F$ et pour tout ensemble $S$ d'\oqcs tels 
que 
$$\textstyle F\cap
\bigcap_{U\in S'} U\neq \emptyset\,\hbox{   pour toute partie finie  }\,S'
\,\hbox{  de  }\,S
$$ 
on a aussi
$F\cap \bigcap_{U\in S} U\neq \emptyset$.
\end{itemize}
En présence des quatre premières propriétés la dernière 
peut se
reformuler comme suit (\cite{Hoc1969}). 
\begin{itemize}
\item Tout fermé irréductible\footnote{Un fermé qui n'est pas réunion de deux fermés strictement plus petits} admet un point
générique.
\end{itemize}


\subsubsection*{Points génériques, relation d'ordre}
\addcontentsline{toc}{subsubsection}{Points génériques, relation d'ordre}

On dit qu'un point $x \in X$ d'un espace spectral est le
\textsl{point générique du fermé $F$} \hbox{si $F=\ovs{x }$}.  Ce point
(quand il existe) est nécessairement unique car les espaces spectraux
sont des espaces de Kolmogoroff.  Les fermés $\ovs{x }$
sont exactement tous les fermés irréductibles de $X$.  La relation
d'ordre $y\in\ovs{x}$ sera notée $x\leq_X y$.

Lorsque $X=\Spec\,\gT$ la relation $\fp\leq_X \fq$ est simplement la
relation d'inclusion usuelle entre \ideps du \trdi $\gT$.

Les points fermés de $\Spec\,\gT$ sont les \idemas de $\gT$.

\medskip 
On appelle \textsl{espace de Stone}\footnote{La terminologie ne semble pas entièrement fixée.  \cite{BW74} appellent espace de Stone un espace topologique qui est à très peu près un espace spectral. Leur but est une catégorie d'espaces topologiques antiéquivalente à celle des \trdis \gui{non bornés}, i.e., sans $0$ et $1$.}  un espace spectral dont le treillis des \oqcs est une \agBz\footnote{Il est homéomorphe à un espace $\Spec\,\gB$ pour une \agB \(\gB\)}.   Il est bien connu que les espaces de
Stone peuvent être caractérisés comme les espaces compacts
totalement discontinus.

\subsubsection*{En \comaz}
\addcontentsline{toc}{subsubsection}{En \comaz}

 D'un point de vue \cofz, $\gT$ est une version \gui{sans points} de
 $\Spec\,\gT$.  En d'autres termes, à défaut d'avoir accès aux
 points de $\Spec\,\gT$, on peut se contenter de l'ensemble de ses
 \oqcsz, qui sont directement visibles (sans recours à l'axiome du
 choix ni au principe du tiers exclu).  La version sans points est
 plus facile à appréhender.  Au contraire les points de
 $\Spec\,\gT$ ne sont pas en général des objets accessibles sans
 recours à des principes non constructifs.

En \coma on a à priori plusieurs possibilités pour définir le
spectre d'un \trdi (toutes équivalentes en \clamaz).
Le plus raisonnable semble de définir $\Spec\,\gT$ comme l'ensemble
des filtres premiers de $\gT$, \cad les filtres pour lesquels on a
$$x\vi y\in\fF\quad \Longrightarrow \quad x\in\fF \;\;\mathrm{ou}\;\;
y\in\fF $$ avec un \gui{ou} explicite.  Mais de tels espaces
$\Spec\,\gT$ n'ont pas toujours suffisamment de points{\footnote{~On
peut par exemple définir un \trdi infini dénombrable explicite qui
ne possède pas d'\ideps récursifs.  Pour un tel \trdiz, il ne peut
pas y avoir de \prco que $\Spec\,\gT$ est non vide}} et on ne peut pas
affirmer \cot que les deux catégories sont antiéquivalentes, du
moins si l'on définit les morphismes entre espaces spectraux comme des
applications, car les applications nécessitent des points.

Une solution alternative satisfaisante (mais un peu troublante au
premier abord) est de considérer $\Spec\,\gT$ comme un \gui{espace
topologique sans points}, \cad un espace topologique défini
uniquement à travers sa base d'ouverts $\DT(a)$ (où $a$ parcourt
$\gT$).  Les morphismes sont alors définis de manière purement
formelle comme donnés par les morphismes des treillis
correspondants, en renversant le sens des flèches.  De ce point de
vue l'antiéquivalence de la catégorie spectrale et de la
catégorie des \trdis devient une pure tautologie définitionnelle.

En tout état de cause, bien que la catégorie spectrale reste utile
pour l'intuition, tout le travail se fait dans la catégorie des
\trdisz.  L'avantage est naturellement que l'on obtient des \thos
\cofsz.

Dans cet article les spectres seront étudiés uniquement du point
de vue des \clamaz, comme source d'inspiration importante pour de
bonnes notions concernant les \trdisz.

\subsubsection*{Espaces spectraux noethériens}
\addcontentsline{toc}{subsubsection}{Espaces spectraux noethériens}
Un espace topologique $X$ est dit \textsl{noethérien} si toute suite
croissante d'ouverts est stationnaire.  Il revient au même de dire
que tout ouvert est \qcz. Pour un espace spectral, il est équivalent de dire que le treillis $\OQC(X)$ est noethérien.  Dans un espace spectral noethérien tout ouvert est un $\DT(a)$ et tout fermé un $\VT(b)$.

\hum{Je serais curieux de connaître une formulation sans point de
la propriété suivante, plus faible que la noethérianité: toute
suite croissante pour $\leq _X$ est stationnaire.}

\subsubsection*{Deux autres topologies intéressantes sur 
$\Spec\,\gT$}
\addcontentsline{toc}{subsubsection}{Deux autres topologies intéressantes sur 
$\Spec\,\gT$}

\hum{Par rapport à l'article original, l'explication ci-dessous est allongée.}
En \clama on a une bijection canonique entre les ensembles
sous-jacents aux espaces $\Spec\,\gT$ et $\Spec\,\gT\cir$: à un
\idep de $\gT$ on associe le filtre premier complémentaire, qui est un \idep de $\gT\cir$. Cela permet d'identifier ces
deux ensembles, même si parfois l'effet n'est pas très heureux.  
Une fois les ensembles sous-jacents identifiés, la topologie n'est pas la même.  Les ouverts de base de $\Spec\,\gT\cir$ sont les $\DTo(a)=\VT(a)$.  
Modulo cette identification, pour $X=\Spec\,\gT$ et
$X'=\Spec\,\gT\cir$, la relation d'ordre $\leq_{X'}$ est la relation
opposée~à~$\leq_X$ (l'ordre est renversé), mais ce qui se
passe pour la topologie est plus compliqué.

\smallskip 
On doit également considérer la \textsl{topologie constructible}
(en anglais: patch topology) dont les ouverts de base sont les $\DT(a) \cap
\VT(b)$.  Cela donne un espace compact naturellement homémorphe à
$\Spec\,\gT^{\rm bool}$ où $\gT^{\rm bool}$ est le treillis
booléen engendré par $\gT$.  En \clama on obtient $\gT^{\rm bool}$
comme la sous-\agB de l'ensemble des parties de $\Spec\,\gT$
engendrée par les $\DT(a)$.  Ce treillis peut aussi être décrit
\cot comme suit (cf.  \cite{CC00}).  On considère une copie disjointe
de $\gT$, que l'on note $\dot{\gT}$.  Alors $\gT^{\rm bool}$ est un \trdi
défini par \gtrs et relations.  Les \gtrs sont les \elts de
l'ensemble $T_1=\gT\cup\dot{\gT}$ et les relations sont obtenues comme
suit~: si $A,F,B,E$ sont quatre parties finies de $\gT$ on a
$$
  \Vi A \vi \Vi E\leq_\gT\Vu B\vu \Vu F
\quad \Longrightarrow \quad
\Vi A \vi \Vi \dot{F} \leq_{T_1} \Vu B \vu \Vu \dot{E}
$$
On montre que
  $\gT$ et $\dot{\gT}$ s'injectent naturellement dans $\gT^{\rm bool}$
et que l'implication ci-dessus est en fait une équivalence.
On obtient par dualité deux applications spectrales bijectives
$\Spec\,\gT^{\rm bool}\to\Spec\,\gT$ et $\Spec\,\gT^{\rm
bool}\to\Spec\,\gT\cir$.

\subsubsection*{Espaces spectraux finis}
\addcontentsline{toc}{subsubsection}{Espaces spectraux finis}
 
 Dans le cas des \trdis \textsl{finis} on obtient les espaces spectraux
 finis, qui ne sont rien d'autre que les ensembles ordonnés finis,
 (car il suffit de connaître l'adhérence des points pour
 connaître la topologie) avec pour base d'ouverts les $\dar a$. 
 Les ouverts sont tous \qcsz, ce sont les parties initiales, et les
 fermés sont les parties finales.  Enfin, une application entre
 espaces spectraux finis est spectrale \ssi elle est croissante (pour
 les relations d'ordre associées).

La notion d'espace spectral apparaît ainsi comme une
généralisation pertinente de la notion d'ensemble ordonné fini
au cas infini. Voir \cite[Théorème XI-5.6, dualité entre ensembles ordonnés finis et \trdis finis]{ACMC}.

Dans le cas fini, si l'on identifie les ensembles sous-jacents à $\Spec\,\gT$ et
$\Spec\,\gT\cir$ les deux spectres sont presque les mêmes: c'est le
même ensemble ordonné au renversement près de la relation
d'ordre.  En outre les ouverts et les fermés sont simplement échangés.

\subsection{Treillis quotients et sous-espaces spectraux}\label{secSESP}

\subsubsection*{Caractérisation des sous-espaces spectraux}
\addcontentsline{toc}{subsubsection}{Caractérisation des sous-espaces spectraux}

En utilisant l'antiéquivalence des catégories, on  pourrait définir
directement la notion de \textsl{sous-espace spectral} comme la notion duale de la notion de treillis quotient.  Le \tho \ref{propSESP} explique cela en détail. 

Nous commençons par un lemme facile, qui caractérise les points
de $\Spec\,\gT$ qui \gui{sont des \elts de $\Spec\,\gT'$} lorsque
$\gT'$ est un quotient de $\gT$.

\begin{lemma}
\label{lemSESP}
Soit $\gT'$ un treillis quotient de $\gT$ et $\pi:\gT\to\gT'$ la
projection canonique.  Notons $X=\Spec\,\gT'$, $Y=\Spec\,\gT$ et
$\pi^\star:X\to Y$ l'injection duale de $\pi$.  Rappelons que pour un
\idep $\fp$ de $\gT$ nous notons $\theta_\fp:\gT\to\Deux$ l'\homo
correspondant de noyau $\fp$.
\Propeq
\begin{itemize}
\item $\fp\in\pi^\star(\Spec\,\gT').$
\item $\theta_\fp$ se factorise par $\gT'.$
\item $\Tt a,b\in\gT\;((a\preceq b,\,b\in\fp)\Rightarrow a\in\fp).$
\end{itemize}
Cela peut se reformuler comme suit. Si le treillis quotient $\gT'$ 
est défini
par un système $R$ de relations $x_i=y_i$, \propeq
\begin{itemize}
\item $\fp\in\pi^\star(\Spec\,\gT').$
\item $\theta_\fp$ \gui{réalise un modèle de $R$}, \cad $\Tt i\;\;
\theta_\fp(x_i)=\theta_\fp(y_i).$
\item  $\Tt i\;\; (x_i\in\fp\;\Leftrightarrow\; y_i\in\fp).$
\end{itemize}
\end{lemma}

Dans le \tho suivant  nous identifions $\Spec\,\gT'$ à une partie de
$\Spec\,\gT$ au moyen de l'injection $\pi^\star$.
Des résultats analogues énoncés dans un langage un peu différent se trouvent dans \cite[section~3]{Esc2001}\footnote{Escard{\'o} écrit son article dans le langage des locales. Il parle de patch topology plutôt que de topologie constructible. Si $Y=\Spec\,\gT$, il note $\Patch\, Y$ pour l'espace de Stone $\Spec\,\gT^{\rm bool}$.}.

\begin{theorem}
\label{propSESP} \emph{(Définition et caractérisations des sous-espaces 
spectraux)}
\begin{enumerate}
\item Avec les notations du lemme \ref{lemSESP},  $X$ est un 
\textsl{sous-espace
topologique} de $Y$. En outre $\OQC(X)=\sotq{U\cap X}{U\in\OQC(Y)}$. 
On dit que
\textsl{$X$ est un sous-espace spectral de~$Y$.}
\item Pour qu'une partie $X$ d'un espace spectral $Y$ soit un sous-espace
spectral il faut et suffit que les conditions suivantes soient 
vérifiées: \\
-- La topologie induite par $Y$  fait de $X$ un espace spectral, et\\
--  $\OQC(X)=\sotq{U\cap X}{U\in\OQC(Y)}$.
\item Une partie $X$ d'un espace spectral $Y$ est un sous-espace 
spectral \ssi
elle est fermée pour la topologie constructible.
\item Si $Z$ est une partie arbitraire d'un espace spectral 
$Y=\Spec\,\gT$ son
adhérence pour la topologie constructible est égale à 
$X=\Spec\,\gT'$ où
$\gT'$ est le treillis quotient de~$\gT$ défini par la relation de 
préordre
$\preceq$ suivante:
\begin{equation} \label{eqSSES}
a\preceq b\quad \Longleftrightarrow\quad (\DT(a)\cap Z)\subseteq 
(\DT(b)\cap Z)
\end{equation}
En outre, $X$ est le plus petit sous-espace spectral de $Y$ contenant 
$Z$.
\end{enumerate}
\end{theorem}
\begin{proof}
Le point \textsl{1} est facile, et définit la notion de sous-espace 
spectral. Le point \textsl{2} en résulte. Le point~\textsl{3} résulte  des points \textsl{2} et \textsl{4}. Montrons le 
point \textsl{4}. \\
Remarquons tout d'abord que la relation (\ref{eqSSES}) définit bien 
un treillis quotient $\gT'$ car les relations~(\ref{eqPreceq}) sont 
trivialement vérifiées si on tient compte des relations (\ref{eqDa}).\\
Montrons que $X=\Spec\,\gT'$ est le plus petit sous-espace spectral 
de $Y$ contenant $Z$.\\
Tout d'abord $Z\subseteq X$: soit $\fp\in Z$, nous voulons montrer 
que si
$b\in\fp$ et $a\preceq b$ alors $a\in\fp$. Si $\DT(a)\cap Z \subseteq 
\DT(b)$ et
$b\in\fp$ alors $\fp\notin \DT(b)$ donc  $\fp\notin \DT(a)\cap Z$ donc
$\fp\notin \DT(a)$ \cad $a\in\fp$.

\noindent Par ailleurs $X$ est minimal. En effet effet si 
$X_1=\Spec\,\gT_1$ est
un sous-espace spectral de $Y$ contenant~$Z$, on a
$a\leq_{\gT_1}b\;\Leftrightarrow\;(\DT(a)\cap X_1)\subseteq 
(\DT(b)\cap X_1)$ ce
qui implique $(\DT(a)\cap Z)\subseteq (\DT(b)\cap Z)$ et donc   
$a\leq_{\gT'}b$, 
d'où $X\subseteq X_1$.\\
  Il reste à montrer que $X$ est l'adhérence de $Z$ pour la 
topologie
constructible. Notons $\wi Z$ cette adhérence.
Nous voulons donc démontrer pour tout $\fp\in\Spec\,\gT$ 
l'équivalence des
deux propriétés suivantes:\\
(1) $\fp\in\wi Z$, \cad:  $\Tt a,b\in\gT ,\;(\fp\in \DT(a)\cap
\VT(b)\,\Rightarrow\, \DT(a)\cap \VT(b)\cap Z\neq \emptyset) $,\\
(2)  $\fp\in\Spec\,\gT'$.\\
Or (2) équivaut successivement à 

\vspace{-1.3em}
$$\begin{array}{lcr}
\Tt a,b\in\gT\;((a\preceq b,\,b\in\fp)\Rightarrow a\in\fp) &\quad    
& (3)
\\[1mm]
\Tt a,b\in\gT\;(a\preceq b,\,b\in\fp,\, a\notin\fp) \;\mathrm{sont\;
incompatibles} &\quad    & (4)  \\[1mm]
\Tt a,b\in\gT\;\;\;\DT(a)\cap Z\subseteq  \DT(b) 
\;\,\mathrm{et}\,\;\fp\in
\DT(a)\cap \VT(b) \,\;\mathrm{sont\; incompatibles} &\quad    & (5) 
\\[1mm]
\Tt a,b\in\gT\;\;\;\DT(a)\cap \VT(b)\cap Z=\emptyset 
\;\,\mathrm{et}\,\;\fp\in
\DT(a)\cap \VT(b) \,\;\mathrm{sont\; incompatibles} &\quad    & (6)
\end{array}$$

\vspace{-.4em}
\noindent et (6)  est clairement équivalent à (1).
\end{proof}

\begin{corollary}
\label{corpropSESP}
Toute réunion finie et toute intersection de sous-espaces spectraux 
de
$X=\Spec\,\gT$ est un sous-espace spectral.
\begin{itemize}
\item Si $X_i=\Spec\,\gT_i$ pour un quotient $\pi_i:\gT\to\gT_i$ alors 
$\bigcap_iX_i$
correspond au quotient engendré par toutes les relations 
$\pi_i(x)=\pi_i(y)$.
\item  Si la famille est finie alors $\bigcup_iX_i$ correspond au quotient par la 
relation
$\&_i(\pi_i(x)=\pi_i(y))$.
%
\end{itemize}

\end{corollary}

\begin{proposition}
\label{propositionOFBSES} \emph{(Ouverts et fermés de base)} Soit $\gT$ un \trdi et $X=\Spec\,\gT$.
\begin{enumerate}
\item $\DT(a)$ est un sous-espace spectral de $X$ canoniquement 
homéomorphe
à $
\Spec(\gT/(a=1))$.
\item  $\VT(b)$ est un sous-espace spectral de $X$ canoniquement 
homéomorphe
à $
\Spec(\gT/(b=0))$.
\end{enumerate}
\end{proposition}
\begin{proof}
Soit $x\preceq y$ l'ordre partiel correspondant au sous-espace 
spectral
$\DT(a)$. On a donc:
$$x\preceq y\;\Leftrightarrow\;\DT(x) \cap \DT(a)\subseteq \DT(y) \cap
\DT(a)\;\Leftrightarrow\;\DT(x\vi a)\subseteq \DT(y\vi 
a)\;\Leftrightarrow\;x\vi
a\leq y\vi a
$$
et ceci est bien la relation de préordre correspondant au quotient
$\Spec(\gT/(a=1)).$\\
Soit maintenant $x\preceq' y$ l'ordre partiel correspondant au 
sous-espace spectral $\VT(b)$. On a:
$$\begin{array}{rcccl}
x\preceq' y& \Longleftrightarrow  &\DT(x) \cap \VT(b)\subseteq \DT(y) 
\cap
\VT(b)   & \Longleftrightarrow   &  \DT(x) \cup \DT(b)\subseteq 
\DT(y) \cup
\DT(b) \\[1mm]
& \Longleftrightarrow  &  \DT(x\vu b)\subseteq \DT(y\vu b) & 
\Longleftrightarrow
&   x\vu b\leq y\vu b
\end{array}$$
et ceci est bien la relation de préordre correspondant au quotient
$\Spec(\gT/(b=0)).$
\end{proof}
\subsubsection*{Fermés de $\Spec\,\gT$}
\addcontentsline{toc}{subsubsection}{Fermés de $\Spec\,\gT$}

Dans ce paragraphe $\gT$ est un \trdi fixé et $X=\Spec\,\gT$.
Si $Z\subseteq X$ on notera $\ov{Z}$ l'ahérence de $Z$ pour la 
topologie
usuelle de $X$.

\begin{proposition}
\label{propositionFSES} \emph{(Fermés de $\Spec\,\gT$)}
\begin{enumerate}
\item Un fermé arbitraire de $\Spec\,\gT$ est de la forme
$\VT(\fJ)=\bigcap_{x\in \fJ}\VT(x)$ où $\fJ$ est un \id arbitraire 
de $\gT$.
C'est un sous-espace spectral et il correspond au quotient 
$\gT/(\fJ=0)$.
\item
L'intersection d'une famille de fermés correspond au sup de la 
famille
d'idéaux. La réunion de deux fermés correspond à 
l'intersection des deux
idéaux.
\item\label{enumtra}
Le treillis $\gT/((a:b)=0)$ est le quotient correspondant à 
$\ov{\VT(a)\cap
\DT(b)}$.
\item
  Donc $\gT$ est une \agH \ssi $X$ vérifie la propriété 
suivante: pour tous
\oqcs $U_1$ et $U_2$, l'adhérence de $U_1\setminus U_2$ est 
le complémentaire d'un 
\oqcz.
\item \`A l'adhérence de $\DT(x)$ correspond le quotient 
$\gT/((0:x)=0)$.
\item Donc à la frontière de $\DT(x)$ correspond le quotient
$\gT\ul x=\gT/(\rK_\gT^x=0)$, où
\begin{equation} \label{eqbordsup}
\rK_\gT^x\,=\,\dar x \,\vu\, (0:x)
\end{equation}
Le treillis $\gT\ul x$ sera appelé \emph{le bord supérieur (de 
Krull) de $x$
dans $\gT$}. On dira aussi que $\rK_\gT^x$ est \emph{l'\id bord de 
Krull de $x$
dans $\gT$.}\\
Lorsque $\gT$ est une \agHz, $\rK_\gT^x\,=\,\dar (x \,\vu\, \lnot x)$ 
et
$\gT\ul x\simeq \uar (x \,\vu\, \lnot x)$ avec l'\homo surjectif
$\pi\ul x:
\left|
\begin{array}{rcl}
\gT& \to  & \uar (x \,\vu\, \lnot x)  \\
y&  \mapsto  & y \,\vu\, x \,\vu\, \lnot x
\end{array}
  \right.
$.
\end{enumerate}
\end{proposition}
\begin{proof}
Pour le seul point délicat (le point \textsl{\ref{enumtra}}), on dit: puisque
$(a:b)=\sotq{x}{x\vi b\leq a}$, la variété associée
$\VT(a:b)$ est l'intersection des $\VT(x)$ tels que $\VT(a)\subseteq 
\VT(x)\cup
\VT(b)$, \cade tels que  $\VT(a)\cap \DT(b)\subseteq \VT(x)$. Or tout 
fermé de
$\Spec\, \gT$ est une intersection de fermés de base $\VT(x)$, donc 
on obtient
bien l'adhérence de $\VT(a)\cap \DT(b)$.
\end{proof}

\rems

\noindent 1) On notera qu'un ouvert arbitraire de $X$ n'est pas en 
général
un sous-espace spectral.

\noindent 2) La \dfn que nous avons donnée pour le treillis bord 
$\gT\ul x$
est \cov 
car elle ne nécessite pas que l'espace $\Spec\,\gT$ ait suffisamment 
de points.
Elle est motivée par l'interprétation qu'on donne de $\Spec\,\gT$ en
\clamaz.

\medskip Le lemme suivant permet de mieux cerner l'\id 
bord de Krull de $x$ dans~$\gT$.

\begin{lemma}
\label{lemBKReg}
Pour tout $x\in\gT$ l'\id $\fj$ bord de Krull de $x$ dans $\gT$ est
\emph{régulier}, \cad \gui{son annulateur est réduit à $0$}. I.e.
$0:\fj=0$.
\end{lemma}
\begin{proof}
Soit $u\in(0:\fj)$. Puisque $\fj=\dar x\vu (0:x)$ on a $u\vi x=0$ et, 
pour tout
$z\in(0:x)$, $u\vi z=0$. En particulier $u\vi u=0$.
\end{proof}

Notez que dans le cas d'une \agH il ne s'agit de rien d'autre que de 
la loi
découverte par Brouwer: $\lnot(x \,\vu\, \lnot x)=0$.

\smallskip La proposition qui suit est la version duale, \covz, 
\gui{sans
points} du fait topologique suivant: si $A$ et $B$  sont fermés, la 
réunion
des bords de $A\cup B$  et de  $A\cap B$  est égale à celle des 
bords de $A$
et $B$.
\begin{proposition}
\label{propBordKUnion}
Pour tous $x,y\in\gT$ on a:
$\;\;\rK_\gT^{x}\cap \rK_\gT^{y}\;\;=\;\; \rK_\gT^{x\vu y}\cap 
\rK_\gT^{x\vi
y}.$
\end{proposition}
\begin{proof}
Soit $z\in\rK_\gT^{x}\cap \rK_\gT^{y}$, autrement dit il existe $u$ 
et $v$ tels
que $z\leq x\vu u$ et $u\vi x=0$, $z\leq y\vu v$ et $v\vi y=0$.
Alors $z\leq (x\vu (u\vu v))\vi (y\vu (u\vu v)) = (x\vi y)\vu(u\vu 
v)$ avec
$(u\vu v)\vi(x\vi y)=(u\vi(x\vi y))\vu (v\vi(x\vi y))=0\vu 0=0$ donc
$z\in\rK_\gT^{x\vi y}$. \\
De même $z\leq (x\vu y)\vu(u\vi v)$ avec $(u\vi v)\vi(x\vu y)=0$ 
donc
$z\in\rK_\gT^{x\vu y}$.\\
Enfin supposons  $z\in\rK_\gT^{x\vu y}\cap \rK_\gT^{x\vi y}$, 
autrement dit il
existe $u$ et $v$ tels que $z\leq x\vu y\vu u$ et $u\vi (x\vu y)=0$,
$z\leq (x\vi y)\vu v$ et $v\vi x\vi y=0$.
Soit $u_1=(y\vu u)\vi v$. On a  $z\leq (x\vi y)\vu v\leq x\vu v$ et 
$z\leq x\vu
(y\vu u)$ donc $z\leq x\vu u_1$. Par ailleurs $x\vi u_1=x\vi (y\vu u) 
\vi v\leq
x\vi y\vi v =0$ et donc $z\in\rK_\gT^{x}$.
\end{proof}

\subsubsection*{Fermés de $\Spec\,\gT\cir$}
\addcontentsline{toc}{subsubsection}{Fermés de $\Spec\,\gT\cir$}

\hum{Dans l'article original, on avait partout $\DT(\cdot)$ à la place de $\VTo(\cdot)$, mais cette notation semble très peu naturelle pour $\DT(F)$.}

Notons qu'il est naturel de noter $\VTo(a)$ pour $\DT(a)$.
Introduisons alors la notation suivante, pour $F\subseteq \gT$:
$\VTo(F)=\bigcap_{a\in\fF} \VTo(a)$. Si $\fF$ est le filtre engendré 
par $F$, on
a $ \VTo(F)=\VTo(\fF)$.

La proposition suivante découle de la proposition 
\ref{propositionFSES} par
renversement de l'ordre (on n'a réécrit que les points \textsl{1} et \textsl{6}) 
modulo
l'identification des ensembles sous-jacents à $\Spec\,\gT$ et
$\Spec\,\gT\cir$.

  Rappelons que la notion opposée à l'idéal $a:b$ est le filtre 
$a\setminus
b\eqdefi \sotq{z}{z\vu b\geq a}$. 

\begin{proposition}
\label{FdSES} \emph{(Fermés de $\Spec\,\gT\cir$)}
\begin{enumerate}
\item Un fermé arbitraire de $\Spec\,\gT\cir$ est de la forme
$\bigcap_{x\in\fF} \VTo(x)$ où $\fF$ est un filtre arbitraire 
de~$\gT$.
C'est le sous-espace spectral qui correspond au quotient 
$\gT/(\fF=1)$.

\item On définit le quotient
$\gT\bal x=\gT/(\rK^\gT_x=1)$, où $\rK^\gT_x$ est le filtre
\begin{equation} \label{eqbordinf}
\rK^\gT_x\,=\,\uar x \,\vi\, (1\setminus x)
\end{equation}
Le treillis $\gT\bal x$ sera appelé \emph{le bord inférieur (de 
Krull) de
$x$ dans $\gT$}.
On dira aussi que $\rK^\gT_x$ est \emph{le filtre bord de Krull de 
$x$}.\\
Lorsque $\gT$ est une \alg de Brouwer, $\rK^\gT_x\,=\,
\uar (x \,\vi\, (1- x))$ et
$\gT\bal x\simeq \dar (x \,\vi\, (1- x))$ avec l'\homo surjectif
$\pi\bal x:
\left|
\begin{array}{rcl}
\gT&\to&\dar (x \,\vi\, (1- x))\\
  y&\mapsto& y \,\vi\, x \,\vi\, (1- x)
\end{array}
  \right.
$.\\
Le quotient $\gT\bal x$ correspond à la 
frontière de $\VTo(x)$ pour la topologie de 
$\Spec\,\gT\cir$\footnote{Il s'agit de la notion \gui{opposée} à celle de frontière. L'intersection des adhérences de $\VTo(x)=\DT(x)$ et $\DTo(x)=\VT(x)$ est remplacée par la réunion de leurs intérieurs. Dans $\Spec\,\gT$, c'est donc le complémentaire de la frontière de~$\DT(x)$.}.
\end{enumerate}
\end{proposition}

\subsubsection*{Recollement d'espaces spectraux}
\addcontentsline{toc}{subsubsection}{Recollement d'espaces spectraux}

\hum{Ajout de la phrase suivante.}
Comme la théorie des \trdis est purement équationnelle, la catégorie possède des limites inductives et projectives arbitraires. Les limites projectives et les limites inductives filtrantes sont conservées par le foncteur d'oubli dans la catégorie des ensembles. Les propriétés duales sont donc satisfaites dans la catégorie antiéquivalente des espaces spectraux et morphismes spectraux. Dans certains cas
ces limites correspondent à celles obtenues dans la catégorie des espaces topologiques et applications continues.

Voici ce que donne par dualité la proposition~\ref{propRecolTD}
(on laisse \alec la traduction du fait~\ref{factRecolTD}).

\begin{proposition}
\label{propRecolSpec}\emph{(Recollement d'une famille finie d'espaces spectraux le long  d'\oqcsz)}
\begin{enumerate}
\item Soit $(X_i)_{1\leq i\leq n}$ une famille finie d'espaces 
spectraux, et
pour chaque $i\neq  j$ un \oqc $X_{ij}$ de $X_i$ avec 
un \iso
$\varphi_{ij}: X_{ij}\to X_{ji}$. On suppose que $\varphi_{ij}= \varphi_{ji}^{-1}$ pour tous $i,j$ et que les relations de compatibilité 
naturelles \gui{trois par trois}
sont vérifiées: si $x=\varphi_{ji}(y)=\varphi_{ki}(z)$ alors $\varphi_{jk}(y)=z$. Alors la limite inductive du diagramme dans la 
catégorie des espaces spectraux est un espace $X$ pour lequel chacun des $X_i$ s'identifie à 
un \oqc via le morphisme $X_i\to X$. 
\item En remplaçant \gui{ouvert quasi-compact} par \gui{fermé de base} le résultat analogue est \egmt valable.
\end{enumerate}

\hum{Dans le point 1, $X_i$ s'identifie à un \oqc de $X$ via le morphisme $X_i\to X$. Dans l'article original il y avait seulement écrit \gui{$X_i$ est un sous-espace spectral}. Le commentaire qui suit est 
en outre plus étoffé que dans l'article original.}

\end{proposition}

\comm Notez que $X$ est aussi la limite 
inductive du diagramme formé par les $X_i$, les inclusions $f_{ij}:X_{ij}\to X_i$ et les isomorphismes $\varphi_{ij}$,
dans la catégorie des ensembles et dans celle des espaces topologiques.
En langage plus imagé: l'espace topologique~$X$ s'obtient comme recollement des espaces $X_i$ le long des $X_{ij}$ en identifiant $x\in X_{ij}$ à $\varphi_{ij}(x)\in X_{ji}$.

En \clama on aurait plutôt tendance à déduire la proposition~\ref{propRecolTD}  de la proposition~\ref{propRecolSpec}
car cette dernière a une démonstration directe facile. Cependant ce raccourci élégant ne permet pas d'obtenir la \prco de la proposition~\ref{propRecolTD}. 

Dans le point \textsl{2} de la proposition~\ref{propRecolSpec}, si on prend pour $X_{ij}$ des fermés arbitraires, (qui sont des sous-espaces spectraux) au lieu de
 fermés de base, le recollement aura lieu en tant qu'espaces topologiques mais ne fournirait pas nécessairement un espace spectral.
 Dans le point \textsl{1} si on prend une infinité d'ouverts de base, le recollement aura lieu en tant qu'espaces topologiques mais ne fournirait pas nécessairement un espace spectral. \eoe

\subsection{Spectre maximal et spectre de Heitmann}

Dans l'article remarquable \cite[\textsl{Generating non-Noetherian modules efficiently}]{Hei84} Raymond Heitmann explique que la notion
usuelle de j-spectrum pour un anneau commutatif n'est pas la bonne 
dans le cas
non noethérien car elle ne correspond pas à un espace spectral au sens de
Stone. Il introduit la modification suivante de la \dfn usuelle: au lieu
de considérer l'ensemble des \ideps qui sont intersections 
d'\idemas il propose de considérer l'adhérence du spectre maximal dans le 
spectre premier, adhérence à prendre au sens de la topologie 
constructible (la patch topology).

\begin{definitions}
\label{defHspec1}
Soit $\gT$ un \trdiz.
\begin{enumerate}
\item  On note $\Max\,\gT$ le sous-espace topologique de $\Spec\,\gT$ formé par les \idemas de~$\gT$. 
On l'appelle le \textsl{spectre maximal de~$\gT$}.
\item  On note $\jspec\,\gT$ le sous-espace topologique de 
$\Spec\,\gT$ formé
par les $\fp$ qui vérifient l'égalité \hbox{$\JT(\fp)=\fp$}, \cad les  \ideps $\fp$ 
qui sont
intersections d'\idemas (c'est le j-spectrum \gui{usuel}).
\item On appelle \textsl{$\rJ$-spectre de Heitmann de $\gT$} et on note
$\Jspec\,\gT$ l'adhérence du spectre maximal dans $\Spec\,\gT$, 
adhérence à prendre au sens de la topologie constructible. Cet ensemble est muni de la topologie induite par $\Spec\,\gT$.
\item  On note $\Min\,\gT$ le sous-espace topologique de $\Spec\,\gT$ 
formé par les \ideps minimaux de $\gT$. On l'appelle le \textsl{spectre  minimal de $\gT$}.
\end{enumerate}
\end{definitions}

Notez que malgré leurs dénominations, les espaces topologiques $\Max\,\gT$,  $\jspec\,\gT$ 
et
$\Min\,\gT$ ne sont pas en général des espaces spectraux.

\begin{theorem}
\label{thDK3}
Soit $\gT$ un \trdiz. L'espace  $\Jspec\,\gT$ est un sous-espace 
spectral de
$\Spec\,\gT$ canoniquement homéomorphe à $\Spec\,\He(\gT)$. Plus
précisément, si $M=\Max\,\gT$,  on a pour
$a,b\in\gT$:
\begin{equation} \label{eqthDK3}
\DT(a)\cap M\,\subseteq\, \DT(b)\cap M \quad \Longleftrightarrow\quad
a\preceq_{\He(\gT)}b.
\end{equation}
\end{theorem}
\begin{proof}
Le treillis $\He(\gT)$ a été défini page
\pageref{defHeT}.\\
La deuxième affirmation implique la première. En effet 
$\Jspec\,\gT$,
d'après le \tho \ref{propSESP}\,(4), est le spectre du treillis 
quotient
$\gT'$ correspondant à la relation de préordre $a\leq _{\gT'}b$ 
définie
par $\DT(a)\cap M\,\subseteq\, \DT(b)\cap M$.\\
Pour le deuxième point on remarque que \propeq
$$\begin{array}{lcl}
\DT(a)\cap M\subseteq \DT(b)\cap M& \;  &  (1)  \\ [1mm]
\DT(a)\cap \VT(b)\cap M =\emptyset  &&  (2)  \\[1mm]
\Tt \fm\in M\; (b\notin\fm \mathrm{\;ou\;} a\in\fm)  &&  (3)  \\[1mm]
  \Tt \fm\in M\; (b\in\fm\Rightarrow  a\in\fm) &&  (4)
\end{array}$$
Et l'assertion (4) revient à dire que, vu dans le treillis quotient 
$\gT/(b=0)$, $a$
appartient au radical de Jacobson. Cela signifie $a\in \JT(b)$,
\cad $a\preceq_{\He(\gT)}b$.
\end{proof}

Quelques points de comparaison.

\begin{fact}
\label{factSpec=Jspec} ~
\begin{enumerate}
\item $\Spec\,\gT=\Jspec\,\gT$ \ssi $\gT=\He(\gT)$, \cad si $\gT$ est 
faiblement
Jacobson.
\item $\Max\,\gT=\Jspec\,\gT$ \ssi $\He(\gT)$ est une \agBz.
\item Si $\gT$ possède un complément de Brouwer, $\Max\,\gT$ est 
un fermé
de $\Spec\,\gT$, correspondant à l'\id $\Imax(\gT)$. C'est un 
espace de Stone,
il est égal à $\Jspec\,\gT$.
\item $\Min\,\gT=\Jspec\,\gT\cir$ \ssi $\He(\gT\cir)$ est une \agBz.
\item Si $\gT$ possède une négation, $\Min\,\gT$  est un fermé 
de
$\Spec\,\gT\cir$, correspondant au filtre $\Fmin(\gT)$. C'est un 
espace de
Stone, il est égal à $\Jspec\,\gT\cir$.
\end{enumerate}
\end{fact}
\begin{proof}
Le point \textsl{1} résulte du \tho \ref{thDK3}. Pour le point \textsl{2}, (on est en 
\clamaz)
on remarque qu'un treillis distributif est une \agB \ssi ses \ideps 
sont tous
maximaux. Le point \textsl{3} résulte du point \textsl{2} et du fait 
\ref{factSpecMax}. Les
points \textsl{4} et \textsl{5} s'obtiennent à partir des points \textsl{2} et \textsl{3} en passant au 
treillis
opposé.
\end{proof}

La proposition suivante est signalée par Heitmann dans 
\cite{Hei84}. L'hypothèse dans le point \textsl{2} est que l'espace $M=\Max\,\gT$ est noethérien, \cad que tout ouvert est \qcz. Comme la toplogie de $M$ est induite par celle de $\Spec\,\gT$, les ouverts $\fD_\gT(a)\cap M$ forment une base de la topologie. Par ailleurs on a  $\fD_\gT(a_1 \vu\dots\vu a_n)=\fD_\gT(a_1)\cup\dots\cup\fD_\gT(a_n)$. Donc, lorsque~$M$ est noethérien, tout ouvert de~$M$
est de la forme $\fD_\gT(a)\cap M$ et tout fermé est un fermé de base $\fV_\gT(a)\cap M$.
\begin{proposition}
\label{propJspecjspec} \emph{(Comparaison de $\Jspec$ et de $\jspec$)}
\begin{enumerate}
\item On a toujours $\jspec\,\gT\subseteq \Jspec\,\gT$.
\item Si $M=\Max\,\gT$ est noethérien, a fortiori si $\Spec\,\gT$ est noethérien, on a 
$\jspec\,\gT=\Jspec\,\gT$.
\end{enumerate}
\end{proposition}
\begin{proof}
On considère un \elt $\fp$ fixé de $\Spec\,\gT$.
Tout d'abord les propriétés suivantes sont successivement 
équivalentes
$$\begin{array}{rcl}
\fp\in\jspec\,\gT&   &     \\[1mm]
\Tt a\in \gT\quad [\;\fp\in \DT(a)&  \Rightarrow  &
\Ex \fm\in (M\,\cap \DT(a)),\, \fp\subseteq \fm \;]  \\[1mm]
\Tt a\in \gT\quad [\;\fp\in \DT(a)&  \Rightarrow  &
\Ex \fm\in (M\,\cap  \DT(a)),\, \Tt b\in\gT\,(\fm\in \DT(b)
\Rightarrow \fp\in \DT(b))\;] \\[1mm]
\Tt a\in \gT\quad [\;\fp\in \DT(a)&  \Rightarrow  &
\Ex \fm\in (M\,\cap \DT(a)),\, \Tt b\in\gT\,(\fp\in \VT(b)
\Rightarrow \fm\in \VT(b))\;]\\[1mm]
\Tt a\in \gT\quad [\;\fp\in \DT(a)&  \Rightarrow  &
\Ex \fm\in M,\, \Tt b\in\gT\,(\fp\in \VT(b)
\Rightarrow \fm\in  M\cap \DT(a) \cap \VT(b) )\;]\qquad (*)
\end{array}$$
De même sont équivalentes les propriétés
$$\begin{array}{rcl}
\fp\in\Jspec\,\gT&   &     \\[1mm]
\Tt a,b\in \gT\quad [\;\fp\in (\DT(a)\cap \VT(b))&  \Rightarrow  &
\Ex \fm\in M\cap \DT(a) \cap \VT(b)\;]
\qquad\qquad  (**)
\end{array}$$
\textsl{1}. Ainsi on voit que $(*)$ est plus fort que $(**)$, puisque dans $(**)$,
$\fm$ peut dépendre de $a$ et $b$ tandis que dans $(*)$ il ne doit 
dépendre
que de $a$.\\
\textsl{2}. Si $M=\Max\,\gT$ est noethérien le fermé 
$\bigcap\limits_{\VT(b)\ni\fp}
(M\cap \VT(b))$ est égal
à un fermé de base $M\cap \VT(b_0)$ et  $(**)$ avec ce $b_0$ donne $(*)$.
\end{proof}

\comm
Comme le fait remarquer Heitmann,  les \thos connus utilisant le 
$\jspec$ le
font toujours sous l'hypothèse \gui{$\Max$ noethérien}. Il est 
donc
probable que  $\Jspec$ soit la seule notion vraiment intéressante.
Notez que $\jspec\,\gT$ n'est un sous-espace spectral de $\Spec\,\gT$ 
que lorsqu'il est égal à $\Jspec\,\gT$.\eoe

\hum{

\textsl{1}. La notion opposée du $\Jspec$ n'est pas sans intérêt en 
\alg
commutative:
par renversement, le radical de Jacobson est remplacé par le filtre
complémentaire de la réunion des \ideps minimaux.


\textsl{2}. Pour le point \emph{2}, l'hypothèse fonctionnerait-elle avec $\Jspec$ noethérien ou $\Max$  noethérien~? Quel est le rapport entre les trois conditions de noethérianité ($\jspec\,\gT$, $\Jspec\,\gT$ et $\Max\,\gT$)~? En fait Heitmann affirme au début de l'article comme une évidence que l'on a $\Jspec=\jspec$ quand le spectre maximal est noethérien.
Il discute un peu plus cette question à la fin de l'article.

\textsl{3}. Il serait intéressant de formuler l'hypothèse \gui{$\Max\,\gT$ est noethérien} sans utiliser les points de $\Spec\,\gT$, \cad comme une propriété du \trdi $\gT$ formulée \cotz.  
}

\Hum {Ajout non retenu dans la version finale: cette section 2.4 

\subsection{Quelques points de l'anti\eqvc de catégories}
\label{subsecAntiEquiv}

Signalons sans démonstration les résultats généraux suivants concernant les morphismes dans la catégorie des \trdis et dans celle des espaces spectraux (voir le théorème de Krull page~\pageref{ThKrull}, \cite[Theorem~IV-2.6]{BW74} et \cite{Lom2020}). 

Le contexte est le suivant: soit $f:\gT\to\gT'$ un morphisme de treillis distributifs et $\Spec(f)$, noté $f^\star$, le morphisme dual, de $X=\Spec(\gT')$ vers $Y=\Spec(\gT)$, dans la catégorie des espaces spectraux.

Rappelons quelques \dfns usuelles en \clamaz.
\begin{itemize}
\item Le morphisme $f$ est dit \textsl{lying over} (en français, il poss\`ede la propriété de relèvement) lorsque~$f^\star$ est surjectif: tout \idep de $\gT$ est image réciproque d'un \idep de~$\gT'$.
\item  Le morphisme $f$ est dit \textsl{going up} (en français, il poss\`ede la propriété de montée pour les cha\^{\i}nes d'\idepsz) lorsque l'on a: \textsl{si $\fq\in\Spec(\gT')$, $f^\star(\fq)=\fp$, et $\fp\subseteq\fp_2$ dans $\Spec(\gT)$, il existe  $\fq_2\in\Spec(\gT')$ tel que
$\fq\subseteq\fq_2$ et $f^\star(\fq_2)=\fp_2$}.
\item  De même $f$ est dit \textsl{going down} (en français, il poss\`ede la propriété de descente pour les cha\^{\i}nes d'\idepsz) lorsque l'on a: \textsl{si $\fq\in\Spec(\gT')$, $f^\star(\fq)=\fp$, et $\fp\supseteq\fp_2$ dans $\Spec(\gT)$, il existe  $\fq_2\in\Spec(\gT')$ tel que $\fq\supseteq\fq_2$ et  $f^\star(\fq_2)=\fp_2$}.
\item  On dit que le morphisme $f$ \textsl{possède la propriété d'incomparabilité} lorsque ses  \gui{fibres} sont formées d'\ideps deux à deux incomparables: si $\fq_1\subseteq \fq_2\in X$ et $f^\star(\fq_1)=f^\star(\fq_2)$ dans~$Y$ alors $\fq_1= \fq_2$.
\item  L'espace spectral $\Spec(\gT)$ est dit \textsl{normal} si tout point est majoré par un unique point fermé (tout \idep de $\gT$ est contenu dans un unique \idemaz).
%
\end{itemize}

\begin{theorem} \label{th-dico-trdi-spec-mor}
 On a les équivalences suivantes.   
\begin{enumerate}
\item $f$ est injectif $\Leftrightarrow$ $f$ est un monomorphisme $\Leftrightarrow$ $f^\star$ est un épimorphisme $\Leftrightarrow$ $f^\star$ est surjectif (lying over).
\item $f$ est un épimorphisme $\Leftrightarrow$ $f^\star$ est un monomorphisme $\Leftrightarrow$ $f^\star$ est injectif.
\item $f$ est surjectif\footnote{Autrement dit, puisque c'est une structure équationnelle, $f$  est un morphisme de passage au quotient.} $\Leftrightarrow$ $f^\star$ est un isomorphisme sur son image, qui est un
sous-espace spectral de $Y$ (voir la section \ref{secSESP}, en particulier le lemme \ref{lemSESP}, le théorème~\ref{propSESP} et la proposition~\ref{propositionFSES}).
\item $f$ est going up $\Leftrightarrow$  pour
tous $a,c\in\gT$ et $y\in\gT'$ on a
$$
f(a)\leq f(c)\vu y \;\Rightarrow\;\exists x\in\gT\; (a\leq c \vu x \hbox{ et } f(x)\leq y).
$$ 
\item $f$ est going down $\Leftrightarrow$  pour
tous $a,c\in\gT$ et $y\in\gT'$ on a
$$
f(a)\geq f(c)\vi y \;\Rightarrow\;\exists x\in\gT\; (a\geq c \vi x \hbox{ et } f(x)\geq y).
$$
\item $f$ possède la propriété d'incomparabilité $\Leftrightarrow$ $f$ est zéro-dimensionnel\footnote{Voir ci-dessous \dots.}.   
\end{enumerate}

\end{theorem}

\smallskip Il peut cependant y avoir des épimorphismes de \trdis non surjectifs (voir \cite[\hbox{section V-8}]{BW74}).
Cela correspond à la possibilité d'un morphisme bijectif entre espaces spectraux qui ne soit pas un isomorphisme. Par exemple le morphisme spectral bijectif  $\Spec\,\gT^{\rm
bool}\to\Spec\,\gT$ n'est pas (en général) un \iso et le morphisme de treillis $\gT\to\gT^{\rm bool}$ est un \gui{épimono} qui n'est pas (en général) surjectif.
}

\section[Dimensions de Krull et de Heitmann: \trdisz]{Dimensions de 
Krull et de
Heitmann pour un \trdiz}
\label{secHtrdi}

  Nous arrivons dans cette section au coeur de l'article.
Nous reprenons le point de vue des \comaz.
Les seules preuves non \covs sont celles qui font le lien entre une 
notion
classique et sa reformulation \covz.
En \clama la \textsl{dimension de Krull} d'un \trdi est définie comme 
en \alg
commutative: c'est la borne supérieure des longueurs des 
chaînes
strictement croissantes d'idéaux premiers.

\subsection{Dimension et bords de Krull}

Nous rappelons maintenant une version constructive \elr de la 
dimension de Krull (\cite{CL2003,CLR05}) en nous appuyant sur l'intuition suivante: une  variété est de dimension $\leq k$ \ssi le bord de toute sous-variété est de dimension $\leq k-1$. Des approches \covs sensiblement équivalentes sont dans \cite{Esp78,Esp82,Esp83}.

Le \tho suivant en \clama nous donne une bonne signification 
intuitive de la
dimension de Krull dans un \trdiz.

\begin{theorem}
\label{thDK1} Soit un \trdi $\gT$ engendré par une partie $S$ et 
$\ell$
un entier positif ou nul.
Les conditions suivantes sont équivalentes.
\begin{enumerate}
\item  Le treillis $\gT$ est de dimension $\leq \ell.$
\item  Pour tout $x\in S$ le
bord  $\gT\ul{x}$ est de dimension $\leq \ell-1$.
\item  Pour tout $x\in S$ le
bord  $\gT\bal{x}$ est de dimension $\leq \ell-1$.
\item 
Pour tous $x_0,\ldots,x_\ell\in S$
il existe $a_0$,\ldots,  $a_\ell\in \gT$ vérifiant:
\begin{equation}
     a_0 \vi x_0  \leq  0\,,\;\;\;
     a_1 \vi x_1 \leq   a_0 \vu x_0\,,\;\;\; \dots\;\;,\;\;\;
     a_{\ell} \vi x_{\ell} \leq     a_{\ell-1} \vu x_{\ell-1}\,,\;\;\;
     1  \leq  a_{\ell} \vu x_{\ell}.
\end{equation}
\end{enumerate}
Lorsque $\gT$ une \agH les conditions précédentes sont 
aussi
équivalentes à
\begin{enumerate}\setcounter{enumi}{4}
\item   Pour toute
suite  $x_0,\dots,x_{\ell}$  dans  $S$  on a l'\egt
\begin{equation}
     1 = x_{\ell}\vu (x_{\ell}\im( \cdots (x_1 \vu (x_1 \im (x_0\vu
\neg x_0)))\cdots))
\end{equation}
\end{enumerate}
Lorsque $\gT$ une \alg de Brouwer les conditions 
précédentes sont
aussi
équivalentes à
\begin{enumerate}\setcounter{enumi}{5}
\item   Pour toute
suite  $x_0,\dots,x_{\ell}$  dans  $S$  on a l'\egt
\begin{equation}
     0 = x_{0}\vi (x_{0}-(x_1 \vi (x_1 - ( \cdots (x_\ell\vi
(1- x_\ell))))\cdots))
\end{equation}
\end{enumerate}
\end{theorem}

En particulier un treillis est de dimension $\leq 0$
\ssi  c'est une \agBz.

L'équivalence entre les points \textsl{1}, \textsl{4} et \textsl{5} a été établie, sans utiliser la notion de bord, dans  \cite{CL2003}, en poursuivant la problématique
d'André Joyal \cite{Joy71,Joy76} et de Luis Espa\~nol \cite{Esp78,Esp82,Esp83,Esp86,Esp88}.
Citons aussi l'article récent \cite{Esp2010} sur le sujet.

On notera aussi que la théorie des \agHs de dimension $\leq k$ est une théorie équationnelle.

\begin{Proof}{Démonstration
du \tho \ref{thDK1}}
On commence par noter que le quotient $\gT\ul{x}=\gT/\rK_\gT^x$  peut 
aussi
être vu comme l'ensemble ordonné obtenu à partir de la 
relation
de préordre $\leq ^x$ définie sur $\gT$ comme suit:
\begin{equation} \label{eqBordSup}
a\leq ^x b \qquad \Longleftrightarrow\qquad \exists y\in 
\gT\;\;(\,x\vi
y=0\;\;\& \;\;a\leq  x\vu y \vu b\,)
\end{equation}
\textsl{1} $\Leftrightarrow$ \textsl{2}:
Nous montrons tout d'abord que tout filtre maximal $\ff $ de $\gT$ 
devient
trivial dans~$\gT\ul{x}$, \cad qu'il contient $0$.
Autrement dit on doit trouver $a$ dans $\ff $ tel \hbox{que $a\leq ^x0$}.
Si $x\in \ff $ alors on prend $a=x$ et $y=0$ dans (\ref{eqBordSup}).
Si $x\notin \ff $ il existe $z\in \ff $ tel \hbox{que $x\vi z= 0$} (puisque 
le
filtre engendré par $\ff $ et $x$ est trivial dans $\gT$) et on 
prend~\hbox{$a=y=z$} dans (\ref{eqBordSup}).
Ceci montre que la dimension de  $\gT\ul{x}$ chute de au moins une
unité par rapport à celle de~$\gT$ (supposée finie).\\
Ensuite nous montrons que si on a deux filtres premiers $\ff '\subset
\ff $, $\ff $ maximal et $x\in \ff \setminus \ff '$ alors~$\ff '$ ne 
devient pas
trivial dans  $\gT\ul{x}$ (ceci montre que la dimension de  
$\gT\ul{x}$
chute de seulement une unité si $x$ est bien choisi). En effet,
dans le cas contraire, on aurait \hbox{un $z\in \ff '$} tel \hbox{que $z\vi x=0$},
mais comme $z$ et $x\in \ff $ cela ferait $0\in \ff $, ce qui est
absurde.\\
Enfin nous remarquons que si $\ff '\subset \ff $ sont des filtres 
premiers
distincts et si $S$ engendre $\gT$ on peut trouver $x\in S$ tel que
$x\in \ff \setminus \ff '$.

\noindent \textsl{1} $\Leftrightarrow$ \textsl{3}: conséquence de \textsl{1} $\Leftrightarrow$ \textsl{2} par renversement de l'ordre.

\noindent \textsl{2} $\Leftrightarrow$ \textsl{4}:  par récurrence sur $\ell$, vue 
la
\dfn du bord.

\noindent \textsl{2} $\Leftrightarrow$ \textsl{5}:  par récurrence sur $\ell$, vue 
la
\dfn du bord
dans le cas d'une \agHz.

\noindent \textsl{3} $\Leftrightarrow$ \textsl{6}:  par récurrence sur $\ell$, vue 
la
\dfn du bord
dans le cas d'une \alg de Brouwer.
\end{Proof}

On a montré (en \clamaz) que si $X$ est le spectre d'un \trdi $\gT$ 
et $x\in
\gT$, alors la \textsl{frontière} de l'\oqc $\DT(x)$ de $X$ est 
(canoniquement
isomorphe à) $\Spec(\gT\ul{x})$.
On obtient ainsi comme corolaire du \tho \ref{thDK1} le \tho suivant 
concernant
la dimension  des espaces spectraux (\cad la longueur maximale des 
chaînes
de fermés irréductibles, ou encore la longeur maximale des 
chaînes de
points pour la relation $\fp\leq_X\fq$). Rappelons que l'unique 
espace spectral
de dimension~$-1$ est le vide.

\begin{theorem}
\label{thDK2}
Soit $k$ un entier $\geq 0$. Un espace spectral $X$ est de dimension 
$\leq  k$
\ssi tout \oqc de $X$ a une frontière (un bord) de dimension~\hbox{$\leq  
k-1$}.
\end{theorem}

Concernant la dimension de Krull, on choisit en \coma la \dfn 
suivante:

\begin{definition}
\label{defDK0} \textsl{(Définition \cov de la dimension de Krull)}\\
La dimension de Krull (notée $\Kdim$) des \trdis est définie 
comme suit.
\begin{enumerate}
\item $\Kdim(\gT)=-1$ \ssi $1=_\gT0$ (i.e. le treillis est réduit 
à un
point).
\item Pour $\ell\geq 0$ on définit $\Kdim(\gT)\leq \ell$ par les 
conditions
équivalentes suivantes:
\begin{enumerate}
\item \label{l1}$\forall x\in \gT,\; \Kdim(\gT\ul x)\leq \ell-1$
\item \label{l2}$\forall x\in \gT,\;\Kdim(\gT\bal{x})\leq \ell-1$
\item \label{l3}$\forall x_0,\ldots,x_\ell\in \gT$
$\Ex a_0,\ldots,  a_\ell\in \gT$ vérifiant:
$$    a_0 \vi x_0  \leq  0\,,\;\;\;
     a_1 \vi x_1 \leq   a_0 \vu x_0\,,\dots\,,\;\;\;
     a_{\ell} \vi x_{\ell} \leq     a_{\ell-1} \vu x_{\ell-1}\,,\;\;\;
     1  \leq  a_{\ell} \vu x_{\ell}.
$$
\end{enumerate}
\end{enumerate}
\end{definition}

Notez qu'il y a en fait trois \dfns possibles ci-dessus pour 
$\Kdim(\gT)\leq
\ell$. Les \dfns basées sur \textsl{\ref{l1}} et \textsl{\ref{l2}} sont inductives, 
tandis que
la \dfn basée sur \textsl{\ref{l3}} est globale. L'équivalence des \dfns 
basées sur
\textsl{\ref{l1}} et \textsl{\ref{l3}} est immédiate par induction (même chose 
pour \textsl{\ref{l2}}
et \textsl{\ref{l3}}).

Par exemple, pour $\ell=2$ les inégalités dans le point \textsl{2c} correspondent au dessin suivant dans~$\gT$.
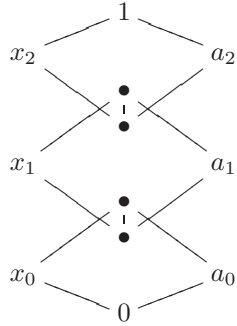
$$\SCO{x_0}{x_1}{x_2}{a_0}{a_1}{a_2}$$

\medskip \rem On peut illustrer le point \textsl{2c} dans la \dfn ci-dessus.
Nous introduisons \gui{l'idéal bord de Krull itéré}.
Pour $x_1,\ldots ,x_n\in\gT$ nous notons 
$$
\gT_\rK[x_1]=\gT\ul{x_1},\,\gT_\rK[x_1,x_2]=(\gT\ul{x_1})\ul{x_2}
,\,\gT_\rK[x_1,x_2,x_3]=((\gT\ul{x_1})\ul{x_2})\ul{x_3}, 
\hbox{ etc}\ldots\,\,
$$  
les treillis bords quotients successifs, et $\rK[\gT;x_1,\ldots
,x_k]=\rK_\gT[x_1,\ldots ,x_k]$ désigne le noyau de la projection 
canonique
$\gT\to \gT_\rK[x_1,\ldots ,x_k]$.  Alors on a 
$y\in\rK_\gT[x_0,\ldots ,x_\ell]$
\ssi $\Ex a_0,
\ldots  a_\ell\in \gT$ vérifiant:
$$    a_0 \vi x_0  \leq  0\,,\;\;\;
     a_1 \vi x_1 \leq   a_0 \vu x_0\,,\dots\,,\;\;\;
     a_{\ell} \vi x_{\ell} \leq     a_{\ell-1} \vu x_{\ell-1}\,,\;\;\;
     y  \leq  a_{\ell} \vu x_{\ell}.
$$
Si $\gT$ est une \agH on a:
$$
\rK_\gT[x_0,\ldots ,x_\ell]=\dar(x_{\ell}\vu (x_{\ell}\im( \cdots (x_1 
\vu (x_1 \im (x_0\vu
\neg x_0)))\cdots)))
$$
La dimension de Krull du treillis est $\leq \ell$ \ssi 
$1\in\rK_\gT[x_0,\ldots ,x_\ell]$
pour tous $x_0,\ldots ,x_\ell$.

\medskip En \coma la dimension de Krull de $\gT$ n'est pas à priori 
un \elt bien
défini de $\NN\cup\so{-1}\cup\so{\infty}$. En \clama cet \elt est défini 
comme la borne
inférieure des entiers $\ell$ tels que  $\Kdim(\gT)\leq \ell$.
On utilise en \coma les \textsl{notations} suivantes{\footnote{~En fait 
si on
regarde $\Kdim(\gT)$ comme l'ensemble des $\ell$ pour lesquels  
$\Kdim(\gT)\leq
\ell$, on raisonne avec des parties finales (éventuellement vides) 
de
$\NN\cup\so{-1}$, la relation d'ordre est alors l'inclusion 
renversée, la
borne supérieure l'intersection et la borne inférieure la 
réunion.}}, pour
se rapprocher du langage classique:

\begin{notation}
\label{notaKdiminf}
{\rm  Soient $\gT$, $\gL$, $\gT_i$ des \trdisz.
\begin{itemize}
\item  $\Kdim\,\gL\leq \Kdim\,\gT$ signifie: $\Tt\ell\geq -1\; 
(\Kdim\,\gT\leq
\ell\;\Rightarrow \Kdim\,\gL\leq \ell)$.
\item  $\Kdim\,\gL= \Kdim\,\gT$ signifie:  $\Kdim\,\gL\leq  
\Kdim\,\gT$ et
$\Kdim\,\gL\geq  \Kdim\,\gT$.
\item  $\Kdim\,\gT\leq  \sup_i\Kdim\,\gT_i$ signifie: $\Tt\ell\geq 
-1\; (\&_i
(\Kdim\,\gT_i\leq \ell)\;\Rightarrow\Kdim\,\gT\leq \ell) $.
\item  $\Kdim\,\gT=  \sup_i\Kdim\,\gT_i$ signifie: $\Tt\ell\geq -1\; 
(\&_i
(\Kdim\,\gT_i\leq \ell)\;\Leftrightarrow\Kdim\,\gT\leq \ell) $.
\end{itemize}
}
\end{notation}

Le fait que la \dfn fonctionne aussi bien avec le bord supérieur 
qu'avec le
bord inférieur donne la constatation suivante.
\begin{fact}
\label{corTTO}
Un \trdi et le treillis opposé ont même dimension.
\end{fact}

En \clama on peut s'en rendre compte directement en considérant les
chaînes d'\ideps qui ont pour complémentaires des chaînes 
de filtres
premiers (et vice-versa): si on identifie les ensembles sous-jacents 
à
$X=\Spec\,\gT$ et  $X'=\Spec\,\gT\cir$ les relations d'ordre $\leq_X$ 
et
$\leq_{X'}$ sont opposées.

\smallskip Notons $\Bd(V,X)$ le bord de $V$ dans $X$ ($X$ est un 
espace
topologique et  $V$  est une partie de $X$). 
Alors si $Y$ est un sous-espace de
$X$ on a $\Bd(V\cap Y,Y)\subseteq \Bd(V,X)\cap Y$, avec égalité 
si $Y$ est
un ouvert. La proposition suivante donne une version duale, \covz, 
sans points,
de cette affirmation.
\begin{proposition}
\label{propTquoBord} \emph{(Bord de Krull d'un treillis quotient)}\\
Soit $\gL$ un treillis quotient d'un \trdi $\gT$. Par abus, nous 
notons $x$
l'image de $x\in\gT$ dans $\gL$. Alors  $\gL\ul x$ est un quotient 
de  $\gT\ul
x$ et  $\gL\bal x$ est un quotient de  $\gT\bal x$. En outre si $\gL$ 
est le
quotient de $\gT$ par un filtre $\ff$,   $\gL\ul x$ est le quotient 
de  $\gT\ul
x$  par le filtre image de $\ff$ dans  $\gT\ul x$.
\end{proposition}
\begin{proof}
Soient $a,b,x\in\gT$, si $a\leq_{\gT\ul x}b$ il existe $z\in\gT$ tel 
que
$x\vi z\leq_{\gT} 0$ et $a\leq _{\gT}x\vu z\vu b$. Puisque  $\gL$ est 
un
quotient de $\gT$, on a fortiori $x\vi z\leq_{\gL}0$ et 
$a\leq_{\gL}x\vu z\vu b$
et donc
  $a\leq_{\gL\ul x}b$.
Voyons le deuxième point. Notons  $\pi:\gT\to\gL$, $\pi\ul 
x:\gT\to\gT\ul x$
et $\theta:\gT\ul x\to\gL\ul x$ les projections. Il est clair que 
$\theta(\pi\ul
x(\ff))=\so{1}$  de sorte l'on a une factorisation de $\theta$ via 
$\gT\ul
x/(\pi\ul x(\ff)=1)$. Inversement soient $a,b\in\gT$
tels que $a\leq _{\gL\ul x} b$. Nous voulons montrer que $a\leq 
_{\gT\ul x
/(\pi\ul x(\ff)=1)} b$. Par hypothèse il existe $z\in\gT$ tel que
$a \leq_{\gL}x\vu z\vu b$ et $x\vi z\leq_{\gL}0$. Cela signifie qu'il 
existe
$f_1$ et $f_2\in\ff$  tels que $a\vi f_1\leq_{\gT} b\vu  x\vu z$ et 
$x\vi z\vi
f_2\leq_{\gT} 0$.
En prenant $f=f_1\vi f_2$ et $z'=z\vi f_2$ on obtient $a\vi 
f\leq_{\gT} b\vu
x\vu z'$ et $x\vi z'\leq_{\gT}0$, \cad  $a\vi f\leq_{\gT\ul x}  b$.
\end{proof}

Le corolaire suivant donne une version duale, \covz, sans points, du 
fait
suivant: la dimension d'un sous-espace spectral est toujours 
inférieure ou
égale à celle de l'espace entier.
\begin{corollary}
\label{corpropTquoBord}
Si $\gL$ est un treillis quotient de $\gT$ on a $\Kdim\,\gL\leq 
\Kdim\,\gT$.
\end{corollary}

Dans la proposition suivante le point \textsl{2} est une version duale, \covz, 
sans
points, du fait topologique suivant: la notion de bord est locale. 
C'est surtout
le point \textsl{1} qui nous sera utile dans la suite. 

Par ailleurs, en renversant la relation d'ordre on 
aurait un
énoncé analogue pour l'autre bord.

\begin{proposition}
\label{propLocBord} \emph{(caractère local du bord de Krull)}
\begin{enumerate}
\item Soit $(\fa_i)_{1\leq i\leq m}$ une famille finie d'\ids de 
$\gT$, avec
$\bigcap_{i=1}^m\fa_i=\so{0}$.
Pour $x\in\gT$ notons encore $x$ son image dans $\gT_i=\gT/(\fa_i=0)$.
Le bord ${\gT\!_i}\ul x$ peut être vu comme le quotient de  $\gT 
\ul x$
par un \id $\fb_i$ et on a: $\bigcap_{i=1}^m\fb_i=\so{0}$.
\item Soit $(\ff_i)_{1\leq i\leq m}$ une famille finie de filtres de 
$\gT$, avec
$\bigcap_{i=1}^m\ff_i=\so{1}$.
Pour $x\in\gT$ notons encore $x$ son image dans $\gT_i=\gT/(\ff_i=1)$.
Le bord ${\gT\!_i}\ul x$ peut être vu comme le quotient de  $\gT 
\ul x$
par un filtre $\ffg_i$ et on a: $\bigcap_{i=1}^m\ffg_i=\so{1}$.
\end{enumerate}
\end{proposition}
\begin{proof}
Voyons le point \emph{1}. Considérons la projection $\pi_i:\gT\to\gT_i$  
puis la
projection $\pi'_i:\gT_i\to{\gT\!_i} \ul x$.
La composée $\gT\to{\gT\!_i} \ul x$ montre que ${\gT\!_i}\ul x\simeq
\gT/(\pi_i^{-1}(\rK_{\gT_i}^x)=0)$, et l'\id
$\pi_i^{-1}(\rK_{\gT_i}^x)$ contient
$\rK_\gT^x$. Ceci prouve la première affirmation.
Soit maintenant $y\in\gT$ tel que pour chaque $i$, 
$\pi_i(y)\in\rK_{\gT_i}^x$.
Cela revient à dire qu'il existe $b_i$ tel que
$\pi_i(y)\leq \pi_i(x)\vu\pi_i(b_i)$ et 
$\pi_i(b_i)\vi\pi_i(x)=\pi_i(0)$,
\cad pour un certain $a_i\in\fa_i$: $y\leq x\vu a_i\vu b_i$ et $x\vi 
b_i\in
\fa_i$. En prenant $c_i=a_i\vu b_i$ cela fait $y\leq x\vu c_i$ et 
$x\vi
c_i\in\fa_i$. Enfin avec $c=c_1\vi\cdots \vi c_m$ on obtient $y\leq 
x\vu c$ avec
$c\vi x=0$. Donc $y\in\rK_\gT^x$, et cela montre qu'un $z\in\gT \ul x$
qui est dans tous les $\fb_i$ est nul (car il est la classe 
d'un~$y$).\\
Pour le point \textsl{2} c'est une conséquence immédiate de la proposition
\ref{propTquoBord} et du fait \ref{factIdDansQuo}, qui affirme que 
tout passage
au quotient est un homomorphisme pour les treillis des filtres (en 
particulier
si une intersection finie de filtres est égale à $1$ cela reste 
vrai après
passage au quotient).
\end{proof}

\begin{corollary}
\label{corpropLocBord}
\emph{(Caractère local de la dimension de Krull)}
\begin{enumerate}
\item Soit $(\fa_i)_{1\leq i\leq m}$ une famille finie d'\ids de 
$\gT$ et
$\fa=\bigcap_{i=1}^m\fa_i$. \\ Alors
$\Kdim(\gT/(\fa=0))=\sup_i\Kdim(\gT/(\fa_i=0))$.
\item Soit $(\ff_i)_{1\leq i\leq m}$ une famille finie de filtres de 
$\gT$ et
$\ff=\bigcap_{i=1}^m\ff_i$. \\ Alors
$\Kdim(\gT/(\ff=1))=\sup_i\Kdim(\gT/(\ff_i=1))$.
\end{enumerate}
\end{corollary}
\begin{proof}
Il suffit de prouver le point \textsl{1}. En remplaçant $\gT$ par 
$\gT/(\fa=0)$ on se
ramène au cas \hbox{où $\fa=0$}. Le résultat est clair pour $\ell=-1$. 
Et la
preuve par induction fonctionne gr\^{a}ce à la
proposition~\ref{propLocBord}.\\
On peut aussi donner une preuve directe basée sur la 
caractérisation 2~(c)
dans la \dfnz~\ref{defDK0}.
\end{proof}

Notez qu'en \clama la caractère local de la dimension de Krull est
en général énoncé sous la forme
$$
\Kdim(\gT)=\sup\sotq{\Kdim(\gT/(\ff=1))}{\ff\hbox{ filtre premier 
minimal}}.
$$  
Il s'agit d'une conséquence directe de la \dfn de la
dimension en \clamaz.  On peut ensuite en déduire le corolaire
\ref{corpropLocBord} mais la preuve qu'on obtient ainsi n'est pas
\covz.

\smallskip On obtient \egmt comme dans le \tho \ref{thDK1}, mais avec
une \prco le résultat important suivant, qui nous dit qu'on peut se
limiter aux bords construits avec les \elts d'un système \gtr de
$\gT$.

\begin{proposition}
\label{propDK1}
Soit un \trdi $\gT$ engendré par une partie $S$ et $\ell$
un entier positif ou nul.
Les conditions suivantes sont équivalentes.
\begin{enumerate}
\item  Le treillis $\gT$ est de dimension $\leq \ell$
\item  Pour tout $x\in S$ le
bord  $\gT\ul{x}$ est de dimension $\leq \ell-1$.
\item  Pour tout $x\in S$ le
bord  $\gT\bal{x}$ est de dimension $\leq \ell-1$.
\item 
Pour tous $x_0,\ldots,x_\ell\in S$
il existe $a_0$,\ldots,  $a_\ell\in \gT$ vérifiant:
\begin{equation}
     a_0\vi x_0  \leq  0\,,\;\;\;
     a_1\vi x_1 \leq   a_0\vu x_0\,,\dots\,,\;\;\;
     a_{\ell}\vi x_{\ell} \leq     a_{\ell-1}\vu x_{\ell-1}\,,\;\;\;
     1  \leq  a_{\ell}\vu x_{\ell}.
\end{equation}
\end{enumerate}
Lorsque $\gT$ est une \agH les conditions précédentes 
sont aussi
équivalentes~à
\begin{enumerate}\setcounter{enumi}{4}
\item   Pour toute
suite  $x_0,\dots,x_{\ell}$  dans  $S$  on a l'\egt
\begin{equation}
     1 = x_{\ell}\vu (x_{\ell}\im( \cdots (x_1 \vu (x_1 \im (x_0\vu
\neg x_0)))\cdots))
\end{equation}
\end{enumerate}
Lorsque $\gT$ une \alg de Brouwer les conditions 
précédentes sont
aussi
équivalentes à
\begin{enumerate}\setcounter{enumi}{5}
\item   Pour toute
suite  $x_0,\dots,x_{\ell}$  dans  $S$  on a l'\egt
\begin{equation}
     0 = x_{0}\vi (x_{0}-(x_1 \vi (x_1 - ( \cdots (x_\ell\vi
(1- x_\ell))))\cdots))
\end{equation}
\end{enumerate}
\end{proposition}
\begin{proof}
\noindent $(2) \Leftrightarrow (4)$  par récurrence sur $\ell$, vue 
la
\dfn du bord.

\noindent \textsl{3} $\Leftrightarrow$ \textsl{4}:  par récurrence sur $\ell$, vue 
la
\dfn du bord.

\noindent \textsl{2} $\Leftrightarrow$ \textsl{5}:  par récurrence sur $\ell$, vue 
la
\dfn du bord
dans le cas d'une \agHz.

\noindent \textsl{3} $\Leftrightarrow$ \textsl{6}:  par récurrence sur $\ell$, vue 
la
\dfn du bord
dans le cas d'une \alg de Brouwer.

\noindent Il reste donc à voir que si \textsl{2} est vrai pour  $x\in S$, 
alors  \textsl{2}
est vrai pour tout $x\in\gT$.
Cela résulte de la proposition \ref{propBordKUnion}, et des 
corolaires
\ref{corpropTquoBord} et  \ref{corpropLocBord}: par exemple pour tous 
$x,y\in
\gT$, $\gT\ul{x\vu y}$ est un quotient
de $\gT/(\rK_\gT^x\cap\rK_\gT^y)$ donc $\Kdim(\gT\ul{x\vu y})\leq \sup
(\Kdim\,\gT\ul x, \Kdim\,\gT\ul y)$.
\end{proof}

\subsection{Dimensions et bord de Heitmann}

\subsubsection*{J-dimension de Heitmann pour un \trdiz}
\addcontentsline{toc}{subsubsection}{$\rJ$-dimension de Heitmann}

Nous donnons maintenant la \dfn \cov de la $\Jdim$
de Heitmann, que nous appelons \textsl{$\rJ$-dimension de Heitmann} du 
treillis
$\gT$ (ou de l'espace spectral $\Spec\,\gT$).

\begin{definition}
\label{defHdimTr}
Soit $\gT$ un \trdiz.
La \textsl{$\rJ$-dimension de Heitmann de $\gT$}, notée $\Jdim\,\gT$, est la 
dimension de Krull du treillis de Heitmann $\He(\gT)$ (cf. \dfn\ref{defHeT}).
\end{definition}

En fait d'un point de vue \cof on se contente de définir, pour tout
entier $\ell\geq -1$, la phrase   \gui{$\Jdim\,\gT\leq \ell$} par
\gui{$\Kdim\,\He(\gT)\leq \ell$}.

Du point de vue classique on peut donner la \dfn directement, à la 
Heitmann,
pour un espace spectral $X$ comme suit: si $M_X$ est l'ensemble des 
points
fermés et $J_X$ l'adhérence de $M_X$ pour la topologie 
constructible, alors
$\Jdim\,X=\Kdim\,J_X$.

\begin{fact}
\label{factKJdim} Soit $\gT$ un \trdiz, $\gT'=\gT/(\JT(0)=0)$ et 
$\fa$ un
idéal de $\gT$.
\begin{enumerate}
\item $\Jdim(\He(\gT))\;=\;\Jdim\,\gT'\;=\;\Jdim\,\gT\;\leq\; 
\Kdim(\gT')
\;\leq\; \Kdim\,\gT$.
\item Si $\gL=\gT/(\fa=0)$ est le quotient de $\gT$ par \textsl{l'\idz} 
$\fa$,
alors $\Jdim\,\gL\leq\Jdim\,\gT$.
\end{enumerate}
\end{fact}
\begin{proof}
Le point \textsl{1} est conséquence du point \textsl{3} dans le fait~\ref{factHeHe} 
et le point \textsl{2} conséquence du point~\textsl{4}.
\end{proof}

\rem L'opération $\gT\mapsto \JT(0)$ n'est pas fonctorielle. Même 
chose pour
$\gT\mapsto \He\,\gT$. En particulier l'item \textsl{2} dans le fait \ref{factKJdim} ne fonctionne plus à priori pour
un quotient plus général, par exemple pour un quotient par un 
filtre.
Contrairement à la $\Kdim$, la $\Jdim$ peut augmenter par passage 
à un quotient (on a des exemples simples en \alg commutative).

\medskip
\exl
Voici un exemple d\^{u} à Heitmann d'un espace spectral pour lequel
$\Jdim(\gT)<\Kdim(\gT/(\JT(0)=0))$. 
\begin{figure}[ht]   
\begin{center}
\includegraphics{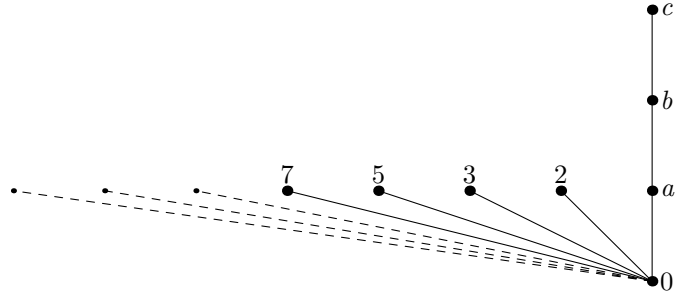}
\end{center}
\caption[]
{\label{figSpec} Exemple de Heitmann }  
\end{figure}  
On considère $X=\Spec\,\ZZ$ et $Y=\bf n$
($n\geq 3$) on recolle ces deux espaces spectraux en identifiant les 
deux \elts
minimaux (le singleton correspondant est bien dans 
les deux cas un sous-espace
spectral). On obtient un espace $Z=\Spec\,\gT$ avec $M_Z=\Max\,\gT$ 
fermé donc
égal à $J_Z=\Jspec\,\gT$ et zéro dimensionnel, tandis que 
l'unique
élément minimal est le seul minorant de $M_Z$, donc $\JT(0)=0$ et
$\Kdim(\gT/(\JT(0)=0))=\Kdim(\gT)=n-2$.\eoe

\medskip
\rem
Donnons la \dfn de la $\Jdim$ complètement mise à plat.

\begin{itemize}
\item $\Jdim\,\gT\leq \ell$ signifie:
$\forall x_0,\ldots,x_\ell\in \gT\;$
$\Ex a_0,
\ldots,  a_\ell\in \gT$ vérifiant:
$$    a_0\vi x_0  \leq_{\He(\gT)} 0\,,\;\;\;
     a_1\vi x_1 \leq_{\He(\gT)}  a_0\vu  x_0\,,\dots\,,\;\;\;
     a_{\ell}\vi x_{\ell} \leq_{\He(\gT)}    a_{\ell-1}\vu   
x_{\ell-1}\,,\;\;\;
     1  \leq_{\He(\gT)} a_{\ell}\vu   x_{\ell}.
$$
\cade
$\forall x_0,\ldots,x_\ell\in \gT\;$
$\Ex a_0,\ldots,  a_\ell\in \gT\;$ $\Tt y\in\gT$
$$\begin{array}{rcl}
(a_0 \vi x_0) \vu y = 1 &  \Rightarrow  &  y=1    \\
(a_1 \vi x_1)   \vu y=1&  \Rightarrow &  a_0 \vu x_0 \vu y=1   \\
\vdots \qquad  &  \vdots &   \qquad \vdots   \\
(a_{\ell}\vi x_{\ell})\vu y=1&  \Rightarrow &  a_{\ell-1}\vu 
x_{\ell-1} \vu y=1
\\
   &   &  a_{\ell}\vu x_{\ell}=1
  \end{array}$$
\item En particulier $\Jdim\,\gT\leq 0$ signifie:\\
  $\forall x_0\in \gT$ $\Ex a_0\in \gT$ $\Tt y\in\gT, 
  (\,(a_0 \vi x_0) \vu y = 1
\Rightarrow    y=1)\; $ et
$\; a_0 \vu x_0=1\,)$.
\item Et $\Jdim\,\gT\leq 1$ signifie:
  $\forall x_0,x_1\in \gT\;$ $\Ex a_0,a_1\in \gT\;$ $\Tt y\in\gT$\,,\\
~~~~~~  $((a_0 \vi x_0) \vu y = 1   \Rightarrow    y=1)\,$ et
$((a_1 \vi x_1) \vu y = 1   \Rightarrow a_0 \vu x_0 \vu   y=1)\,$ et
$ \,a_1 \vu x_1=1$.
\end{itemize}

\subsubsection*{La dimension de Heitmann pour un \trdiz}
\addcontentsline{toc}{subsubsection}{Dimension de Heitmann}

Bien que Heitmann définisse et utilise dans \cite{Hei84} la 
dimension
$\Jdim\,X$  où $X$ est le spectre d'un anneau commutatif, ses 
preuves sont en
fait implicitement basées sur une notion voisine, mais non 
équivalente, que
nous appellerons la \textsl{dimension de Heitmann} et que nous noterons 
$\Hdim$.

\smallskip Nous présentons cette notion directement au niveau des 
\trdisz,
où les choses s'expliquent plus simplement.

La dimension $\Hdim\,\gT$ est toujours inférieure ou égale à
$\Jdim\,\gT$, ce qui fait que les \thos établis pour la $\Hdim$ 
seront a
fortiori vrais avec la $\Jdim$ et avec la $\Kdim$.

\begin{definition}
\label{defBHeit}
Soit $\gT$ un \trdi et $x\in\gT$.
On appelle \textsl{bord de Heitmann de $x$ dans~$\gT$} le treillis quotient
$\gT_\rH^x\eqdefi\gT/(\rH_\gT^x=0)$, où
\begin{equation} \label{eqbordHeit}
\rH_\gT^x\,=\,\dar x \,\vu\, (\JT(0):x)
\end{equation}
On dira aussi que $\rH_\gT^x$ est \textsl{l'\id bord de Heitmann de $x$ 
dans
$\gT$.}
\end{definition}

\begin{lemma}
\label{lemTBKBH} \emph{(Bord de Krull et bord de Heitmann)}\\
Soit $\gT$ un \trdiz, $\gT'=\gT/(\JT(0)=0)$, $\pi:\gT\to\gT'$ la 
projection
canonique,
$x\in\gT$ et $\ov{x}=\pi(x)$. Alors on a:
\begin{enumerate}
\item $\rH^{\ov{x}}_{\gT'} =\rK^{\ov{x}}_{\gT'}.$
\item  $\pi^{-1}(\rK^{\ov{x}}_{\gT'}) =\rH^x_\gT$
et ${{\gT'}_\rH^{\ov{x}}}\simeq {{\gT'}\ul{\ov{x}}}\simeq\gT_\rH^x$.
\end{enumerate}
\end{lemma}
\begin{proof}
Clair d'après les \dfnsz.
\end{proof}

\begin{lemma}
\label{lemTquoBordH}
Soit $\gL=\gT/(\fa=0)$ un treillis quotient d'un \trdi $\gT$ par un 
\idz. Par
abus, nous notons $x$ l'image de $x\in\gT$ dans $\gL$. Alors 
$\gL_\rH^x$ est un
quotient de  $\gT_\rH^x.$
\end{lemma}
\begin{proof}
Notons $\pi:\gT\to\gL$ la projection canonique. On veut montrer que si
$z\in\rH_\gT^x$ alors $\pi(z)\in\rH_\gL^{\pi(x)}$. On suppose donc 
$z\leq_\gT
x\vu u$ avec $x \vi u  \in\JT(0)$. Comme 
$\pi(\JT(0))\subseteq\rJ_\gL(0)$, cela
donne $\pi(z)\leq_\gL \pi(x)\vu \pi(u)$ avec $\pi(x) \vi \pi(u) 
\in\rJ_\gL(0)$,
ce qui implique  $\pi(z)\in\rH_\gL^{\pi(x)}$.
\end{proof}

\rem Le lemme précédent serait faux pour un quotient plus 
général, par
exemple pour un quotient par un filtre. Il reste vrai chaque fois que
$\pi(\JT(0))\subseteq \rJ_\gL(0)$.

\begin{proposition}
\label{propBHeitHeyt} \emph{(Comparaison de deux bords à la 
Heitmann)}\\
On considère pour $x\in\gT$ son image $\hat{x}\in\He(\gT)$.
\begin{enumerate}
\item On compare les deux quotients de $\gT$ que sont 
$\He(\gT_\rH^x)$ et
$(\He(\gT))\ul{\hat{x}}$: le premier est un quotient du second.
\item Lorsque $\He(\gT)$ est une \agHz, on a \egtz.
\end{enumerate}
\end{proposition}
\begin{proof}
Le treillis $(\He(\gT))\ul{\hat{x}}$  est un quotient de $\gT$  dont 
la relation
de préordre $a \preceq  b$
peut etre décrite de la manière suivante:
$$\Ex y\;\;\;(x \vi  y \in \JT (0) \;\;\mathrm{et} \;\;\Tt z\;[\;a 
\vu  z = 1
\;\Rightarrow\; x \vu  y \vu  b  \vu  z = 1\;])\;\;\;(*)$$
Considérons
  le préordre qui définit $\He(\gT _\rH^x)$:
$$\Tt u\;\;\;(1 \leq_{\gT _\rH^x}  a \vu  u \;\;\Rightarrow\;\;
  1 \leq_{\gT _\rH^x}  b \vu  u)\;\;\;(**) $$
Prouvons que la relation de préordre $\preceq$ entraine la relation 
de
préordre $(**)$.\\
On a un $y$ vérifiant $(*)$.
On considère un $u$ tel que $1 \leq  a \vu  u$ dans $\gT _\rH^x$ et 
on cherche
à montrer que
$1 \leq  b \vu  u$  dans $\gT _\rH^x$.\\
La relation $1 \leq  a \vu  u$ dans $\gT _\rH^x$ s'ecrit $1 \leq  a 
\vu  u \vu
x \vu  y'$
pour un $y'$ tel que $y' \vi  x \in \JT (0)$.\\
On pose $z = u \vu  x \vu  y'$, on applique $(*) $ et on obtient
   $x \vu  b \vu  u \vu  (y \vu  y') = 1.$ \\
Mais $(y \vu  y') \vi  x  \in \JT (0)$  donc en posant $y'' = y \vu  
y'$
on a
$1 \leq  b \vu  u \vu  x \vu  y''$    avec $x \vi  y'' \in \JT (0)$
\cad  $1 \leq  b \vu  u$  dans $\gT _\rH^x$.\\
Voyons le point \textsl{2}. Nous notons $\pi:\gT\to\He(\gT)$ la projection 
canonique.\\
Rappelons (fait~\ref{factHeHe}) que $\pi^{-1}(0)=\JT(0)$ et
$\pi^{-1}(1)=\so{1}$. Nous supposons que $\He(\gT)$ est une \agHz. 
Notons
$\wi{x}$ un \elt de $\gT$ tel que $\pi(\wi{x})\,=\,\pi(x)\im 0\,$ dans
$\He(\gT)$.
Alors on peut réécrire $(*)$ sous la forme
$$
\Tt z\;[\;a \vu  z = 1  \;\Rightarrow\; x \vu  \wi x \vu  b  \vu  z =
1\;]\;\;\;(*')
$$
De même  $1 \leq_{\gT _\rH^x}  a \vu  u$, qui signifie
$$\Ex y'\;(x\vi y'\in \JT (0) \;\mathrm{et}\; 1 \leq  a \vu  u \vu  x 
\vu
y'),$$
  se réécrit $1 \leq  a \vu  u \vu  x \vu  \wi x$. En 
conséquence $(**)$ se
réécrit
$$
\Tt u\;[\;a \vu  u \vu  x \vu  \wi x= 1  \;\Rightarrow\; b \vu  u 
\vu  x \vu
\wi x\;]\;\;\;(**')
$$
et il est clair que $(*')$ et $(**')$ sont équivalents.
\end{proof}

\begin{definition}
\label{defHdim}
La \textsl{dimension de Heitmann d'un \trdi $\gT$}, notée 
$\Hdim\,\gT$, est
définie de manière inductive comme suit:
\begin{itemize}
\item $\Hdim\,\gT=-1$ \ssi $\gT=\Un$.
\item Pour $\ell\geq 0$, $\Hdim\,\gT\leq \ell$ \ssi pour tout 
$x\in\gT$,
$\Hdim(\gT_\rH^x)\leq \ell-1$.
\end{itemize}
\end{definition}

Du lemme \ref{lemTquoBordH} on déduit par \recu sur $\ell$,
avec la même convention de notation qu'en \ref{notaKdiminf} le 
lemme
suivant.

\begin{lemma}
\label{lemDimHquo}
Si $\gL$ est le quotient de $\gT$ par un \idz, alors $\Hdim\,\gL\leq
\Hdim\,\gT$.
\end{lemma}

\begin{proposition}
\label{propJdimHdim}
\emph{(Comparaison des dimensions $\Jdim$ et $\Hdim$)}
%
\begin{enumerate}
\item On a toujours $\Hdim\,\gT\leq  \Jdim\,\gT$.
\item Lorsque $\He(\gT)$ est une \agHz, on a \egtz: $\Hdim\,\gT=  
\Jdim\,\gT$.
\end{enumerate}

\noindent En \clama l'\egt a donc lieu lorsque $\He(\gT)$ est noethérien.
\end{proposition}
\begin{proof}
Le point \textsl{1} se démontre par \recu sur $\Jdim\,\gT$ à partir du 
point \textsl{1} de
la proposition~\ref{propBHeitHeyt}.

\noindent Pour le point \textsl{2}, on suppose que  $\He(\gT)$ est une \agH et 
l'on fait
une \recu en utilisant le point \textsl{2} de la proposition 
\ref{propBHeitHeyt}.
Pour que la \recu fonctionne il faut montrer 
que~\hbox{$(\He(\gT))\ul{\hat{x}}=\He(\rH_\gT^x)$} est également une \agHz.
Or cela résulte de ce que~\hbox{$(\He(\gT))\ul{\hat{x}}$} est un quotient 
de~$\He(\gT)$ par un \id principal (puisque $\He(\gT)$ est une \agHz) et 
du
fait~\ref{factQuoAgH2}.
\end{proof}

\begin{proposition}
\label{propHdim0} Notons $\gT'$ pour le quotient $\gT/(\JT(0)=0)$.
\begin{enumerate}
\item On a toujours $\Hdim\,\gT=\Hdim(\gT')$.
\item On a les équivalences $\Jdim\,\gT\leq 0\;\Longleftrightarrow\; \Hdim\,\gT\leq 
0\;\Longleftrightarrow\;
\Kdim(\gT')\leq 0$, (autrement dit $\gT'$ est une \agBz). C'est le cas 
lorsque  $\He(\gT)$
est fini.
\end{enumerate}
\end{proposition}
\begin{proof}
Le premier point résulte du lemme \ref{lemTBKBH}, point \textsl{2}.\\
Pour le deuxième point on sait déjà que $\Hdim\,\gT\leq 
\Jdim\,\gT\leq
\Kdim\,\gT'$. Le lemme \ref{lemTBKBH}, point \textsl{1}, prouve que 
$\Hdim\,\gT\leq 0$
implique  $\Kdim\,\gT'\leq 0$.\\
Lorsque $\gT$ est fini l'espace $\Jspec\,\gT$ est simplement 
l'ensemble des
\idemas avec pour topologie toutes les parties (puisque les points 
sont
fermés).
Par ailleurs le fait \ref{factHeHe}.\textsl{1} permet de conclure aussi lorsque
$\He(\gT)$ est fini.
\end{proof}

\hum{l'ancienne preuve était bizarrement nettement plus 
compliquée}

\rem
En \clama  $\He(\gT)$ est fini
  \ssi l'ensemble $M$ des
  \idemas est fini (le cas semi-local en \alg
  commutative). 

\medskip La proposition suivante est l'analogue de la proposition
\ref{propBordKUnion}, en remplaçant le bord de Krull par le bord 
de
Heitmann.

\begin{proposition}
\label{propBordHUnion}
Pour tous $x,y\in\gT$ les \ids bords de Heitmann de $x\vi y$ et $x\vu 
y$
contiennent
l'intersection des  \ids bords de Heitmann de $x$ et $y$, plus 
précisément:
$$\rH_\gT^{x}\cap \rH_\gT^{y}\;\;=\;\; \rH_\gT^{x\vu y}\cap 
\rH_\gT^{x\vi y}$$
\end{proposition}

\begin{proof}
Résulte de la proposition \ref{propBordKUnion} et du lemme 
\ref{lemTBKBH}.
\end{proof}

La proposition suivante est l'analogue du point \textsl{1} de la
proposition~\ref{propLocBord}.
\begin{proposition}
\label{propLocBordH}
Soit $(\fa_i)_{1\leq i\leq m}$ une famille finie d'\ids de $\gT$, avec
$\bigcap_{i=1}^m\JT(\fa_i)\subseteq \JT(0)$ (c'est le cas en 
particulier si
$\bigcap_{i=1}^m\fa_i=\{0\}$).
Pour $x\in\gT$ notons encore $x$ son image dans 
le quotient~\hbox{$\gT_{\!i}=\gT/(\fa_i=0)$}.
Alors le bord ${\gT_{\!i}}_\rH^x$ peut être vu comme le quotient 
de 
$\gT_\rH^x$
par un \idz~$\fb_i$ et l'on a: $\bigcap_{i=1}^m\fb_i=\so{0}$.
\end{proposition}
\begin{proof}
Résulte du lemme \ref{lemTBKBH} et de la proposition 
\ref{propLocBord},
appliquée au treillis
$\gT'=\gT/(\JT(0)=0)$ et aux idéaux images des $\JT(\fa_i)$ dans 
$\gT'$.
\end{proof}

Le corolaire suivant est l'analogue du 
corolaire~\ref{corpropLocBord} (avec
les mêmes notations qu'en~\ref{notaKdiminf}).

\begin{corollary}
\label{corpropLocBordH}
Soit $(\fa_i)_{1\leq i\leq m}$ une famille finie d'\ids de $\gT$ et
$\fa=\bigcap_{i=1}^m\fa_i$. \\
Alors
$\Hdim(\gT/(\fa=0))=\sup_i\Hdim(\gT/(\fa_i=0))$.
\end{corollary}
\begin{proof}
  En remplaçant $\gT$ par $\gT/(\fa=0)$ on se ramène au cas 
où $\fa=0$.
La chose est claire pour $\ell=-1$. Et la preuve par \recu fonctionne
gr\^{a}ce à la proposition~\ref{propLocBordH}.
\end{proof}

\begin{proposition}
\label{propHdimgen}
Soit $S$ un système \gtr du \trdi $\gT$ et $\ell\geq 0$.
\Propeq
\begin{enumerate}
\item Pour tout $x\in\gT$,
$\Hdim(\gT_\rH^x)\leq \ell-1$.
\item Pour tout $x\in S$,
$\Hdim(\gT_\rH^x)\leq \ell-1$.
\end{enumerate}
\end{proposition}
\begin{proof}
Cela résulte de la proposition \ref{propBordHUnion}, du lemme
\ref{lemTquoBordH} et du corolaire \ref{corpropLocBordH}: par 
exemple pour \hbox{tous
$x,y\in \gT$}, puisque $\rH_\gT^{x\vu 
y}\subseteq\rH_\gT^x\cap\rH_\gT^y$, le
treillis $\gT_\rH^{x\vu y}$ est un quotient
de $\gT/(\rH_\gT^x\cap\rH_\gT^y)$ par un idéal, donc 
$\Hdim(\gT_\rH^{x\vu
y})\leq \sup (\Hdim\,\gT_\rH^x, \Hdim\,\gT_\rH^y)$.
\end{proof}

  \rem
Explicitons encore un peu plus la dimension de Heitmann.
Pour ceci nous introduisons \gui{l'idéal bord de Heitmann 
itéré}.
Pour $x_1$, \ldots, $x_n\in\gT$ nous notons 
$$
\gT_\rH[x_1]=\gT_\rH^{x_1},\,
\gT_\rH[x_1,x_2]=(\gT_\rH^{x_1})_\rH^{x_2},\,
\gT_\rH[x_1,x_2,x_3]=((\gT_\rH^{x_1})_\rH^{x_2})_\rH^{x_3},\,\ldots
$$ 
les treillis bords quotients successifs. Et $\rH[\gT;x_1,\ldots
,x_k]=\rH_\gT[x_1,\ldots ,x_k]$ désigne le noyau de la projection 
canonique
$\gT\to \gT_\rH[x_1,\ldots ,x_k]$.

\noindent Dire que $\Hdim\,\gT\leq \ell$ revient à dire que pour 
tous
$x_0,\ldots ,x_\ell\in\gT$ on a $1\in\rH_\gT[x_0,\ldots ,x_\ell]$.
Il nous faut donc expliciter les \ids $\rH_\gT[x_0,\ldots 
,x_\ell]$.
Pour ceci nous devons expliciter $\pi^{-1}(\rH[\gT/(\fa=0);\pi(x)])$
(que nous noterons $\rH[\gT,\fa;x]$) lorsqu'on a une projection
$\pi:\gT\to\gT/(\fa=0)$.

\noindent Par \dfn on a $y\in\rH[\gT,\fa;x]$ \ssi $y\leq x\vu 
z\;\mod\;\fa$ pour
un $z$ qui vérifie $\pi(z\vi x)\in \rJ_{\gT/(\fa=0)}(0)$.
Cette dernière condition signifie: $\Tt u\in\gT,\; (\pi((z\vi x)\vu
u=\pi(1)\;\Rightarrow \;\pi(u)=\pi(1)$.
Et ceci s'écrit encore
$$\Tt u\in\gT,\; ((\Ex a\in\fa\;(z\vi x)\vu u\vu a=1)\;\Rightarrow 
\;(\Ex
b\in\fa\; u\vu b=1))$$
Par ailleurs, $y\leq x\vu z\;\mod\;\fa$ signifie $\Ex a'\in\fa\;y\leq 
x\vu z\vu
a'$ et la condition   $\pi(z\vi x)\in \rJ_{\gT/(\fa=0)}(0)$ n'est pas 
changée
si on remplace $z$ par $z\vu a'$. En bref nous obtenons la condition 
suivante
\hbox{pour  $y\in\rH[\gT,\fa;x]$}:
$$ \Ex z\in\gT\; [\,y\leq x\vu z\;\mathrm{et}\;
\Tt u\in\gT,\; ((\Ex a\in\fa\;(z\vi x)\vu u\vu a=1)\;\Rightarrow 
\;(\Ex
b\in\fa\; u\vu b=1))\,]$$
On voit donc apparaître une formule d'une complexité logique 
redoutable.
Surtout si on songe que $\Ex a\in\fa$ et $\Ex b\in\fa$ devront 
être
explicités avec $\fa=\rH_\gT[x_1,\ldots ,x_k]$ si on veut obtenir un
expression pour
$y\in\rH_\gT[x_1,\ldots ,x_k,x]$.
Contrairement à l'expression pour $\Jdim\,\gT\leq \ell$ qui ne 
comportait que
deux alternances de quantificateurs quel que soit l'entier $\ell$, on 
voit pour
$\Hdim\leq \ell$ des expressions de plus en plus imbriquées au fur 
et à
mesure que $\ell$ augmente. En fait il se trouve que pour les anneaux
commutatifs, c'est la $\Hdim$ qui fait fonctionner les preuves par 
\recu dans
les sections~\ref{secKroBass} et suivantes, et c'est la vraie raison 
pour
laquelle on est amené à introduire cette dimension. Comme elle 
vérifie
$\Hdim\,\gT\leq \Jdim\,\gT\leq \Kdim(\gT/(\JT(0)=0))$
on a quand même des moyens raisonnables pour la majorer.
Mais on manque d'exemples avec une majoration meilleure que celle par
$\Kdim(\gT/(\JT(0)=0))$.

\hum{~

1. On peut se poser la question de comparer $\Hdim\,\gT$ et  
$\Hdim\,\gT\cir$.

\smallskip
2. On pourrait peut être parler des treillis locaux\,? 
Constructivement on peut
distinguer entre $\Tt x, y\in\gT\;(x\vu y=1)\Rightarrow 
(x=1\;\mathrm{ou}\;y=1)$
(treillis local) et
$\gT/(\JT(0)=0)\simeq\Deux$ (treillis local résiduellement discret)}
\section[Dimensions de Krull et Heitmann: anneaux 
commutatifs]{Dimensions de
Krull et Heitmann pour les anneaux commutatifs}
\label{secBPA}
Dans cette section, $\gA$ désigne toujours un anneau commutatif.

\subsection{Le treillis de Zariski}

Nous rappelons ici l'approche \cov de  \cite{Joy76} pour le spectre
d'un anneau commutatif.

Si $J\subseteq \gA$, nous notons $\cI_\gA(J)$ ou $\gen{J}_\gA$ (ou
$\gen{J}$ si le contexte est clair) l'\id engendré par~$J$; nous
notons $\DA(J)$ (ou $\rD(J)$ si le contexte est clair) le
nilradical de l'\id $\gen{J}$:
\begin{equation} \label{eqZar}
\begin{array}{rclcl}
\DA(J)&  = & \sqrt[\gA]{\gen{J}} &=&\sotq{x\in\gA}{\Ex m\in\NN\;\; 
x^m\in\gen{J}}
\end{array}
\end{equation}
Lorsque $J=\so{x_1,\ldots ,x_n}$ nous notons
$\DA(x_1,\ldots ,x_n)$ pour $\DA(J)$.
Si le contexte est clair, nous abrégeons $\DA(x)$ en $\wi{x}$.
Les  $\wi{x}$  forment un système \gtr de $\ZarA$, stable par $\vi$.

Par \dfn le {\sl treillis de Zariski} de $\gA$, noté $\ZarA$, a pour
\elts les radicaux d'\itfsz: ce sont donc les $\DA(x_1,\ldots ,x_n)$,
\cad les $\DA(\fa)$ pour les \itfs $\fa$.  La relation d'ordre est
l'inclusion, le inf et le sup sont donnés par
$$
\DA(\fa_1)\vi\DA(\fa_2)=\DA(\fa_1\fa_2)\quad \mathrm{et} \quad
\DA(\fa_1)\vu\DA(\fa_2)=\DA(\fa_1+\fa_2).
$$
Le treillis de Zariski de $\gA$
est un \trdiz, et $\DA(x_1,\ldots ,x_n)=
\wi{x_1}\vu\cdots \vu\wi{x_n}.$

Si $J\subseteq \gA$  nous notons $\wi{J}=\sotq{\wi{x}}{x\in J}
\subseteq\ZarA$.

Soient $U$ et  $J$  deux familles finies dans $\gA$, on a les équivalences
$$ \Vi\wi{U}\leq_{\ZarA} \Vu\wi{J}
\quad\Longleftrightarrow \quad
\prod\nolimits_{u\in U} u  \in \sqrt{\gen{J}}
\quad\Longleftrightarrow \quad
\cM(U)\cap \gen{J}\neq \emptyset
$$
où $\cM(U)$ est le \mo multiplicatif engendré par $U$.

Cela suffit à décrire le treillis $\ZarA$. Plus précisément 
on a
(cf. \cite{CC00,CL2003}):
\begin{proposition}
\label{propZar} Le treillis $\ZarA$  est
(à \iso près) le treillis engendré par des symboles $\DA(a)$
soumis aux relations suivantes
$$\begin{array}{cccc}
\DA(0_\gA) =0   \;,\; \DA(1_\gA)= 1 \;,\;
   \DA(x+y) \leq \DA(x)\vu\DA(y) \;,\; \DA(xy) = \DA(x)\vi \DA(y)
\end{array}$$
\end{proposition}

L'opération $\Zar$ est un foncteur de la catégorie des anneaux 
commutatifs
vers celle des \trdisz.
Notez que via ce foncteur la projection $\gA\to\gA/\DA(0)$ donne un 
\iso
$\ZarA\simeq \Zar(\gA/\DA(0))$. On a $\ZarA=\Un$ \ssi $1_\gA=0_\gA$.

\smallskip Un \tho important de  \cite{Hoc1969} affirme que tout 
espace
spectral
est homéomorphe au spectre d'un anneau commutatif.
Voici une version sans point du \tho de Hochster: 

\smallskip \noindent \textsl{Tout \trdi est 
isomorphe au
treillis de Zariski d'un anneau commutatif} 

\smallskip \noindent 
Pour une preuve non \cov  voir \cite{Ban96}.

\subsection{Idéaux, filtres et quotients de $\ZarA$}

Rappelons qu'en \clama le \textsl{spectre de Zariski} $\Spec\,\gA$ d'un anneau
commutatif est un espace topologique dont les points sont les 
\ideps de l'anneau et dont la topologie est définie par la base 
d'ouverts formée par les $\fD_\gA(a)=\sotq{\fp\in\Spec\,\gA}{a\notin\fp}$.
On note aussi $\fD_\gA(x_1,\ldots ,x_n)$ pour 
$\fD_\gA(x_1)\cup\cdots\cup \fD_\gA(x_n)$.

\subsubsection*{Idéaux de $\gA$ et de $\ZarA$}

Nous disons qu'un \id $I$ de $\gA$ est \textsl{radical} si $I=\sqrt[\gA]{I}.$
En \clamaz, tout \id radical est l'intersection des \ideps qui  le 
contiennent.

Introduisons la notation (lorsque $J\subseteq\gA$)
$$
\IZA(J)=\cI_{\ZarA}(\wi{J})
$$
pour l'\id de $\ZarA$ engendré par $\wi{J}$.  En particulier
$$
\IZA(\so{x_1,\ldots ,x_n})=\dar\DA(x_1,\ldots ,x_n)=
\dar(\wi{x_1}\vu\cdots \vu\wi{x_n})\,.
$$

On a $\IZA(J)=\IZA(\sqrt{\gen{J}})$, et on établit facilement le
fait fondamental suivant.

\begin{fact}
\label{factSpecAzarA} ~  
\begin{itemize}
\item L'application $\fa\mapsto \IZA(\fa)$ définit un \iso du
treillis des \ids radicaux de $\gA$ vers le treillis des \ids de
$\ZarA$.  
\item Par restriction les \ideps (resp.  les \idemasz) de
l'anneau $\gA$ et ceux du \trdi $\ZarA$ sont également en
correspondance naturelle bijective.  
\item Pour tout anneau commutatif
$\gA$, $\Spec\,\gA$ (au sens des anneaux commutatifs) s'identifie à
$\Spec(\ZarA)$ (au sens des \trdisz).
\end{itemize}
\end{fact}

\rem En \clama on a un \iso entre le treillis $\ZarA$ et le treillis des
\oqcs de $\Spec\,\gA$.  On peut alors identifier 
\begin{itemize}
\item $\DA(x_1,\ldots ,x_n)$, qui est un \elt de $\ZarA$,
\item  $\fD_{\!\ZarA}(\DA(x_1,\ldots ,x_n))$, qui est un \oqc de $\Spec(\ZarA)$,
\item    et  $\fD_\gA(x_1,\ldots ,x_n)$, qui est un \oqc de $\Spec\,\gA$.
\end{itemize}

\noindent Du point de vue \cofz, on considère $\Spec\,\gA$ comme un \gui{espace
topologique sans point}, \cad un espace défini uniquement à
travers une base d'ouverts, et 
les \gui{identifications} ci-dessus sont de pures identités.

\smallskip On a aussi facilement:
\begin{fact}
\label{factQuoAT} \emph{(Quotients)}\\
Si $J\subseteq\gA$, alors $\Zar(\aqo{\gA}{J}) \simeq \Zar(\gA/\DA(J))
\simeq \Zar(\gA)/(\IZA(J)=0)$.
\end{fact}

\begin{fact}
\label{factTransporteurs} \emph{(Transporteurs)}\\
Soient $\fA$ et $\fB$ des \ids de $\gA$, $\fa=\DA(\fA)$ et
$\fb=\DA(\fB)$.  
Alors $\fa:\fb=\fa:\fB$ est un \id radical de $\gA$
et dans $\ZarA$ on a $\IZA(\fa):\IZA(\fb)=\IZA(\fa:\fb)$.
\end{fact}

\begin{fact}
\label{factRecouvI} \emph{(Recouvrement par des \idsz)}\\
Soit $\fa_i$ une famille finie d'\ids de $\gA$.  Les $\IZA(\fa_i)$
recouvrent $\ZarA$ (\cad leur intersection est réduite à $0$) \ssi
$\bigcap_i\fa_i\subseteq \DA(0)$.
\end{fact}

En \clama le treillis $\Zar\,\gA$ est noethérien (ce qui 
revient à
dire que~$\Spec\,\gA$ est noethérien) \ssi tout \id radical est
\gui{radicalement de type fini}, \cad est un \elt de~$\Zar\,\gA$.

Par ailleurs,  $\Zar\,\gA$ est une \agH \ssi
\hbox{$\Tt \fa,\fb\in\Zar\,\gA\;(\fa:\fb)\in\Zar\,\gA$}.

Le résultat suivant est important en \comaz.

\begin{proposition}
\label{propZarHeyt} (cf. \cite{CL2003})
Si $\gA$ est un anneau noethérien cohérent $\Zar\,\gA$ est une 
\agHz.
Si en outre $\gA$ est fortement discret, la relation d'ordre dans  
$\Zar\,\gA$
est décidable. On dit alors que le treillis est \emph{discret}.
\end{proposition}

\subsubsection*{Filtres de $\gA$ et  de $\ZarA$}

Un \textsl{filtre} dans un anneau commutatif est un \mo $\fF$ qui 
vérifie
$xy\in\fF\Rightarrow x\in\fF$. Un \textsl{filtre premier} est un filtre 
qui
vérifie $x+y\in\fF\Rightarrow x\in\fF\;\mathrm{ou}\;y\in\fF$ (c'est 
le
complémentaire d'un \idepz).

Pour $x\in \gA$ le filtre $\uar \wi x$ de $\ZarA$ est noté 
$\FZ_\gA(x)$.
Plus généralement pour  $S\subseteq\gA$ on note
$\FZ_\gA(S)$ le filtre de $\ZarA$:
$$
\FZ_\gA(S)=
\bigcup\nolimits_{x\in\cM(S)}\uar \wi x.
$$
On a aussi $\FZ_\gA(S) =\FZ_\gA(\fF)=\bigcup\nolimits_{x\in\fF}\!\uar \wi x$, où $\fF$ est le filtre de  $\gA$ engendré par $S$.

Les faits suivants sont faciles.

\begin{fact}
\label{factFZ}
L'application $\ff\mapsto \FZ_\gA(\ff)$ établit une correspondance
\emph{injective} croissante des filtres de $\gA$ vers les filtres de 
$\ZarA$, et
un  sup fini (le sup de $\ff_1$ et $\ff_2$  est engendré par les 
$f_1f_2$ où
$f_i\in\ff_i$) donne pour image le sup fini des filtres images.
Cette correspondance $\FZ_\gA$  se restreint en une bijection entre 
les filtres
premiers de $\gA$ et ceux de $\ZarA$.
\end{fact}

Notez cependant que le filtre principal de $\ZarA$ engendré par
$\wi{a_1}\vu\cdots \vu \wi{a_n}$ (\cad l'intersection des filtres 
$\uar\wi
{a_i}$), ne correspond en général à aucun filtre de $\gA$.

\begin{fact}
\label{factLocalises} \emph{(Localisés)} \\
Soit $S$ un \mo de $\gA$, $\fF$ le filtre engendré par $S$, et
$\ff=\FZ_\gA(S)=\FZ_\gA(\fF)$. \\
Alors $S^{-1}\gA=\gA_S=\gA_\fF$ et $\Zar(\gA_S)\simeq\Zar(\gA)/(\ff=1)$.
\end{fact}

\begin{fact}
\label{factComplement} \emph{(Filtre complémentaire)}\\
Soit $x\in\gA$, le filtre $1_\ZarA\setminus \FZ_\gA(x)$
est égal à $\FZ_\gA(1+x\gA)$.
\end{fact}

\begin{fact}
\label{factRecouvF} \emph{(Recouvrement par des filtres)}\\
Soit $(S_i)_{1\leq i\leq n}$ une famille finie de \mos de $\gA$. Les 
filtres
$\FZ_\gA(S_i)$ recouvrent $\ZarA$ (\cad leur intersection est 
réduite à
$\so{1}$) \ssi les \mos $S_i$ sont \com \cad que pour tous $x_i\in 
S_i$
on a $\gen{x_1,\ldots ,x_n}=\gen{1}$. \\
Plus généralement on a $\,\FZ_\gA(S_1)\,\cap \cdots
\cap\,\FZ_\gA(S_n)\subseteq \FZ_\gA(S)$ \ssi pour  tous $x_i\in S_i$
il existe $x\in S$ tel que $x\in\gen{x_1,\ldots ,x_n}$.
\end{fact}

\subsection{Le treillis de Heitmann}

Dans un anneau commutatif, le \textsl{radical de Jacobson d'un \id 
$\fJ$} est (du point de vue des \clamaz) l'intersection des \idemas qui contiennent $\fJ$. On
le note $\JA(\fJ)=\rJ(\gA,\fJ)$, ou encore $\rJ(\fJ)$ si le contexte est 
clair.
En \coma on utilise la \dfn suivante, classiquement équivalente:
\begin{equation} \label{eqRadJac}
\JA(\fJ)\eqdefi\sotq{x\in\gA}{\Tt y\in\gA,\;\; 1+xy
\hbox{ est inversible modulo } \fJ}
\end{equation}
On notera $\JA(x_1,\ldots ,x_n)=\rJ(\gA,x_1,\ldots ,x_n)$ pour
$\JA(\gen{x_1,\ldots ,x_n})$.  L'\id  $\JA(0)$ est appelé le
\textsl{radical de Jacobson de l'anneau $\gA$}.

\begin{definition}
\label{defHeitA}
On appelle \textsl{treillis de Heitmann} d'un anneau commutatif $\gA$ 
le treillis
$\He(\ZarA)$, on le note $\HeA$.
\end{definition}

En \clamaz, vu le fait \ref{factSpecAzarA} et vue la \dfn du radical 
de Jacobson
via les intersections d'\ids maximaux, le fait suivant,
qui conduit à une interprétation simple du treillis $\HeA$, est 
évident.
Nous sommes néanmoins intéressés par une \prco directe.

\begin{fact}
\label{factRadJac} \emph{(Radical de Jacobson)}\\
La correspondance bijective $\IZA$ préserve le passage au radical 
de Jacobson.
Autrement dit si $\fJ$ est un \id de $\gA$ et $\fj=\IZA(\fJ)$, alors
$\rJ_\ZarA(\fj)=\IZA(\JA(\fJ))$.
\end{fact}
\begin{proof}
Il suffit de montrer que pour tout $x\in\gA$, $\wi 
x\in\rJ_\ZarA(\fj)$ \ssi $\wi
x\in\IZA(\JA(\fJ))$. Puisque $\JA(\fJ)$ est un \id radical, on 
cherche donc à
montrer l'équivalence
$$\Tt x\in\gA\quad (\,\wi x\in\rJ_\ZarA(\fj)\;\Leftrightarrow\; 
x\in\JA(\fJ)\,)
$$
Par \dfn $\wi x\in\rJ_\ZarA(\fj)$ signifie
$$\Tt y\in\ZarA\quad (\wi x\vu y=1_\ZarA\;\Rightarrow\;\Ex 
z\in\fj\;\; z\vu
y=1_\ZarA)$$
\cade puisque tout $y\in\ZarA$ est de la forme $\DA(y_1,\ldots ,y_k)$,
$$\Tt y_1,\ldots ,y_k\in\gA\quad (\gen{x,y_1,\ldots ,y_k}=1_\gA\;
\Rightarrow\;\Ex z\in\fj\;\; z\vu \DA(y_1,\ldots ,y_k)=1_\ZarA)$$
ceci est immédiatement équivalent à
$$\Tt y_1,\ldots ,y_k\in\gA\quad (\gen{x,y_1,\ldots ,y_k}=1_\gA\;
\Rightarrow\;\Ex u\in\fJ \;\;\gen{u,y_1,\ldots ,y_k}=1_\gA)$$
puis à
$$\Tt y\in\gA\;(\gen{x,y}=1_\gA\;
\Rightarrow\;\Ex u\in\fJ \;\;\gen{u,y}=1_\gA)$$
ou encore à: tout  $y\in\gA$ de la forme $1+xa$ est inversible 
modulo $\fJ$.
\Cad $x\in\JA(\fJ)$.
\end{proof}

\begin{corollary}
\label{propHeitA}
Soient $\fj_1$ et $\fj_2$ deux \itfs de $\gA$. Les \elts $\DA(\fj_1)$ 
et
$\DA(\fj_2)$ de $\ZarA$ sont égaux dans le quotient $\HeA$ \ssi
$\JA(\fj_1)=\JA(\fj_2)$. En conséquence $\HeA$ s'identifie
à l'ensemble des $\JA(x_1,\ldots ,x_n)$, avec
$\JA(\fj_1)\vi\JA(\fj_2)=\JA(\fj_1\fj_2)$ et
$\JA(\fj_1)\vu\JA(\fj_2)=\JA(\fj_1+\fj_2)$.
\end{corollary}

\rems ~\\
1. Vu les bonnes propriétés de la correspondance $\IZA$ on
a, avec $\gT=\ZarA$, $\Zar(\gA/\JA(0))\simeq\gT/(\JT(0)=0)$. Par 
contre il ne
semble pas qu'il y ait une $\gA$-algèbre $\gB$ naturellement 
attachée à~$\gA$ pour laquelle on ait $\Zar\,\gB\simeq\He(\ZarA)$.\\
2. Notez que, en général $\JA(x_1,\ldots ,x_n)$ est un \id 
radical mais pas
le radical d'un \itfz. \\
3. On voit aussi facilement que
$\JA(\fj_1)\vi\JA(\fj_2)=\JA(\fj_1)\cap\JA(\fj_2)=\JA(\fj_1\cap\fj_2)$ 
(cela
résulte d'ailleurs du lemme \ref{lemJacInter}). Il peut sembler 
surprenant que
$\JA(\fj_1)\cap\JA(\fj_2)=\JA(\fj_1\fj_2)$ (c'est à priori moins 
clair que pour
les $\DA$). Voici le calcul élémentaire qui (re)démontre ce 
fait.
On a $x\in\JA(\fj_1)$ \ssi $\Tt y\;(1+xy)$ est inversible modulo 
$\fj_1,$ et
$x\in\JA(\fj_2)$ \ssi $\Tt y\;(1+xy)$ est inversible modulo $\fj_2$.
Mais si $a=1+xy$ est inversible modulo $\fj_1$ et $\fj_2$, il est 
inversible
modulo leur produit: en effet $1+aa_1\in\fj_1$ et  $1+aa_2\in\fj_2$  
impliquent
que $(1+aa_1)(1+aa_2)$, qui se réécrit $1+aa'$, est dans 
$\fj_1\fj_2$.

\subsection{Dimension et  bords de Krull}

En \coma on donne la \dfn suivante.

\begin{definition}
\label{defKdimA}
La dimension de Krull d'un anneau commutatif
est la dimension de Krull de son treillis de Zariski.
\end{definition}

Vus le fait \ref{factSpecAzarA} et le \tho \ref{thDK1}, il s'agit 
d'une \dfn
équivalente à la \dfn usuelle en \clamaz.

\begin{definition}
\label{defZar2} Soit  $\gA$  un anneau commutatif, $x\in\gA$ et $\fj$ 
un \itfz.
\begin{itemize}
\item [$(1)$] Le \textsl{bord supérieur de Krull de $\fj$ dans $\gA$} 
est
l'anneau quotient
$\gA\ul{\fj}:=\gA/\rK_\gA(\fj)$ où
\begin{equation} \label{eqBKAC}
\rK_\gA(\fj):=\fj+(\DA(0):\fj)
\end{equation}
  On note aussi $\gA\ul{x}=\gA\ul{x\gA}$ et on l'appelle le \textsl{bord
supérieur de $x$ dans $\gA$}. On dira aussi que $\rK_\gA(\fj)$ est 
\textsl{l'\id
bord de Krull de $\fj$.}  On notera aussi  $\rK_\gA(y_1,\ldots ,y_n)$ 
pour
$\rK_\gA(\gen{y_1,\ldots ,y_n})$ et $\rK_\gA^x$ pour $\rK_\gA(x)$.
\item [$(2)$] Le \textsl{bord inférieur de Krull de $x$ dans $\gA$} 
est l'anneau
localisé
$\gA\bal{x}:=\gA_{\rS\bal{x}}$ où $\rS\bal{x}=x^\NN(1+x\gA)$.
On dira aussi que le \mo $x^\NN(1+x\gA)$ est le \textsl{\mo bord de 
Krull de~$x$.}
\end{itemize}
\end{definition}

Ainsi un \elt arbitraire de  $\rK_\gA(y_1,\ldots ,y_n)$ s'écrit 
$\sum_i
a_iy_i+b$ avec tous les $by_i$ nilpotents.

\begin{proposition}
\label{propZar2} Soit   $x\in \gA$ et  $\fj=\gen{j_1,\ldots ,j_n}$ un 
\itfz.
Considérons $\wi x\in\ZarA$ et $\fa=\DA(\fj)=\wi{j_1}\vu\cdots
\vu\wi{j_n}\in\ZarA$. Alors
\begin{enumerate}
\item L'idéal bord de Krull de $\fa=\DA(\fj)$ dans $\ZarA$,
$\rK_\ZarA^\fa$, est égal à
$\IZA(\rK_\gA(\fj))$. En conséquence  $(\ZarA)\ul{\fa}$ 
s'identifie
naturellement avec $\Zar(\gA\ul{\fj})$.
\item Le filtre  bord de Krull de $\wi x$ dans $\ZarA$, 
$\rK^\ZarA_{\tilde x}$, est
égal à
$\FZ(\rS\bal{x})$. En conséquence   $(\ZarA)\bal{\tilde x}$ s'identifie
naturellement avec
$\Zar(\gA\bal{x})$.
\end{enumerate}
\end{proposition}
\begin{proof}
Pour l'idéal bord, on a par \dfn
  $$\rK_\ZarA^\fa=\rK_\ZarA(\DA(\fa))=\DA(\fa) \vu 
(\DA(0):\DA(\fa) ).
$$
D'après les faits \ref{factSpecAzarA} et \ref{factTransporteurs},
il est égal à $\IZA(\fa + (\DA(0):\fa ))$ et aussi à
  $\IZA(\fj + (\DA(0):\fj ))=\IZA(\rK_\gA^\fj)$. Enfin pour les 
passages au
quotient on applique le fait~\ref{factQuoAT}.\\
Pour le filtre bord de Krull, cela fonctionne de la même 
manière en
utilisant les faits \ref{factFZ} et \ref{factComplement} puis en 
passant au
treillis quotient avec le fait~\ref{factLocalises}.
\end{proof}

\hum{pour le bord inférieur avec un \itf à la place d'un \elt il 
semble qu'on
obtient non pas un anneau, mais un schéma de Grothendieck obtenu en 
recollant
un nombre fini d'anneaux le long de localisations convenables}

Comme corolaire des propositions \ref{propDK1} et \ref{propZar2}  on 
obtient
l'analogue suivant du \tho \ref{thDK1}, dans une version 
entièrement \covz.
Rappelons que la dimension de Krull d'un anneau \hbox{égale $-1$} \ssi 
l'anneau est
trivial (i.e., $1_\gA=0_\gA$).

\begin{theorem}
\label{thDKA} Pour un anneau commutatif  $\gA$ et un entier $\ell\geq 
0$ \propeq
\begin{enumerate}
\item La dimension de Krull de $\gA$ est $\leq \ell$.
\item Pour tout $x\in \gA$ la dimension de Krull de $\gA\ul x$ est 
$\leq \ell-
1$.
\item Pour tout \itf $\fj$ de $\gA$ la dimension de Krull de 
$\gA\ul\fj$ est
$\leq \ell-1$.
\item Pour tout $x\in \gA$ la dimension de Krull de $\gA\bal x$ est 
$\leq \ell-
1$.
\end{enumerate}
\end{theorem}

Ce \tho nous donne une bonne signification intuitive de la dimension 
de Krull.

Avec le fait \ref{factSpecAzarA} on obtient le même \tho en 
\clamaz.

Vu son importance, nous allons donner des preuves directes simples des équivalences entre les points \textsl{1}, \textsl{2} et \textsl{4} en \clamaz.

\begin{Proof}{Démonstration directe en \clamaz}
Montrons d'abord l'équivalence des points \textsl{1} et \textsl{2}.
Rappelons que les \ideps de $S^{-1}\gA$ sont de la forme $S^{-1}\fp$ 
où $\fp$
est un \idep
de $\gA$ qui ne coupe pas $S$.
L'équivalence résulte alors clairement des deux affirmations 
suivantes.\\
(a) Soit $x\in\gA$, si $\fm$ est un idéal maximal de $\gA$ il coupe 
toujours
$\rS\bal{x}$. En effet si $x\in\fm$ c'est clair et sinon, $x$ est 
inversible
modulo $\fm$ ce qui signifie que $1+x\gA$ coupe $\fm$.\\
(b) Si $\fm$ est un idéal maximal de $\gA,$  et si 
$x\in\fm\setminus\fp$ où
$\fp$ est un \idep contenu dans $\fm$,
alors $\fp\cap \rS\bal{x}=\emptyset$: en effet si $x(1+xy)\in\fp$ 
alors,
puisque $x\notin\fp$ on a $1+xy\in\fp\subset\fm$, ce qui donne
la contradiction $1\in\fm$ (puisque $x\in\fm$).\\
Ainsi, si $\fp_0\subsetneq \cdots \subsetneq \fp_\ell$
est une chaîne avec
$\fp_\ell$ maximal, elle est raccourcie d'au moins son dernier terme 
lorsqu'on
localise en $\rS\bal{x}$, et elle n'est raccourcie
que de son dernier terme si
$x\in\fp_\ell\setminus\fp_{\ell-1}$.\\
L'équivalence des points \textsl{1} et \textsl{4} se démontre de manière 
\gui{opposée},
en remplaçant les idéaux premiers par les filtres premiers.
On remarque d'abord que les filtres premiers de $\gA/\fJ$ sont de la 
forme
$(S+\fJ)/\fJ$, où $S$ est un filtre premier de $\gA$ qui ne coupe 
pas $\fJ$.
Il suffit alors de démontrer les deux affirmations \gui{opposées} 
de (a) et
(b) qui sont les suivantes:\\
(a') Soit $x\in\gA$, si $S$ est un filtre maximal de $\gA$ il coupe 
toujours
$\rK_\gA^x$. En effet si $x\in S$ c'est clair et sinon,
puisque $S$ est maximal $Sx^\NN$ contient $0$, ce qui signifie qu'il
y a un entier $n$ et un \elt $s$ de $S$ tels que $sx^n=0$.
Alors $(sx)^n=0$ et $s\in (\sqrt{0}:x)\subset \rK_\gA^x$.\\
(b') Si $S$ est un filtre maximal de $\gA,$  et si $x\in S\setminus 
S'$ où
$S'\subset S$ est un filtre premier,
alors $S'\cap \rK_\gA^x=\emptyset$. En effet si $ax+b\in S'$ avec 
$(bx)^n=0$
alors, puisque $x\notin S'$ on a $ax\notin S'$ et, vu que $S'$ est 
premier,
$b\in S'\subset S$, mais comme $x\in S$, $(bx)^n=0\in S$ ce qui est 
absurde.
\end{Proof}

En outre le \tho \ref{thDKA} implique la caractérisation \cov \elr 
de cette
dimension  en terme d'identités algébriques, décrite 
dans~\cite{Lom02,CL2003}, comme suit.

\begin{corollary}
\label{corthDKA} \Propeq
\begin{itemize}
\item  [$(1)$] La dimension de Krull de $\gA$ est $\leq \ell$
\item  [$(2)$] Pour tous $x_0,\ldots ,x_\ell\in\gA$ il existe 
$b_0,\ldots,b_\ell\in \gA$ tels que 
\begin{equation} \label{eqCG}
\left.
\begin{array}{rcl} 
\DA(b_0x_0)& =  &\DA(0)    \\ 
\DA(b_1x_1)& \leq  & \DA(b_0,x_0)  \\
\vdots\quad& \vdots  &\quad  \vdots \\
\DA(b_\ell x_\ell )& \leq  & \DA(b_{\ell -1},x_{\ell -1})  \\
\DA(1)& =  &  \DA(b_\ell,x_\ell )
\end{array}
\right\}
\end{equation}
\item  [$(3)$] Pour tous $x_0,\ldots ,x_\ell\in\gA$ il existe 
$a_0,\ldots,a_\ell\in \gA$ et
$m_0,\ldots,m_\ell\in\NN$ tels que
$$ x_0^{m_0}(x_1^{m_1}\cdots(x_\ell^{m_\ell} (1+a_\ell x_\ell) + 
\cdots+a_1x_1)
+ a_0x_0) =0
$$
\end{itemize}
\end{corollary}
\begin{proof}
Montrons l'equivalence de (1) et (3). 
Utilisons par exemple pour (1) la caractérisation via les localisés  $\gA\bal x$.
L'équivalence pour la dimension $0$ est claire. Supposons la chose 
établie
pour la dimension $\leq \ell$. On voit alors que $S^{-1}\gA$ est de 
dimension
  $\leq \ell$ \ssi pour tous
$x_0,\ldots ,x_\ell\in \gA$ il existe $a_0,\ldots,a_\ell\in \gA$,
$m_0,\ldots,m_\ell\in\NN$  et $s\in S$ tels que
$$ x_0^{m_0}(x_1^{m_1}\cdots(x_\ell^{m_\ell} (s+a_\ell x_\ell) + 
\cdots+a_1x_1) +
a_0x_0)=0.$$
Il reste donc à remplacer $s$ par un \elt arbitraire de la forme
$x_{\ell+1}^{m_{\ell+1}} (1+a_{\ell+1} x_{\ell+1})$. \\
On a  $(3) \Rightarrow (2)$ en prenant:
$b_\ell=1+a_\ell x_\ell$, et 
$b_{k -1} = x_k^{m_k} b_k+ a_{k -1}x_{k -1}$, pour $k=\ell,$
\ldots, $1$. \\
On a  $(2) \Rightarrow (1)$ en considérant la caractérisation (4) 
de la dimension de Krull
d'un \trdi donnée dans la proposition \ref{propDK1} 
et en l'appliquant au treillis de Zariski $\ZarA$ avec 
$S=\sotq{\DA(x)}{x\in\gA}$.
On pourrait aussi v'erifier par un calcul direct que $(2) \Rightarrow 
(3)$.
\end{proof}

\rem Le système d'inégalités (\ref{eqCG}) dans le point (2) du 
corolaire précédent établit une relation intéressante et
symétrique entre les deux suites 
$(b_0,\ldots ,b_\ell)$ et $(x_0,\ldots ,x_\ell)$. Lorsque $\ell=0$, 
cela signifie $\DA(b_0)\vi\DA(x_0)=0$ et $\DA(b_0)\vu\DA(x_0)=1$, 
\cad 
que les deux \eltsz~$\DA(b_0)$ et $\DA(x_0)$ sont compléments l'un de l'autre
dans $\ZarA$.
Dans $\Spec\,\gA$ cela signifie que les ouverts de base 
correspondants sont complémentaires. Nous introduisons donc  
la terminologie suivante: lorsque deux suites 
$(b_0,\ldots ,b_\ell)$ et $(x_0,\ldots ,x_\ell)$ vérifient les 
inégalités (\ref{eqCG}) nous dirons qu'elles sont 
\textsl{complémentaires}.

\medskip 
Signalons aussi qu'il est facile d'établir \cot que
$\Kdim(\gK[X_1,\ldots,X_n])=n$ lorsque~$\gK$ est un corps, ou même
un anneau zéro dimensionnel (cf.~\cite{CL2003}).  On peut aussi traiter
de façon \cov la dimension de Krull des anneaux géométriques
(les $\gK$-algèbres de présentation finie).
Donc les \thos des sections \ref{secKroBass} et suivantes ont un 
contenu algorithmique clair pour ces anneaux (comme pour tout anneau où 
l'on est capable d'expliciter \cot la dimension de Krull).

\medskip
\rem
On a aussi (déjà démontré pour les \trdisz) les résultats 
suivants:
\begin{itemize}
\item si $\gB$ est un quotient ou un localisé de $\gA$,  
$\Kdim\,\gB\leq
\Kdim\,\gA$,
\item si $(\fa_i)_{1\leq i\leq m}$ une famille finie d'\ids de $\gA$ 
et
$\fa=\bigcap_{i=1}^m\fa_i$, alors
$\Kdim(\gA/\fa)=\sup_i\Kdim(\gA/\fa_i)$.
\item si $(S_i)_{1\leq i\leq m}$ une famille finie de \moco de $\gA$ 
alors
$\Kdim(\gA)=\sup_i\Kdim(\gA_{S_i})$.
\item en \clama on a
$\Kdim(\gA)=\sup_\fm\Kdim(\gA_{\fm})$, où $\fm$ partcourt tous les 
\idemasz.
\end{itemize}

\medskip
\rem
On peut illustrer le corolaire \ref{corthDKA} ci-dessus
en introduisant \gui{l'idéal bord de Krull itéré}.
Pour $x_1,\ldots ,x_n\in\gA$ considérons
$(\gA\ul{x_1})\ul{x_2}$, $((\gA\ul{x_1})\ul{x_2})\ul{x_3}$, 
etc\ldots\,\,  les
anneaux bords supérieurs successifs, et notons $\rK_\gA[x_1,\ldots 
,x_\ell]$
le noyau de la projection canonique $\gA\to
(\dots(\gA\ul{x_1})^{\cdots})\ul{x_\ell}$.  Alors on a 
$y\in\rK_\gA[x_0,\ldots
,x_\ell]$ \ssi $\Ex a_0,
\ldots,  a_\ell\in \gT$ et $m_0,\ldots,m_\ell\in\NN$  vérifiant:
$$  x_0^{m_0}(x_1^{m_1}\cdots(x_\ell^{m_\ell} (y+a_\ell x_\ell) + 
\cdots+a_1x_1)
+ a_0x_0) =0.
$$
Et la dimension de Krull est $\leq \ell$ \ssi pour tous $x_0,\ldots
,x_\ell\in\gA$ on a $1\in\rK_\gA[x_0,\ldots ,x_\ell]$.

\subsection{Dimensions de Heitmann}
\subsubsection*{Le spectre de Heitmann}

L'espace spectral que Heitmann a défini pour remplacer le 
j-spectrum,
\cad l'adhérence pour la topologie constructible du spectre maximal 
dans
$\Spec\,\gA$,
  correspond à la \dfn suivante.
\begin{definition}
\label{defJspecA}
On appelle \textsl{spectre de Heitmann} d'un anneau commutatif $\gA$ le 
sous-espace $\Jspec(\ZarA)$ de $\Spec\,\gA$. On le note aussi 
$\Jspec\,\gA$. On note
$\jspec\,\gA$ pour $\jspec(\ZarA)$, \cad le j-spectrum de l'anneau au 
sens
usuel.
\end{definition}

En \clama le \tho \ref{thDK3} donne:
\begin{fact}
\label{factHSpecA}
Pour tout anneau commutatif $\gA$, le spectre de Heitmann de $\gA$ 
s'identifie
à l'espace spectral  $\Spec(\HeA)$ (au sens des \trdisz).
\end{fact}

On a alors la \dfn \cov \elr sans points de la dimension introduite 
par
Heitmann.
\begin{definition}
\label{defJdimA}
La $\rJ$-dimension de Heitmann de $\gA$, notée $\Jdim\,\gA$, est la 
dimension
de Krull de $\Heit(\gA)$, autrement dit c'est la $\Jdim$ de $\ZarA$.
\end{definition}

En \clama $\Jdim\,\gA$ est égal à la dimension de l'espace 
spectral
$\Jspec\,\gA$, définie de manière abstraite \gui{avec points}.

On peut aussi noter $\jdim\,\gA$ pour la dimension de $\jspec\,\gA$, 
qui n'est
pas un espace spectral (et nous ne proposons pas de \dfn \cov sans 
point pour
cette dimension).

\medskip \rem
Précisons la signification de $\Jdim\,\gA\leq \ell$ dans le cas des 
anneaux commutatifs. Comme il s'agit de la dimension de Krull de 
$\Heit\,\gA$ et que les \elts de $\Heit\,\gA$  s'identifient aux 
radicaux de Jacobson d'\itfs on obtient la caractérisation 
suivante:\\ 
$\forall x_0,\ldots,x_\ell\in \gA\;$ 
$\Ex \fa_0,\ldots,  \fa_\ell,$ \itfs de $\gA$ tels que: 
$$\begin{array}{rcl} 
x_0\,\fa_0  &  \subseteq  &  \JA(0)    \\ 
x_1\,\fa_1 &  \subseteq  &  \JA(\gen{x_0}+\fa_0)   \\
\vdots \;  &  \vdots &   \qquad \vdots   \\
x_{\ell}\,\fa_{\ell} &  \subseteq  &  
\JA(\gen{x_{\ell-1}}+\fa_{\ell-1})   \\
\gen{1} &  =  &  \JA(\gen{x_{\ell}}+\fa_{\ell})   
\end{array}$$
On ne peut apparemment pas éviter le recours aux \itfs et cela fait 
que l'on n'obtient pas une \dfn \gui{au premier ordre}.\\
Notez que chaque appartenance $x\in\JA(y_1,\ldots ,y_m)$ s'exprime 
elle-même par: $\Tt z\in\gA,$ $1+xz$ est inversible modulo
$\gen{y_1,\ldots ,y_m}$, \cade 

\medskip  \hspace*{4em}
$\Tt z\in\gA \;\Ex t,u_1,\ldots ,u_m\in\gA, 
\;\;1=(1+xz)t+u_1y_1+\cdots 
+u_my_m.$

\subsubsection*{Dimension  et  bord de Heitmann}

\begin{definition}
\label{defHdimA}
La dimension de Heitmann d'un anneau commutatif
est la dimension de Heitmann de son treillis de Zariski.
\end{definition}

\begin{definition}
\label{defHei2} Soit  $\gA$  un anneau commutatif, $x\in\gA$ et $\fj$ 
un \itfz.
Le \textsl{bord de Heitmann de $\fj$ dans $\gA$} est l'anneau quotient
$\gA/\rH_\gA(\fj)$ avec
          $$\rH_\gA(\fj):=\fj+(\JA(0):\fj)$$
qui est aussi appelé \textsl{l'\id bord de Heitmann de $\fj$ dans 
$\gA$}.
On notera aussi  $\rH_\gA(y_1,\ldots ,y_n)$ pour  
$\rH_\gA(\gen{y_1,\ldots
,y_n})$, $\rH_\gA^x$ pour $\rH_\gA(x)$ et $\gA_\rH^x$ pour 
$\gA/\rH_\gA^x$.
\end{definition}

Ainsi un \elt arbitraire de  $\rH_\gA(y_1,\ldots ,y_n)$ s'écrit 
$\sum_i
a_iy_i+b$ avec tous les $by_i$ dans $\JA(0)$.

La proposition suivante résulte des bonnes propriétés de la 
correspondance
bijective $\IZA$ (voir les faits \ref{factSpecAzarA}, \ref{factQuoAT}, 
\ref{factTransporteurs}
et \ref{factRadJac}).
\begin{proposition}
\label{propBordH-TA}
Pour un \itf  $\fj$ le bord de Heitmann de $\fj$ au sens des anneaux 
commutatifs
et celui au sens  des \trdis se correspondent. Plus précisément,
avec $j=\DA(\fj)$ et $\gT=\ZarA$, on a: 
$$\IZA(\rH_\gA(\fj))=\rH_\gT(j),\;\;
\emph{et}\;\; \Zar(\gA/\rH_\gA(\fj))\simeq
\gT/(\rH_\gT(j)=0)=\gT_\rH^j.$$
\end{proposition}

Comme corolaire des propositions \ref{propHdimgen} et 
\ref{propBordH-TA}
on obtient le résultat suivant.

\begin{proposition}
\label{lemDHA} Pour un anneau commutatif  $\gA$ et un entier 
$\ell\geq 0$
\propeq
\begin{enumerate}
\item La dimension de Heitmann de $\gA$ est $\leq \ell$.
\item Pour tout $x\in\gA$,
$\Hdim(\gA/\rH_\gA(x))\leq \ell-1$.
\item Pour tout \itf $\fj$ de $\gA$, $\Hdim(\gA/\rH_\gA(\fj)) \leq 
\ell-1$.
\end{enumerate}
\end{proposition}

\rem La dimension de Heitmann de $\gA$ peut donc être définie 
de manière
inductive comme suit:
\begin{itemize}
\item $\Hdim\,\gA=-1$ \ssi $1_\gA=0_\gA$.
\item Pour $\ell\geq 0$, $\Hdim\,\gA\leq \ell$ \ssi pour tout 
$x\in\gA$,
$\Hdim(\gA/\rH_\gA(x))\leq \ell-1$.
\end{itemize}
On  peut illustrer cette \dfnz.
Nous introduisons \gui{l'idéal bord de Heitmann itéré}.
Pour $x_1,\ldots ,x_n\in\gA$ nous notons $\gA_\rH[x_1]=\gA_\rH^{x_1}$,
$\gA_\rH[x_1,x_2]=(\gA_\rH^{x_1})_\rH^{x_2}$,
$\gA_\rH[x_1,x_2,x_3]=((\gA_\rH^{x_1})_\rH^{x_2})_\rH^{x_3}$, 
etc\ldots\,\,  les
anneaux bords de Heitmann successifs, et $\rH[\gA;x_1,\ldots
,x_k]=\rH_\gA[x_1,\ldots ,x_k]$ désigne le noyau de la projection 
canonique
$\gA\to \gA_\rH[x_1,\ldots ,x_k]$.  Pour décrire ces \ids nous 
avons besoin de
la notation $$\bidule{z,x,a,y,b}=1+(1+(z+ax)xy)b.$$
Alors on a:
\begin{itemize}
\item $z\in\rH_\gA[x_0]$ \ssi: $$\Ex a_0  \ \Tt y_0 \ \Ex
b_0,\;\bidule{z,x_0,a_0,y_0,b_0}=0$$
\item  $z\in\rH_\gA[x_0,x_1]$ \ssi: $$\Ex a_1  \ \Tt y_1 \ \Ex b_1\ 
\Ex a_0  \
\Tt y_0 \ \Ex b_0, 
\;\bidule{\bidule{z,x_1,a_1,y_1,b_1},x_0,a_0,y_0,b_0}=0$$
\item  $z\in\rH_\gA[x_0,x_1,x_2]$ \ssi:
$$\Ex a_2  \ \Tt y_2 \ \Ex b_2 \ \Ex a_1  \ \Tt y_1 \ \Ex b_1 \ \Ex 
a_0  \ \Tt
y_0 \ \Ex
b_0,\;\bidule{\bidule{\bidule{z,x_2,a_2,y_2,b_2},x_1,a_1,y_1,b_1},x_0,a_0,y_0,
b_0}=0$$
\end{itemize}
Et ainsi de suite. Et la dimension de Heitmann est $\leq \ell$ \ssi 
pour tous
$x_0,\ldots ,x_\ell\in\gA$ on a $1\in\rH_\gA[x_0,\ldots ,x_\ell]$.

\begin{proposition}
\label{propHei2} Soit $\fj=\gen{j_1,\ldots ,j_n}$ un \itfz.  Notons
$\varphi=\JA(\fj)=\JA(j_1)\vu\cdots \vu\JA(j_n)$.  Alors
$\Heit(\gA/\rH_\gA(\fj))$ s'identifie naturellement avec un quotient
de $(\HeA)\ul{\varphi}$.  Il y a \egt lorsque $\HeA$ est une \agHz,
donc en particulier lorsque $\Jspec\,\gA$ est noethérien.
\end{proposition}
\begin{proof}
Déjà démontré pour un \trdi arbitraire à la place de $\ZarA$
(proposition~\ref{propBHeitHeyt}).
\end{proof}

\rem On a aussi (déjà démontré pour les \trdisz) les 
résultats
suivants:
\begin{itemize}
\item on a toujours  $\Hdim\,\gA\leq \Jdim\,\gA\leq 
\Kdim(\gA/\JA(0))$,
\item si $(\fa_i)_{1\leq i\leq m}$ une famille finie d'\ids de $\gA$ 
et
$\fa=\bigcap_{i=1}^m\fa_i$, alors
$\Hdim(\gA/\fa)=\sup_i\Hdim(\gA/\fa_i)$.
\item si $\HeA$ est une \agH (en particulier si  $\Jspec\,\gA$ est 
noethérien)
on a $\Hdim\,\gA=\Jdim\,\gA$,
\item si $\Max\,\gA$ est noethérien, alors  
$\jspec\,\gA=\Jspec\,\gA$  et
$\Hdim\,\gA=\Jdim\,\gA=\jdim\,\gA$.
\item $\Hdim\,\gA\leq 0\;\Leftrightarrow\;\Jdim\,\gA\leq
0\;\Leftrightarrow\;\Kdim(\gA/\JA(0))\leq 0$.
\end{itemize}

\medskip Notez que le treillis  $\Heit\,\gA$ est une \agH \ssi est
vérifiée la propriété suivante:
$$\Tt \fa,\fb\in\Heit\,\gA\;\Ex\fc\in\Heit\,\gA\;(\fc\fb\subseteq\fa\;\mathrm{et}\;
\Tt x\in\gA\;(x\fb\subseteq\fa\Rightarrow x\in\fc))$$
($\fa=\JA(a_1,\ldots ,a_n)$, $\fb=\JA(b_1,\ldots ,b_m)$, 
$\fc=\JA(c_1,\ldots
,c_\ell)$).

\section[Le \tho de Kronecker et le stable range de Bass]{Le \tho de 
Kronecker
et le stable range de Bass (versions non noethériennes de Heitmann)}
\label{secKroBass}

Dans cette section et la suivante on établit la version 
\elrz,  non
noethérienne et \cov des \thos cités dans les titres. Les démonstrations en \clama  
de ces \thos sont pour l'essentiel données dans \cite{Hei84}.

Pour le \tho de Kronecker en particulier seule la dimension de
Krull intervient, et le \tho se trouve donc dans \cite{Hei84} (voir 
aussi
\cite{Hei76}). 

Pour les autres \thos des sections \ref{secKroBass} et \ref{secSFS}, 
notre \gui{dimension de Heitmann} peut être dans certains cas 
strictement 
plus petite que la
dimension utilisée par Heitmann (celle de son J-spectrum), et donc 
notre
version est dans ce cas \gui{meilleure}. 

Enfin concernant le \tho de Swan et toute la section \ref{secSwan}, 
nos  versions sont \covs et non noethériennes. Nous démontrons les résultats
souhaités pour la dimension de Heitmann, donc aussi pour la \(\Jdim\).
Nous répondons ainsi à une conjecture de Heitmann, qui n'avait le résultat que pour une variante moins forte  (il utilisait une dimension appelée
$\delta$-dim).

\smallskip Les preuves données dans les sections suivantes sont
essentiellement celles de \cite{Coq2004} (pour le \tho de Kronecker)
et de \cite{CLQ2004}, quelquefois améliorées en tenant compte de
\cite{Duc2006}. Nous avons rajouté en outre quelques précisions, ainsi 
que les
résultats de la section~\ref{subsecBKH}.

\smallskip Signalons enfin que le contenu des sections 5, 6 et 7 est traité en détail dans le chapitre XIV du livre \cite{ACMC}.

\subsection{Le \tho de Kronecker}
\label{subsecKro}

Ce \tho a d'abord été démontré par Kronecker sous la forme 
suivante  (\cite{Kro82}): une
variété algébrique dans $\CC^n$ peut toujours être 
définie par $n+1$
équations.

Il a été étendu au cas des anneaux noethériens (par van der 
Waerden,
dans \cite{vW41})
sous la forme suivante: dans un anneau noethérien de dimension de 
Krull $n$,
tout \id  a même radical qu'un \id engendré par au plus $n+1$ 
\eltsz.

La version de Kronecker a été améliorée par divers auteurs  
dans
\cite{EG73a,Stor72} qui ont montré que $n$ équations suffisent en 
général. On
peut trouver une preuve \cov de ce \tho dans \cite{CLS2005}.
Par ailleurs on ne sait toujours pas si toute courbe dans l'espace 
complexe de
dimension 3 est ou non intersection de deux surfaces.

Enfin, Heitmann a généralisé la version de van der 
Waerden au cas non noethérien (\cite{Hei84}).

Le lemme suivant, bien que terriblement anodin, est une clef 
essentielle.
\begin{lemma}
\label{gcd}
Pour $u,v\in\gA$ on a
$$\begin{array}{rcccl}
~~~~~~~~~~~~~~~\DA(u,v)& =  & \DA(u+v,uv)  & =  & \DA(u+v)\vu\DA(uv) 
\quad
\mathrm{et\;donc\;aussi} \\[1mm]
\JA(u,v)& =  & \JA(u+v,uv)  & =  & \JA(u+v)\vu\JA(uv)
\end{array}$$
En particulier
\begin{itemize}
\item  Si $uv\in\DA(0)$, 
alors $\DA(u,v)=\DA(u+v)$
\item  Si $uv\in\JA(0)$, 
alors $\JA(u,v)=\JA(u+v)$
\end{itemize}
\end{lemma}

Rappelons que deux suites qui vérifient les inégalités 
(\ref{eqCG})
dans le corolaire \ref{corthDKA} sont dites complémentaires.
\begin{lemma}
\label{lemKroH} 
Si $(b_1,\ldots ,b_n)$ et  $(x_1\ldots ,x_n)$ sont deux suites 
complémentaires dans $\gA$
alors pour tout $a\in\gA$ on a: 
$$\DA(a,b_1,\dots,b_n) = \DA(b_1+ax_1,\dots,b_n+ax_n),$$
\cad encore: $a\in \DA(b_1+ax_1,\dots,b_n+ax_n)$.
\end{lemma}
\begin{proof}
Pour $n=0$ l'anneau est trivial et $\DA(a)=\DA(\emptyset)$.\\
Faisons aussi le cas $n=1$ bien que ce soit inutile 
(vue la preuve par \recuz). 
On a $b_1x_1\in\DA(0)$ et
$1\in\DA(b_1,x_1)$. A fortiori  $b_1ax_1\in\DA(0)$ et 
$a\in\DA(b_1,ax_1)$. 
Et le lemme \ref{gcd} nous dit que $\DA(b_1,ax_1)=\DA(b_1+ax_1)$.\\
Donnons maintenant la preuve de la \recuz.
On considère l'anneau $\gB=\gA/\DA(b_1,x_1)$. Par \hdr
on a  pour tout $a$,  $$a\in \rD_\gB(b_2+ax_2,\dots,b_n+ax_n).$$ 
Cela signifie dans $\gA$
$$a\in\DA(b_1,x_1,b_2+ax_2,\dots,b_n+ax_n).$$ 
On en déduit  $a\in\DA(b_1,ax_1,b_2+ax_2,\dots,b_n+ax_n)$
(car $a\in\DA(\fJ)$ implique $a\in\DA(a\fJ)$). 
Mais $\DA(b_1x_1)=0$ implique $\DA(b_1ax_1)=0$ et 
le lemme \ref{gcd} nous dit que
$\DA(b_1,ax_1)=\DA(b_1+ax_1)$.
\end{proof}

\begin{theorem}
\label{thKroH} \emph{(de Kronecker, avec la dimension de Krull, sans
noethérianité)}\\
Soit $n\geq 0$. Si $\Kdim\,\gA <n$ et $b_1,\dots,b_n\in\gA$ alors il 
existe
$x_1,\dots,x_n$ tels que pour tout $a\in\gA$
$\DA(a,b_1,\dots,b_n) = \DA(b_1+ax_1,\dots,b_n+ax_n)$.
\\
En conséquence, dans un anneau de dimension de Krull $\leq n$, tout 
\itf a
même radical qu'un \id engendré par au plus $n+1$ \eltsz.
\end{theorem}

\begin{proof}
Le premier point est clair d'après le lemme
\ref{lemKroH} et le corolaire \ref{corthDKA}. La deuxième 
affirmation en découle car il suffit
d'itérer le processus. En fait, si $\Kdim\,\gA\leq n$ et 
$\fJ=\DA(b_1,\ldots
,b_{n+r})$  ($r\geq 2$) on obtient en fin de compte   
$\fJ=\DA(b_1+c_1,\ldots
,b_{n+1}+c_{n+1})$ avec les $c_i\in\gen{b_{n+2},\ldots ,b_{n+r} }$.
\end{proof}

\subsection{Le \tho \texorpdfstring{\gui{stable range}}{"stable 
range"} \,de
Bass}
\label{subsecBass}

Le \tho suivant est d\^u à Bass dans le cas noethérien avec la 
dimension de Krull. La version non noethérienne avec la $\Jdim$ est due à Heitmann. Nous donnons ici la version non noethérienne avec la $\Hdim$, à priori plus générale.

Notez que \textsl{la version non noethérienne avec la dimension de 
Krull résulte directement du premier point dans le \thoz~\ref{thKroH}}.

\begin{theorem}
\label{Bass} \emph{(Théorème de Bass, avec la dimension de Heitmann, sans
noethérianité)}\\
Soit $n\geq 0$. Si $\Hdim\,\gA <n$ et $1 \in\gen{a,b_1,\dots,b_n } $ 
alors il
existe $x_1,\dots,x_n$ tels que $1  \in\gen{b_1+ax_1,\dots,b_n+ax_n}$.
\end{theorem}

\begin{proof}
On remarque que  pour toute liste $L$
$$1 \in\gen{L}\;\Leftrightarrow\; 1\in\DA(L)\;\Leftrightarrow\; 
1\in\JA(L).$$ 
La preuve qui suit peut être considérée comme
une répétition de celle du \tho \ref{thKroH} lorsque $a=1$ 
en remplaçant le bord de Krull par le bord de Heitmann.
  La preuve est par \recu sur $n$. \\
  Lorsque $n=0$ l'anneau est trivial et $\DA(1)=\DA(\emptyset)$.\\
Supposons  $n\geq 1$. Soit $\fj=\rH_\gA(b_n)$ l'\id bord de Heitmann 
de $b_n$.
Puisque $\gA/\fj$ est de dimension de Heitmann $\leq n-2$ l'\hdr nous 
donne
$x_1,\ldots ,x_{n-1}\in\gA$ tels que
$$
1\in \rD(b_1+x_1a,\ldots ,b_{n-1}+ x_{n-1} a) \quad \mathrm{dans} 
\quad
\gA/\fj\,.
$$
Notons  $L$ pour $b_1+x_1a,\ldots ,b_{n-1}+ x_{n-1} a$. Par \dfn du 
bord
de Heitmann, l'\egt ci-dessus implique
qu'il existe $x_n$ tel que
$$
x_nb_n\in\JA(0)\quad  \mathrm{et}\quad
1\in\DA(L,b_n,x_n)=\DA(L,b_n)\vu\DA(x_n).
$$
Puisque
$$
1\in\DA(L,b_n,a)=\DA(L,b_n)\vu\DA(a),
$$
cela implique par distributivité
$$
1\in \DA(L,b_n,x_na)=\JA(L,b_n,x_na).
$$
Puisque $b_n\,x_n\,a\in\JA(0)$  le lemme \ref{gcd} nous dit que
$\JA(b_n,x_na)=\JA(b_n+x_na)$, et donc que
$$\JA(L,b_n+x_na)=\JA(L,b_n,x_na)=\gen{1},
$$
ce qui était le but recherché.
\end{proof}

Rappelons qu'un vecteur $L$ de $\gA^m$ est dit \textsl{unimodulaire} 
lorsque ses
coordonnées engendrent l'\id $\gen{1}$, \cade lorsque $\DA(L)=1$.
\begin{corollary}
\label{corBass}
Soit $n\geq 0$. Si $\Hdim\,\gA \leq  n$ et $V\in\gA^{n+2}$ est 
unimodulaire, il
peut être transformé en le vecteur $(1,0\ldots ,0)$
par des manipulations \elrsz.
\end{corollary}

Il s'ensuit alors directement qu'un module stablement libre de rang 
$\geq n+1$
sur un anneau dont la dimension de Heitmann est $\leq  n$ est libre 
(comme dans
\cite{Lam1978}, p. 28), sans aucune hypothèse noethérienne.

\subsection{Une généralisation de Heitmann}
\label{subsecBKH}

On traite ici \cot le corolaire 2.2 de  \cite{Hei84}, qui
généralise à la fois les résultats de Bass et Kronecker. Cela 
débouche sur une amélioration du \tho de Kronecker. Dans \cite{Hei84} ces résultats sont établis pour la $\Jdim$.

\goodbreak
\begin{lemma} \label{HenriLemma} ~
Soient $a, c_1,\ldots,c_m, u,v,w \in \gA$. Notons $Z=c_1,\ldots,c_m$.
\begin{enumerate}
\item
On a $a \in \DA(Z)  \iff  1 \in \gen {Z}_{\gA[a^{-1}]}$.
\item 
$w\in \JA(\gA[a^{-1}],0)$ et $a\in \DA(Z,w)$
$\Rightarrow$ $a\in \DA(Z)$.
\item 
$uv\in \JA(\gA[a^{-1}],0)$ et $a\in \DA(Z,u,v)$ $\Rightarrow$ $a\in \DA(Z,u+v).$ 
\end{enumerate}
\end{lemma}

\begin{proof}
\textsl{1.} Immédiat.

\noindent\textsl{2.} Supposons $a\in \DA(Z,w)$ et travaillons dans l'anneau
$\gA[a^{-1}]$. On a $1 \in \gen {Z}_{\gA[a^{-1}]} +
\gen{w}_{\gA[a^{-1}]}$ et comme $w$ est dans 
$\JA(\gA[a^{-1}],0)$, cela implique $1 \in \gen {Z}_{\gA[a^{-1}]}$, i.e.
$a\in \DA(Z)$.

\noindent\textsl{3.}
Résulte du point \textsl{2} car  $\DA(Z,u,v)=\DA(Z,u+v,uv)$ (lemme \ref{gcd}).
\end{proof}

\rem On peut se demander si l'idéal $\JA(\gA[a^{-1}],0)$ est le meilleur
possible. La réponse est oui. L'implication du point \textsl {2} est
vérifiée (pour tous $c_1,\dots,c_m$) en remplaçant $\JA(\gA[a^{-1}],0)$
par $\fJ$ \ssi $\fJ\subseteq \JA(\gA[a^{-1}],0)$.

\medskip  Voici maintenant le corolaire 2.2 de  \cite{Hei84}, avec la 
$\Hdim$ qui
remplace la $\Jdim$. 

Si $L=(b_1,\dots,b_n)\in \gA^n$  nous écrivons 
$\DA(L)$ pour
$\DA(b_1,\dots,b_n)$.

\begin{lemma}
\label{thCor2.2Heit}
Si $\Hdim(\gA [a^{-1}])< n$, $L\in \gA^n$ et $\DA(b)\leq\DA(a)\leq \DA(b, L)$
alors il existe $X\in \gA^n$ tel que $\DA(L+bX)=\DA(b, L)$, ce qui équivaut
à $b \in\DA(L+bX)$ ou encore à $a \in\DA(L+bX)$.
En outre nous pouvons prendre $X=aY$ avec $Y\in \gA^n$.
\end{lemma}
\begin{proof}~\\
\textsl{Remarque préliminaire.}
Si $\DA(L+bX)=\DA(b, L)$, on a évidemment $b\in\DA(L+bX)$ mais aussi
$a\in\DA(L+bX)$ car $a\in \DA(b,L)$. Réciproquement, si $b\in\DA(L+bX)$, on
a $\DA(L+bX)=\DA(b, L)$; a fortiori, si $a\in\DA(L+bX)$, on a $b \in
\DA(L+bX)$ (puisque $b\in\DA(a)$), donc aussi $\DA(L+bX)=\DA(b, L)$.

\sni
On raisonne par \recu sur $n$, le cas $n=0$ étant trivial.  On commence par
chercher $X\in \gA^n$.  Soient $\fj=\rH_{\gA[a^{-1}]}(b_n)$ et
$\gA'=\gA/(\fj\cap\gA)$ où $\fj\,\cap\,\gA$ est mis pour
\gui{l'image réciproque de $\fj$ dans $\gA$}.
On a une identification $\gA[a^{-1}]/\fj=\gA'[a^{-1}]$.
Comme $\Hdim(\gA'[a^{-1}])< n-1$, on peut appliquer l'hypothèse de \recu à
$\gA'$ et $a,b,b_1,\ldots ,b_{n-1}$ (en remarquant que $b_n = 0$ dans $\gA'$);
on obtient $x_1,\dots,x_{n-1}\in \gA$ tels que $\rD(b_1+bx_1,\dots,
b_{n-1}+bx_{n-1})=
\rD(b, b_1,\dots,b_{n-1})$ dans $\gA'$. D'après la
remarque préliminaire, cette \egt équivaut, 
en notant $Z$ pour $(b_1+bx_1,\ldots, b_{n-1}+ bx_{n-1})$, à
$$
a \in \rD(Z) \quad \mathrm{dans} \quad\gA'\,.
$$
L'appartenance ci-dessus signifie $1 \in \gen{Z}$ dans $\gA'[a^{-1}]$,
i.e. $1 \in \gen{Z} + \fj$ dans $\gA[a^{-1}]$. Par
\dfn du bord de Heitmann, cela veut dire qu'il existe $x_n$ que l'on peut
choisir dans $\gA$, tel que $x_nb_n\in\JA(\gA[a^{-1}],0)$ et $1 \in \gen{Z,
b_n,x_n}_{\gA[a^{-1}]}$.  On a donc $a\in\DA(Z,b_n,x_n)$. Mais on a aussi $a
\in\DA(Z,b_n,b)$, puisque

\snic {
\gen {Z, b_n, b} = \gen {b_1,, \ldots, b_{n-1}, b_n, b} 
\eqdefi {\gen{L, b}}
}

\sni et que $a \in \DA(L,b)$ par hypothèse.
Bilan: $a \in\DA(Z,b_n,x_n)$, $a \in\DA(Z,b_n,b)$ donc $a \in\DA(Z,b_n,bx_n)$.
L'application du lemme~\ref{HenriLemma}~\textsl{3}  avec
$u = b_n$, $v = bx_n$ fournit $a \in\DA(Z,b_n+bx_n)$, i.e.  $a \in \DA(L +
bX)$ où $X = (x_1, \ldots, x_n)$.
\\
Enfin, si $b^p\in\gen{a}_\gA$, nous pouvons appliquer le résultat avec $b^{p+1}$ à la place de $b$
puisque $\DA(b)=\DA(b^{p+1})$. Alors $L+b^{p+1}X$ se réécrit $L+baY$.
\end{proof}

Pour $a\in\gA$ on a toujours  $\Hdim(\gA[a^{-1}])\leq 
\Kdim(\gA[a^{-1}])\leq
\Kdim(\gA)$. Par conséquent le \tho suivant améliore le \tho de 
Kronecker.

\begin{theorem}
\label{KroH2} \emph{(de Kronecker, variante dimension de Heitmann)}
\begin{enumerate}
\item Soit $n\geq 0$. Si $\Hdim(\gA[a^{-1}]) <n$ et $a,b_1,\dots,b_n\in\gA$
alors il existe $x_1,\dots,x_n\in \gA$ tels que
$\DA(a,b_1,\dots,b_n) =\DA(b_1+ax_1,\dots,b_n+ax_n)$.
\item En conséquence, si $a_1,\ldots ,a_r,b_1,\dots,b_n\in\gA$ et
$\Hdim(\gA[a_i^{-1}]) <n$ pour $1\leq i\leq r$ alors  il existe $y_1,\ldots,y_n\in
\gen{a_1,\ldots,a_r}$ tels que

\snic{\DA(a_1,\ldots ,a_r,b_1,\dots,b_n) =\DA(b_1+y_1,\dots,b_n+y_n).}
\end{enumerate}
\end{theorem}
\begin{proof}
\textsl{1.} Conséquence directe du lemme 
\ref{thCor2.2Heit} en
faisant $a=b$.

\noindent \textsl{2.} Se déduit de \textsl{1} par \recu sur $r$: 
$$\fJ=\DA(a_1,\ldots
,a_r,b_1,\dots,b_n) = \DA(a_1,\ldots ,a_{r-1},b_1,\dots,b_n) \vu
\DA(a_r)=\fII\vu\DA(a_r)$$ 
$\fII= \DA(b_1+z_1,\dots,b_n+z_n)$ avec
$z_1,\dots,z_n\in \gen{a_1,\ldots ,a_{r-1}}$, donc
$\fJ=\DA(a_r,b_1+z_1,\dots,b_n+z_n)$ et on applique une nouvelle fois 
le
résultat.
\end{proof}

\rem Dans le cas d'un anneau local noethérien $\gA$, l'\idema est 
radicalement
engendré par $n=\Kdim\,\gA$ \eltsz, mais pas moins. Ceci montre que 
parmi $n$
\gui{\gtrs radicaux} du maximal, il n'y en a
aucun  vérifiant $\Hdim(\gA[x^{-1}])< n-1$. 
Par ailleurs \hbox{$\Hdim(\gA)=0$} puisque $\gA$ est local.

\section[Le splitting off de Serre et le \tho de Forster]{Le 
splitting off de
Serre et le \tho de Forster, à la Heitmann}
\label{secSFS}

Dans cette section nous donnons une version \cov du travail de 
Heitmann (\cite{Hei84}) concernant le splitting off de Serre, le \tho 
de Forster et une variation sur le \tho de Swan.
Nous remplaçons la $\Jdim$ par la $\Hdim$, qui est plus facile 
à manipuler et donne des résultats à priori plus précis.

Dans la section suivante, nous donnerons des \thos qui sont  
\gui{meilleurs} dans la mesure où ils sont plus faciles à mettre 
en {\oe}uvre (l'hypothèse est moins lourde car elle ne fait pas 
intervenir la dimension de nombreux localisés) et où la dimension 
de Heitmann est mal controlée lors d'une localisation.
\subsection{Manipulations \elrs de colonnes}

\begin{lemma}
\label{main}
Si $a\in \gA$, $\Hdim(\gA[a^{-1}])< n$ et $L\in \gA^n$ il existe $X\in
\gA^n$ tel que $a\in \DA(L-aX)$.  En outre nous pouvons prendre $X=aY$
avec $Y\in \gA^n$.
\end{lemma}

\begin{proof}
Cas particulier du lemme \ref{thCor2.2Heit}: faire $a=b$.
\end{proof}

Pour $n=1$ cela donne: \textsl{si $\Hdim(\gA[a^{-1}])\leq 0$ alors pour tout 
$b$ on peut
trouver $x$ tel \hbox{que $a\in\DA(b+xa)$}}.

\begin{corollary}
\label{cormain}
Soit $M$ une matrice dans $\gA^{n\times n}$ et $\delta$ son
déterminant.  Si $\Hdim(\gA[\delta^{-1}]) < n$ alors pour tout $C\in
\gA^n$ il existe $X\in \gA^n$ tel que $\delta\in \DA(MX-C)$.  En outre
nous pouvons trouver $X$ de la forme $\delta Y,~Y\in \gA^n$.
\end{corollary}

\begin{proof}
La preuve est basée sur les formules de Cramer.  Soit $\wi{M}$ la
matrice adjointe de $M$, et $L = \wi{M}C$.  Nous avons $\widetilde
{M}(MX -C) = \delta X - L$ pour un vecteur colonne arbitraire $X\in
\gA^n$.  Donc l'\id engendré par les coordonnées de $\delta X-L$
est inclus dans celui engendré par les coordonnées de $MX - C$, et
$$
\DA( \delta X - L) \leq \DA(MX -C)
$$
D'après le lemme \ref{main} nous pouvons trouver {\sl un} $X\in
\gA^n$ tel que $\delta\in \DA(\delta X-L)$, et donc $\delta\in
\DA(MX-C)$ comme demandé.
\end{proof}

Pour la fin de cette section nous fixons les notations suivantes:
\begin{notation}
\label{notaMatrix}
{\rm Soit $F$  une matrice dans $\gA^{n\times p}$ avec pour colonnes
$C_0,C_1,\dots,C_p$, et
$G$ la matrice ayant pour colonnes $C_1,C_2,\dots,C_p$.
Soit $\Delta_k(F)$ l'\id déterminantiel d'ordre $k$ de $F$ (l'\id 
engendré
par les mineurs $\nu$ d'ordre $k$ de $F$).
}
\end{notation}

\begin{theorem}
\label{matrix}
Fixons $0< k\leq  p$. Supposons que pour chaque mineur $\nu$ d'ordre 
$k$ de $G$
l'anneau $\gA[\nu^{-1}]$ est de dimension de Heitmann $<k$. Alors il 
existe
$t_1,\dots,t_p$ tels que
$$\DA(C_0,\Delta_k(F))\leq \DA(C_0+t_1C_1+\cdots+t_pC_p),$$
(i.e., $\Delta_k(F)\subseteq \DA(C_0+t_1C_1+\cdots+t_pC_p)$ et 
$\DA(C_0)\leq
\DA(C_0+t_1C_1+\cdots+t_pC_p)$).
\end{theorem}

\hum{L'hypothèse du \tho est-elle que  l'ouvert
$D_\gA(\Delta_k(G))\subseteq\Spec\,\gA$ est un espace spectral de 
dimension de
Heitmann $< k$. ??? C'est sans doute un peu plus compliqué que 
cela. Il semble
que si on raisonnait avec la $\Jdim$, il faudrait recoller plutot les
$\He(\Zar(\gA[\nu^{-1}]))$ que les $\Zar(\gA[\nu^{-1}])$. En outre 
pour la
$\Hdim$, cela pourrait ne correspondre à aucune dimension d'espace 
spectral.}

\begin{proof}
L'inclusion $\Delta_k(F)\subseteq 
\DA(C_0+t_1C_1+\cdots+t_pC_p)=\DA(C'_0)$
signifie que pour tout mineur $\nu$ d'ordre $k$ de $F$ on a 
$\nu\in\DA(C'_0)$.
Soit $\Delta_k(G)$ l'\id déterminantiel d'ordre $k$ de $G$.
Remarquons que si le mineur $\nu$ fait intervenir la colonne $C_0$ et 
les
colonnes $C_{i_j}$ ($1\leq i_1 <\cdots<i_{k-1}\leq p)$ alors on peut 
y remplacer
la colonne $C_0$  par la colonne $C'_0=C_0+t_1C_1+\cdots+t_pC_p$  
sans modifier
sa valeur autrement qu'en ajoutant une combinaison linéaire de 
mineurs d'ordre
$k$ de $G$, de sorte que $\nu\in \DA(C'_0,\Delta_k(G))$.
En conclusion il nous suffit de réaliser
$\Delta_k(G)\subseteq \DA(C'_0)$ pour avoir $\Delta_k(F)\subseteq 
\DA(C'_0)$.

\noindent En fait il suffit de savoir réaliser
$$ \DA(C_0,\nu_1)\leq \DA(C_0+t_1C_1+\cdots+t_pC_p)=\DA(C'_0)$$
   pour \textsl{un}
mineur $\nu_1$ d'ordre $k$ de $G$.
Car alors  nous remplaçons $C_0$ par $C'_0$ dans $F$  (ce qui ne 
change pas
$G$) et nous pouvons passer à un autre mineur $\nu_2$ de $G$ pour 
lequel nous
obtiendrons  $t'_1,\dots,t'_p$ avec
$$\DA(C_0,\nu_1,\nu_2)\leq\DA(C'_0,\nu_2)\leq 
\DA(C'_0+t'_1C_1+\cdots+t'_pC_p)=
\DA(C''_0)$$
avec $C''_0=C_0+t''_1C_1+\cdots+t''_pC_p$ et ainsi de suite.\\
Pour réaliser
$$\nu\in \DA(C_0+t_1C_1+\cdots+t_pC_p)\quad \mathrm{et}\quad 
\DA(C_0)\leq
\DA(C_0+t_1C_1+\cdots+t_pC_p)$$
   pour un
mineur $\nu$ d'ordre $k$ de $G$ nous utilisons le corolaire 
\ref{cormain}:
nous trouvons $t_1,\dots,t_p\in\nu\gA$ (avec $t_i = 0$ pour les 
colonnes qui
n'interviennent pas dans le mineur $\nu$)
  tels que
$$
\nu\in \DA(C_0+t_1C_1+\cdots+t_pC_p)=\DA(C'_0)
$$
Puisque $t_1,\dots,t_p$ sont des multiples de $\nu$ nous avons aussi
$\DA(C_0)\leq \DA(\nu,C'_0)=\DA(C'_0)$
comme demandé.
\end{proof}

Toujours avec les notations \ref{notaMatrix}, nous obtenons le 
corolaire
suivant.
\begin{theorem}
\label{basic}
Supposons que $1 \in \Delta_1(F)$ et que pour   $k=1,\ldots ,p$ et 
chaque mineur
$\nu$ d'ordre $k$ de $G$ l'anneau
$\gA[\nu^{-1}]/\Delta_{k+1}(F)$ est de dimension de Heitmann $< k$. 
Alors il
existe $t_1,\dots,t_p$ tels que le vecteur $C_0+t_1C_1+\cdots+t_pC_p$ 
est
unimodulaire. \\
En particulier si $k\leq p$, $1 \in \Delta_k(F)$ et 
$\Hdim(\gA[\nu^{-1}])<k$ pour tout
$\nu\in\gA$,  il existe $t_1,\dots,t_p$ tels que le vecteur
$C_0+t_1C_1+\cdots+t_pC_p$ est unimodulaire.
\end{theorem}

\hum{L'hypothèse du \tho \ref{basic}
est-elle que  pour chaque $k$ la dimension de
Heitmann du constructible  $V_\gA(\Delta_{k+1}(F))\setminus 
D_\gA(\Delta_k(G))$
est  $< k$. ??? C'est sans doute un peu plus compliqué que cela}

\begin{proof}
En utilisant le \tho \ref{matrix}, nous définissons une suite de 
vecteurs
$C_0^k,~k=1,\ldots ,p$ avec
$C_0^1 = C_0$. Pour $k=1,\ldots ,p$, nous  raisonnons dans
$\gA_k=\gA/\Delta_{k+1}(F)$. Nous construisons $C_0^{k+1}$ de la forme
$C_0^k + u^k_1C_1+\cdots + u^k_pC_p$ tel que
$\rD_{\gA_k}(C_0^k,\Delta_k(F))\leq \rD_{\gA_k}(C_0^{k+1})$.
Cela signifie que nous avons, dans $\gA$
$$\DA(C_0^k,\Delta_k(F))\leq \DA(C_0^{k+1},\Delta_{k+1}(F))$$
D'où le résultat puisque $\Delta_1(F)= 1$ et $\Delta_{p+1}(F) = 
0$.
\end{proof}

Un cas particulier plus simple s'énonce avec la dimension de Krull.  Notez que
la matrice $G$ ne figure plus dans l'hypothèse.  La conclusion
resterait donc valable après avoir permuté des colonnes.

\begin{theorem}
\label{basicK}
Supposons que $1 = \Delta_1(F)$ et que pour  $k=1,\ldots ,p,$  
l'anneau
$\gA/\Delta_{k+1}(F)$ est de dimension de Krull $< k$. Alors il existe
$t_1,\dots,t_p$ tels que le vecteur $C_0+t_1C_1+\cdots+t_pC_p$ est
unimodulaire.\\
En particulier si $\Kdim \,\gA< k\leq p$ et $\Delta_{k}(F)=1$, il 
existe un vecteur
unimodulaire dans le module image de $F$.
\end{theorem}

Notez que la condition  $1 \in \Delta_1(F)$ peut se relire
\gui{$\gA/\Delta_{1}(F)$ est de dimension de Krull $< 0$}. On aurait 
donc pu
énoncer l'hypothèse du \tho sous la forme: pour chaque $k\geq 0$ 
l'anneau
$\gA/\Delta_{k+1}(F)$ est de dimension de Krull $< k$.

\junk{
\begin{corollary}
\label{corbasicK}
Supposons que $1 \in \Delta_1(F)$ et que pour chaque $k>0$ l'anneau
$\gA/(\Delta_{k+1}(F)+\JA(0))$ est de dimension de Krull $< k$. Alors 
il existe
$t_1,\dots,t_p$ tels que le vecteur $C_0+t_1C_1+\cdots+t_pC_p$ est
unimodulaire.\\
En particulier si $\Kdim(\gA/\JA(0))< k$ et $\Delta_{k}(F)=1$, il 
existe un
vecteur unimodulaire dans le module image de $F$.
\end{corollary}

\begin{proof}
Il suffit d'appliquer le \tho \ref{basic} avec l'anneau $\gA/\JA(0).$
On obtient avec cet anneau un vecteur unimodulaire dans le module 
image de $F$.
Mais si un vecteur $(a_1,\ldots ,a_s)$ est unimodulaire modulo 
$\JA(0)$ cela
signifie qu'on a $\sum a_ib_i\equiv 1\mod\;\JA(0)$
pour des $b_i$ convenables, et
tout \elt dans $1+\JA(0)$
est inversible.
\end{proof}
}

Nous donnons maintenant en \clama une reformulation \gui{avec 
\idepsz} du
\thoz~\ref{basic} dans lequel on remplace la dimension de Heitmann 
par 
la dimension de Krull. Dans l'énoncé qui suit, $\DA(\fJ)$ pour un 
\itf $\fJ$ est identifié à l'\oqc correspondant de $\Spec\,\gA$ et
 $\VA(\fJ)$ est son complémentaire. 

\begin{theorem}
\label{thbasic}
Pour tout \idep $\fp$ notons  $r_\fp$ le rang de la matrice $F$ vue 
dans le
corps des fractions de $\gA/\fp$ (on a $0\leq r_\fp\leq p+1$).
On donne deux formulations équivalentes pour les hypothèses du 
\thoz. 
\begin{enumerate}
\item Supposons que pour $k=0,\,\ldots ,\,p$ la dimension de Krull
du sous-espace spectral $\DA(\Delta_k(G))\cap \VA(\Delta_{k+1}(F))$ 
de $\Spec\,\gA$ soit $<k$ (pour $k=0$ cela signifie simplement 
$1\in\Delta_{1}(F)$).
\item Supposons que pour tout $\fp\in\DA(\Delta_{r_\fp-1}(G))$ on a
$\Kdim(A/\fp)<r_\fp.$
\end{enumerate}
Alors il existe un vecteur unimodulaire dans le module image de $F$.
En outre ce vecteur peut être écrit sous la forme donnée dans le
\thoz~\ref{basic}.
\end{theorem}
\begin{proof}
Fixons $k\leq p$ et posons $\gB_k=\gA/\Delta_{k+1}(F)$. Soit $\nu$ un 
mineur d'ordre $k$ de $G$. Les idéaux premiers de
$\gB_k[1/\nu]$ peuvent être identifiés aux idéaux premiers 
$\fp$ 
de $\gA$ qui contiennent $\Delta_{k+1}(F)$ mais ne contiennent pas 
$\nu$.
Donc le maximum des $\Kdim\,\gB_k[1/\nu]$ est la longueur
maximum d'une chaîne d'idéaux premiers qui contiennent 
$\Delta_{k+1}(F)$
mais évitent au moins un mineur d'ordre $k$ de $G$.
C'est donc la dimension de Krull
de $\DA(\Delta_k(G))\cap \VA(\Delta_{k+1}(F))$. Ainsi la 
première formulation de l'hypothèse nous 
ramène aux hypothèses du \thoz~\ref{basic}. \\
Soit $\fp$ un \idep de $\gA$. Dire que $r_\fp\leq k$ signifie que 
modulo $\fp$
tous les mineurs d'ordre $p+1$ sont nuls, \cad 
$\Delta_{k+1}(F)\subseteq\fp$. La deuxième
condition est donc bien équivalente à la première.
\end{proof}

\subsection{Quelques conséquences}

\subsubsection*{Le \tho splitting off de Serre, version non 
noethérienne}
\addcontentsline{toc}{subsubsection}{Le \tho splitting off de Serre}

On déduit directement du \tho \ref{basic}, comme dans 
\cite{EG73b,Hei84},
la version suivante \gui{améliorée} du
Splitting-off de Serre.

\begin{theorem}
\label{thSerre} \emph{(Théorème de Serre non noethérien à la Heitmann, 
version dimension de Krull)}\\
Soit $M$ un \Amo projectif de rang $\geq k$ sur un anneau $\gA$ tel
que $\Kdim\, \gA< k$. Alors  $M\simeq N\oplus \gA$ pour un certain $N$.
\end{theorem}

\begin{proof}
Si $F$ est une matrice de projection dont l'image est (isomorphe à) 
$M$,
l'hypothèse que $M$ est de rang $\geq k$ signifie que 
$\Delta_{k}(F)=\gen{1}$.
Le \tho~\ref{basic} nous donne dans l'image de $F$ un vecteur 
unimodulaire.
C'est un \elt $C$ de $M$ tel que $\gA C$ est facteur direct dans $M$.
\end{proof}

\rem Dans le fonctionnement de cette preuve, on voit qu'on aurait pu 
se 
contenter de supposer que $M$ est un module image d'une matrice $F$ 
telle que
$\Delta_{k}(F)=1$. 

\hum{En fait l'hypothèse peut être un petit peu affaiblie car il 
suffit que
les $\Hdim$ des  $\gA[\nu^{-1}]/\Delta_{\ell+1}$ pour les mineurs 
$\nu$ d'ordre
$\ell \geq k$ soient toutes $< k$. C'est un peu difficile à 
formuler. Cela
semble dépendre de la matrice mais ce n'est peut-être pas le 
cas.
Il semble qu'il doit y avoir aussi une formulation avec uniquement 
des mineurs
principaux.
Mais tout cela semble un peu du pinaillage.
Peut-être un  énoncé \gui{avec points} avec la $\Jdim$ serait 
plus parlant.}

\subsubsection*{Le \tho de Forster, versions non noethériennes}
\addcontentsline{toc}{subsubsection}{Le \tho de Forster}

Voici une version du \tho de Forster  (\cite{For64}).

\begin{theorem}
\label{forster}
Soit $M$ un module \tf sur un anneau $\gA$ et $\fJ$ l'annulateur de 
$M$.
Si la dimension de Krull de $\gA/\fJ$ est $\leq  d$
et si $M$ est localement engendré par $r$ \elts alors $M$ peut 
être engendré
par $d+r$ \eltsz.
\end{theorem}
\begin{proof}
$M$ est un quotient d'un module de présentation finie $M'$ dont 
l'idéal de
Fitting d'ordre $r$ égal à $1$,
et nous pouvons donc supposer $M$ \pfz.
On regarde $M$ comme un $\gA/\fJ$-module.
Soit $m_0,m_1,\dots,m_p$ un système de \gtrs de
$M$ et $F$ une matrice de présentation de $M$ correspondant à ce 
système
\gtrz. Dire que  $M$ est localement engendré par $r$ \elts 
signifie que $1 = \Delta_{p+1-r}(F)$. 
Si $p\geq d+r$ on a $1 = \Delta_{d+1}(F)$ et d'après le \tho
\ref{matrix}, appliqué à la matrice transposée de $F$,
nous pouvons trouver $t_1,\dots,t_p$ tels que $M$ soit engendré
par $m_1-t_1m_0,\dots,m_p-t_pm_0$.
\end{proof}

Dans le cas d'un module \pf le \tho \ref{basic} appliqué
à une matrice de présentation du module nous donne des 
énoncés plus
sophistiqués comme dans~\cite{Hei84}.

Donnons par exemple une formulation avec la dimension de 
Krull, plus facile à formuler.

\begin{theorem}
\label{Forster1} \emph{(Théorème de Forster, version non noethérienne de 
Heitmann)}\\
Soit  $M$  un module \pf sur $\gA$. Notons  $\ff_k$
son \idf d'ordre $k$ et $V_{k}=\VA(\ff_{k})\subseteq\Spec\,\gA$%
.
\begin{itemize}
\item Version \cov (sans points). 
Si $\Kdim(\gA[\nu^{-1}]/\ff_k)<  m-k$ pour $k=0,\ldots ,m$ 
et pour tout générateur $\nu$ de $\ff_{k+1}$, 
alors $M$ peut être engendré
par $m$ \eltsz.
\item Version classique (avec points). 
Pour tout \idep $\fp$ de $\gA$, notons $\mu_\fp$ la dimension de 
l'espace vectoriel $M\otimes_\gA\Frac(\gA/\fp)$ sur le corps 
$\Frac(\gA/\fp)$
(autrement dit $\mu_\fp=k$ \ssi $\fp\in V_{k-1}\setminus V_{k}$). 
\begin{itemize}
\item Si pour  $k=0,\ldots ,m$ et pour tout $\fp\in V_{k}\setminus 
V_{k+1}$ on a
$\Kdim(\gA/\fp)\leq m-k$,\\ ou ce qui revient au même,
\item si pour tout $\fp\in V_{0}$ on a
$\Kdim(\gA/\fp)+\mu_\fp\leq m$,
\end{itemize}
alors $M$ peut être engendré par $m$ \eltsz.
\end{itemize}
\end{theorem}

Rappelons que le support du module $M$ est le fermé $V_{0}$.
La version avec points peut donc être énoncée
en demandant que pour tout $\fp$ dans le support de $M$ on ait
$\Kdim(\gA/\fp)+\mu_\fp\leq m$.

La version sans points reste valable en remplaçant la dimension de
Krull par celle de Heitmann, mais pas la version avec points, car on 
n'obtient pas alors l'équivalent de la version sans points.

La version sans points avec la $\Jdim$ (reformulée avec points) est
la variante du \tho de Swan obtenue par Heitmann dans \cite{Hei84}.

\junk{
\begin{theorem}
\label{Forster2} \emph{(de Swan modifié à la Heitmann, avec la 
$\Hdim$)}
\\
Soit  $M$  un module \pf sur $\gA$ et $\ff_k$
son \idf d'ordre $k$
\begin{itemize}
\item Version \cov (sans points): ;
si $\Hdim(\gA/\ff_k)[1/a]<  m-k$ pour $k$ et pour tout 
$a\in\ff_{k+1}$ ,
alors $M$ peut être engendré par $m$ \eltsz.
\item Version classique (affaiblie, avec points): 
Supposons que pour tout \idep $\fp$ de $\gA$, et pour tout $a\in\gA$
$\Hdim(\gA/\fp[1/a])+r_\fp\leq m$,
alors $M$ peut être engendré par $m$ \eltsz.
\end{itemize}
\end{theorem}

}

\section[\Tho de Swan]{\Tho de Swan, version non noethérienne}
\label{secSwan}

Nous passons maintenant à des versions où la dimension de Heitmann
qui intervient n'utilise plus de localisations.

Ceci nous permet d'améliorer le \tho de Swan lui-même (\cite{Swan67}),
répondant ainsi positivement à une question posée dans 
\cite{Hei84}.

\subsection{Manipulations \elrs de colonnes}

\begin{lemma}
\label{mainlemma}
Soient $L,L_1\dots,L_k\in\gA^{m}$ et $b_1,\dots,b_k\in A$.
Si $\Hdim(\gA) < k$ et 
$$1 = \DA(a,b_1,\dots,b_k)\vu \DA(L)$$
alors il existe $x_1,\dots,x_k\in a\gA$ tels que 
$$1 = \DA(b_1+ax_1,\dots,b_k+ax_k)\vu \DA(L+x_1L_1+\cdots+x_kL_k)$$
\end{lemma}

\begin{proof}
La preuve est par induction sur $k$.  Pour $k=0$, c'est clair.  Si
$k>0$, soit $\fj=\rH_\gA(b_k)$ l'\id bord de Heitmann $b_k$.  On a
$b_k\in \fj$ et $\Hdim(\gA/\fj)< k-1$, donc par induction, on peut
trouver $y_1,\dots,y_{k-1}$ tels que
$$
1 = \rD(b_1+a^2y_1,\dots,b_{k-1}+a^2y_{k-1}) \vu 
(L+ay_1L_1+\cdots+ay_{k-1}L_{k-1})\quad \textsl{dans}\;\gA/\fj\quad 
\quad (\alpha)
$$
Posons $L'=L+ay_1L_1+\cdots+ay_{k-1}L_{k-1}$, 
$X=(b_1+a^2y_1,\dots,b_{k-1}+a^2y_{k-1})$ et remarquons que 
$\DA(L',a)=\DA(L,a)$ et $\DA(X,a)=\DA(a,b_1,\ldots b_{k-1})$. L'\egt 
$(\alpha)$  signifie qu'il existe $y_k$ tel que $b_ky_k\in \JA(0)$ et
$$
1 = \DA(X) \vu \DA(L') \vu \DA(y_k)\vu \DA(b_k)\quad \quad (\beta)
$$
On a
$$
\DA(X)\vu \DA(L'+ay_kL_k) \vu \DA(a^2)\vu \DA(b_k) = 
\DA(a,b_1,\dots,b_k,L) = 1\quad \quad (\gamma) 
$$
et, d'après $(\beta)$,
$$
\DA(X)\vu \DA(L'+ay_kL_k)\vu \DA(y_k)\vu \DA(b_k) = \DA(X) \vu 
\DA(L') \vu \DA(y_k) \vu \DA(b_k) = 1\quad \quad (\delta)
$$
$(\gamma)$ et $(\delta)$ impliquent
$$ \DA(X,L'+ay_kL_k ,b_k,a^2y_k) = 1 = 
\JA(X,L'+ay_kL_k,b_k,a^2y_k)\quad \quad (\eta)$$
et d'après le lemme \ref{gcd} puique $b_ka^2y_k\in \JA(0)$
$$1 = \JA(X,L'+ay_kL_k,b_k+a^2y_k)$$
\cad
$$
1\in\gen{b_1+a^2y_1,\dots,b_{k-1}+a^2y_{k-1},b_k+a^2y_k,L+ay_1L_1+
\cdots+ay_{k-1}L_{k-1}+ay_kL_k}.
$$
\end{proof}

\rems \\
1) Lorsqu'on remplace la dimension de Heitmann par celle de Krull,
le lemme précédent admet une version \gui{uniforme} dans laquelle 
les $x_i$ ne dépendent pas de $a$ (comme dans le 
lemme~\ref{lemKroH}).\\
2) Le lemme ci-dessus peut être vu comme une variante raffinée
du \tho \ref{Bass}.  En fait dans la suite, nous utilisons ce lemme
uniquement avec l'hypothèse renforcée comme suit: Si $\Hdim(\gA) <
k$ et $1 = \DA(a,L)$ \ldots

\medskip On en déduit le corolaire suivant (qui simplifie
l'argument dans \cite{CLQ2004} et a été suggéré par Lionel Ducos
dans \cite{Duc2006}).

\begin{corollary}\label{MAIN}
Supposons que $C,C_1,\dots,C_k$ sont des vecteurs dans $\gA^m$ et
$\nu$ est un mineur d'ordre $k$ de $C_1,\dots,C_k$.  Si $1 =
\DA(\nu)\vee \DA(C)$ et $\Hdim(\gA) < k$ alors il existe
$x_1,\dots,x_k\in A$ tels que $1 = \DA(C+x_1C_1+\dots+x_kC_k)$.
\end{corollary}

\begin{proof}
On prend $a=\nu$ et $b_i$ le mineur obtenu en remplaçant $C_i$ par
$C$ dans $C_1,\dots,C_k$.  On applique le lemme \ref{mainlemma}.  On
remarque alors que $b_i + ax_i$ est aussi le mineur obtenu en
remplaçant $C_i$ par $C + x_1C_1 + \cdots + x_nC_n$ dans $C_1,
\ldots , C_n$ de sorte que $b_i + ax_i\in\DA(C + x_1C_1 + \cdots +
x_nC_n) $ et
$$
1 = \DA(b_1 + ax_1, \ldots , b_n + ax_n) \vu \DA(C + x_1C_1 + 
\cdots + x_nC_n)
$$
implique que $\DA(C + x_1C_1 + \cdots + x_nC_n)=1$.
\end{proof}

Nous allons prouver une variante du \tho \ref{matrix}.
Nous reprenons les notations~\ref{notaMatrix}.

\begin{lemma}\label{MAINCOR}
Si $\DA(C)\vee \Delta_k(G) = 1$ et $\Hdim(\gA)<k$ alors il existe
$t_1,\dots,t_p$ tels que $C+t_1C_1+\dots+t_pC_p$ est unimodulaire.
\end{lemma}

\begin{proof}
Comme $\DA(C)\vee \Delta_k(G) = 1$, on a une famille
$\nu_1,\dots,\nu_p$ de mineurs de $G$ d'ordre $k$ tels que $1 =
\DA(C)\vee \bigvee _i \DA(\nu_i)$.  On applique alors le corolaire
\ref{MAIN} pour avoir $\DA(C)\vee \DA(\nu _k) = 1$ dans $\gA/J_k$ avec
$J_k = \bigvee _{i>k} \DA(\nu_i)$ jusqu'à ce que nous obtenions
$\DA(C) = 1$ dans $\gA$.
\end{proof}

\begin{theorem}
\label{matrix2}
Fixons $k\leq p$. Supposons que la dimension de Heitmann de $\gA$ est
$<k$ et $\Delta_k(F) = 1$. Alors il existe $t_1,\dots,t_p$ tels que
$1 = \DA(C_0+t_1C_1+\cdots+t_pC_p)$.
\end{theorem}

\begin{proof}
Cela résulte directement du lemme \ref{MAINCOR}.
\end{proof}


\subsubsection*{Les \thos  de Serre et Swan avec la dimension de 
Heitmann}
\addcontentsline{toc}{subsubsection}{Les \thos  de Serre et Swan 
avec la dimension de 
Heitmann}

On a alors les résultats suivants, que l'on prouve à partir du
\tho \ref{matrix2} par les mêmes arguments que pour les \thos 
\ref{thSerre}
et \ref{forster}.

\begin{theorem}
\label{thSerre2} \emph{(Splitting Off de Serre, version dimension de Heitmann)}\\
Soit $M$ un \Amo projectif de rang $\geq k$ sur un anneau $\gA$ tel
que $\Hdim \gA< k$.  Alors $M\simeq N\oplus \gA$ pour un certain $N$.
\end{theorem}

\begin{theorem}
\label{thSwan} \emph{(de Swan, version dimension de Heitmann)}\\
Si $\Hdim(\gA)\leq d$ et si $M$ est un module \tf sur $\gA$ localement engendré par $r$ \eltsz, alors $M$ peut être engendré par $d+r$
\eltsz.
\end{theorem}

\subsection{D'autres résultats}

\subsubsection*{Le \tho  de Swan, forme sophistiquée}
\addcontentsline{toc}{subsubsection}{Le \tho  de Swan, forme 
sophistiquée}

 Un raffinement du théorème \ref{matrix2} est le résultat 
suivant.
\begin{theorem}
\label{basic2}
Supposons que $1 \in \Delta_1(F)$ et que  pour  $k=1,\ldots ,p,$  
l'anneau
$\gA/\Delta_{k+1}(F)$ est de dimension de Heitmann $< k$. Alors il 
existe
$t_1,\dots,t_p$ tels que le vecteur $C_0+t_1C_1+\cdots+t_pC_p$ est 
unimodulaire.
\end{theorem}

\begin{proof}
On applique le lemme \ref{MAINCOR} successivement avec  
$\gA/\Delta_2(F)$, $\gA/\Delta_3(F)$, \ldots (raisonner comme dans 
la preuve du \tho~\ref{basic}).
\end{proof}

\begin{theorem}
\label{thSwan2} \emph{(de Swan général, version dimension de
Heitmann)}\\
Soit  $M$  un module \tf sur $\gA$ et $\ff_k$ son \idf d'ordre $k$.
Si $\Hdim(\gA/\ff_k)<  m-k$ pour $k=0,\ldots ,m$, alors $M$ peut 
être engendré
par $m$ \eltsz.
\end{theorem}
\begin{proof}
Si le module est \pf il suffit d'appliquer le \tho \ref{basic2} à
la matrice transposée d'une matrice de présentation de $M$ 
(de la même manière que pour le \thoz~\ref{Forster1}).\\
Le raisonnement dans le cas général est un peu plus subtil. 
Notons $h_1,\ldots ,h_q$ des \gtrs de $M$ et supposons $q>m$.  Notons
$h=(h_1,\ldots ,h_q)$.  Toute relation entre les $h_i$ peut être
écrite comme un vecteur colonne $C\in\gA^q$ qui vérifie $hC=0$. 
L'idéal de Fitting $\ff_{q-i}$ de $M$ est l'idéal $\Delta_i$
engendré par les mineurs d'ordre $i$ des matrices $F\in \gA^{q\times
n}$ qui vérifient $hF=0$ (\cad les matrices dont les colonnes sont
des relations entre les $h_i$).  Nous les appelerons des \gui{matrices
de relation}.\\
Comme $\DA(\Delta_q)=\DA(\ff_0)=\DA(\Ann(M))$ et que nous pouvons
remplacer $\gA$ par $\gA/\Ann(M)$, nous supposons \spdg que
$\DA(\Delta_q)=\DA(0)$.\\
Les hypothèses impliquent que $\Hdim(\gA/\Delta_k)<k-1$ pour
$k=1,\ldots,q$ (ce sont les hypothèses lorsque $q=m+1$).\\
On a donc $\Delta_1=1$, ce qui se constate sur une matrice de
relations $F$.  On peut alors appliquer le lemme \ref{MAINCOR} avec la
matrice $F$, l'entier $k=1$ et l'anneau $\gA/\Delta_2$.  Ceci nous
donne une relation $C$ telle que $\rD_{\gA/\Delta_2}(C)=1$, \cad 
encore
telle qu'il existe $a\in\Delta_2 $ vérifiant $1+a\in\DA(C)$.\\
Le fait que $a\in\Delta_2$ se constate sur une matrice de relations
$G$ et on considère la nouvelle matrice $F=(C|G)$.  On a donc
$\DA(C)\vu\Delta_2(G)=1$ et $\Hdim(\gA/\Delta_3)<2$.  On peut alors
appliquer le lemme \ref{MAINCOR} avec la nouvelle matrice $F$, 
l'entier
$k=2$ et l'anneau $\gA/\Delta_3$.  Ceci nous donne une relation $C$
telle que $\rD_{\gA/\Delta_3}(C)=1$, et ainsi de suite.\\
On obtient en fin de compte une relation $C$ telle que $\DA(C)=1$, 
puisque $\DA(\Delta_q)=0$.
Et $\Hdim\,\gA=\Hdim(\gA/\Delta_q)<q-1$.  En conséquence le 
corolaire
\ref{corBass} du \tho stable range de Bass s'applique et nous pouvons
transformer la colonne $C$ en $(1,0,\ldots ,0)$ par des manipulations
\elrsz, ce qui revient à
remplacer le système \gtr $h$ par une autre système de $q$ \elts 
dont le
premier est nul, bref à réduire de 1 le nombre de \gtrsz.   
\end{proof}

Notons qu'il n'y a pas à priori de version \gui{avec points} du \tho 
précédent car la $\Hdim$ n'est pas (à notre connaissance) la 
dimension de Krull d'un espace spectral. Pour obtenir une version 
avec points, nous devrions considérer la variante (à priori plus 
faible) dans laquelle la $\Hdim$ est replacée par la $\Jdim$. On 
aurait alors sans doute une formulation classique utilisant les 
points du J-spectre de Heitmann.

\subsubsection*{Le \tho  de simplification de Bass, avec la dimension 
de Heitmann, sans hypothèse nothérienne}
\addcontentsline{toc}{subsubsection}{Le \tho  de simplification de 
Bass}

Nous prouvons maintenant le \tho de simplification de Bass.
\begin{lemma}\label{forbass}
Si $\Delta_k(G) = 1 =\DA(C)\vee \DA(a)$ et $\Hdim(\gA)<k$ alors il 
existe
$t_1,\dots,t_p$ tels que $1 = \DA(C+at_1C_1+\dots+at_pC_p).$
\end{lemma}

\begin{proof}
Il suffit d'appliquer le lemme \ref{MAINCOR} aux vecteurs 
$C,aC_1,\dots,aC_n$
(ce qui remplace $G$ par~$aG$).
\end{proof}

\begin{theorem}
\label{thBassCancel2} \emph{(Théorème de simplification de Bass, version 
dimension
de Heitmann)}\\
Soit $M$ un \Amo projectif de rang $\geq k$ sur un anneau $\gA$ tel 
que
$\Hdim \gA< k$, et $N,\,Q$ deux autres \Amos projectifs de type fini. 
Alors
$M\oplus Q\simeq N\oplus Q$ implique  $M\simeq N$.
\end{theorem}

\begin{proof}
Pour être complet, on redonne l'argument de \cite{EG73b}.  On peut supposer
que $Q=\gA$ et que $N\subset \gA^m$ est image d'une matrice
idempotente.  On a un \iso $\phi:M\oplus \gA\to N\oplus \gA\subset
\gA^{m+1}$.  Soit $(C,a)\in\gA^{m+1}$ le vecteur $\phi(0,1)$.  Il nous
faut construire un \auto de $N\oplus \gA$ qui envoie $(C,a)$ sur
$(0,1)$.  \\
Soit $\alpha :M\oplus \gA\to\gA,\,(U,x)\mapsto x$.  La forme
linéaire $\alpha \circ\phi^{-1}$ envoie $(C,a)$ sur $1$, d'où
l'égalité $\DA(C)\vee \DA(a)=1$.  Donc, par le lemme \ref{forbass}, il existe
$C'\in N$ tel que $\DA(C+aC') = 1$.  On a donc une forme linéaire
$\lambda :N\to\gA$ vérifiant $\lambda(C+aC')=1$ (obtenue comme
restriction d'une forme linéaire sur $\gA^m$).  On a alors les trois
\autos suivants de $N\oplus \gA$ qui par composition réalisent le but fixé:
$$
\begin{array}{rcclcrcl} 
\psi_{1}: &(V,x)   &\mapsto   &(V+xC',\lambda (xC-aV))   & , & 
(C,a)&\mapsto &(C+aC',0) \\ 
\psi_{2}: &(V,x)   &\mapsto   &(V,x+\lambda(V))   & , & 
(C+aC',0)&\mapsto 
&(C+aC',1) \\
\psi_{3}: &(V,x)   &\mapsto   &(V-x(C+aC'),x)   & , & 
(C+aC',1)&\mapsto 
&(0,1)  
\end{array}
$$
(on vérifie que $\psi_{1}$ est un \auto en remarquant que si 
$\psi_{1}(V,x)=(W,y)$ alors $x=y+a\lambda(W)$).
\end{proof}

\rems Dans le fonctionnement de cette preuve, on voit qu'on aurait pu se contenter
de supposer que $N$ est un module image d'une matrice $F$ telle que
$\Delta_{k}(F)=1$ et $M\oplus Q\simeq N\oplus Q$ pour un module
projectif $Q$.  \\
Notons aussi que les \autos $\psi_{i}$ sont réalisés en pratique
comme des \autos de~$\gA^{n+1}$ qui fixent $N\oplus \gA$.  De
manière plus générale tous les \thos d'\alg commutative que nous
avons démontrés se ramènent en fin de compte à des \thos
concernant les matrices et leurs manipulations \elrsz.

\medskip 
\noindent 
{\bf Remerciements.} Juillet 2008. Nous remercions Peter Schuster qui a fait un relecture
très attentive de l'article et nous a signalé des erreurs, notamment dans les lemmes 5.6 \hbox{et 5.7}, corrigées dans cette version.

\newpage
\small
\bibliographystyle{apalike-fr}
\bibliography{Heitmannbib}


\normalsize
\newpage
\rdb
\subsection*{Post-Scriptum: Errata dans l'article original.}
\addcontentsline{toc}{section}{Post-Scriptum}

\medskip $\bullet$ Le fait \ref{factRecolTD} et la proposition 
\ref{propRecolTD} ne sont pas corrects dans la formulation générale qui était proposée. Un commentaire à ce sujet est fait après la proposition \ref{propRecolSpec}. On donne maintenant des énoncés corrects (moins forts) avec leurs démonstrations. Un lemme est ajouté, ce qui décale les numéros qui suivent d'une unité.

\medskip $\bullet$ 
Dans la section \ref{subsecAgHagB} \textsl{Algèbres de Heyting, de Brouwer, de Boole,}
il faut lire $a\leq \neg\neg a$ et non pas l'inégalité contraire.

\medskip $\bullet$ 
Dans la section \ref{secESSP} on a remplacé $\rD_\gT(\cdot)$ et $\rV_\gT(\cdot)$ par 
$\DT(\cdot)$ et $\VT(\cdot)$ pour les ouverts et fermés de $\Spec\,\gT$. Concernant $\Spec\,\gT\cir$ on a introduit les notations $\DTo(\cdot)$ et $\VTo(\cdot)$ pour que le propos soit plus clair (proposition \ref{FdSES}).

\medskip $\bullet$ 
Dans le  point \textsl{4} de la proposition \ref{propositionFSES}, on a rectifié comme suit:

\smallskip \noindent  pour tous
ouverts quasi-compacts $U_1$ et $U_2$, l'adhérence de $U_1\setminus U_2$ est le complémentaire d'un 
ouvert quasi-compact.

\medskip $\bullet$ 
La proposition 3.13 de l'article original a été ramenée en \ref{propJspecjspec}, où elle a mieux sa place.
Cette proposition et sa démonstration sont clarifiées. Dans le point \textsl{2} l'hypothèse \gui{$\jspec\,\gT$ noethérien} a été remplacée par \gui{$\Max\,\gT$ noethérien}. 

\medskip $\bullet$ 
Juste avant la proposition \ref{propZar}
l'équivalence 
$$ \widetilde{U}\leq_{\mathrm{Zar}\mathbf{A}} \widetilde{J}
\quad\Longleftrightarrow \quad
\prod\nolimits_{u\in U} u  \in \sqrt{\langle{J}\rangle}
\quad\Longleftrightarrow \quad
{\cal M}(U)\cap \langle{J}\rangle \neq \emptyset
$$

\noindent doit être remplacée par la suivante, dans laquelle
le premier terme est changé

$$ \bigwedge\widetilde{U}\leq_{\mathrm{Zar}\mathbf{A}}\bigvee \widetilde{J}
\quad\Longleftrightarrow \quad
\prod\nolimits_{u\in U} u  \in \sqrt{\langle{J}\rangle}
\quad\Longleftrightarrow \quad
{\cal M}(U)\cap \langle{J}\rangle \neq \emptyset
$$

\medskip $\bullet$ 
Dans la section \ref{subsecBKH}, les lemmes \ref{HenriLemma} et \ref{thCor2.2Heit} sont modifiés, une erreur dans la preuve a été rectifiée.

\medskip Par ailleurs, signalons que nous avons en plusieurs endroits ajouté des références nouvelles et des commentaires variés.

\end{document}